\documentclass[11pt]{amsart}

\usepackage{amsmath, amsthm, amssymb, fullpage}

\usepackage{paralist}
\usepackage{enumerate}

\usepackage{hyperref} 

\newcommand{\rd}{\mathrm{d}}
\newcommand{\rM}{\mathrm{M}}
\newcommand{\cO}{\mathcal{O}}
\newcommand{\bA}{\mathbb{A}}
\newcommand{\bP}{\mathbb{P}}  
\newcommand{\bR}{\mathbb{R}}
\newcommand{\bC}{\mathbb{C}} 
\newcommand{\bT}{\mathbb{T}}  
\newcommand{\bK}{\mathbb{K}}
\newcommand{\bZ}{\mathbb{Z}}
\newcommand{\sP}{\mathsf{P}}
\newcommand{\sA}{\mathsf{A}}
\newcommand{\sH}{\mathsf{H}}
\newcommand{\sR}{\mathsf{R}}
\newcommand{\sF}{\mathsf{F}}
\newcommand{\sJ}{\mathsf{J}} 
\newcommand{\as}{\mathfrak{as}}
\newcommand{\ph}{\mathfrak{ph}}
\DeclareMathOperator{\spe}{sp}
\DeclareMathOperator{\ord}{ord} 
\DeclareMathOperator{\PGL}{PGL}
\DeclareMathOperator{\Aff}{Aff}
\DeclareMathOperator{\dep}{depth} 
\DeclareMathOperator{\hypRes}{hypRes}
\DeclareMathOperator{\ordRes}{ordRes}
\DeclareMathOperator{\Fix}{Fix} 
\DeclareMathOperator{\Id}{id}
\DeclareMathOperator{\MHRL}{MHRL}
\DeclareMathOperator{\Char}{char}
\DeclareMathOperator{\supp}{Supp}
\DeclareMathOperator{\pr}{pr}
\DeclareMathOperator{\Res}{Res}
\DeclareMathOperator{\Spec}{Spec}
\DeclareMathOperator{\Sch}{Sch}
\DeclareMathOperator{\im}{Im}
\DeclareMathOperator{\Rat}{Rat}
\DeclareMathOperator{\SL}{SL}
\DeclareMathOperator{\initial}{in}
\DeclareMathOperator{\gr}{gr}
\DeclareMathOperator{\tr}{tr}
\DeclareMathOperator{\Sym}{Sym}
\DeclareMathOperator{\GL}{GL}
\newcommand{\Vu}{\boldsymbol{u}}
\newcommand{\Vv}{\boldsymbol{v}}
\newcommand{\Vw}{\boldsymbol{w}}
\newcommand{\widevec}[1]{\overrightarrow{#1}}

\theoremstyle{plain}
\newtheorem{theorem}{Theorem}[section]

\newtheorem{proposition}[theorem]{Proposition}

\newtheorem{mainth}{Theorem}

\theoremstyle{definition}
\newtheorem{definition}[theorem]{Definition}
\newtheorem*{acknowledgement}{Acknowledgement}
\newtheorem{remark}[theorem]{Remark}
\theoremstyle{remark}
\newtheorem{fact}[theorem]{Fact}

\theoremstyle{plain}
\newtheorem{apptheorem}{Theorem}[section]

\newtheorem{appproposition}[apptheorem]{Proposition}
\newtheorem{appcorollary}[apptheorem]{Corollary}
\theoremstyle{remark}

\theoremstyle{definition}
\newtheorem{appdefinition}[apptheorem]{Definition}

\numberwithin{equation}{section}

\makeatletter
\@namedef{subjclassname@2020}{\textup{2020} Mathematics Subject Classification}
\makeatother

\begin{document} 

\title[]{Intrinsic semistable reduction loci for the iterations of non-archimedean quadratic rational functions and moduli-theoretic applications}

\author{Y\^usuke Okuyama}
\address{Division of Mathematics, Kyoto Institute of Technology, Sakyo-ku, Kyoto 606-8585 JAPAN}
\email{okuyama@kit.ac.jp}



\date{\today} 
 
\begin{abstract}
We study the intrinsic semistable reduction loci for the iterations of quadratic rational functions over algebraically closed complete non-archimedean fields. Using a reduction-theoretic slope formula for the hyperbolic resultant function and the non-archimedean argument principle, we establish a precise eventual stationarity theorem when the residue characteristic is different from $2$.
In the finite-order case, we define an intrinsic prehole parameter, and for quadratic rational functions with the same reduction order and boundary period, classify the intrinsic reduction data of all iterations at distinguished Fatou boundary points up to simultaneous projective equivalence.
For antispecialization algebraic arcs degenerating to an order $p>1$ boundary point of the compactified quadratic moduli space over an algebraically closed field of characteristic $0$, we define an intrinsic antispecialization parameter using the fixed-point multiplier polynomial. We show that this parameter equals the prehole parameter of every Puiseux lift of an antispecialization arc degenerating to the boundary point and that the $p$th iteration specialization determines all iteration specializations. Over $\bC$, the two intrinsic parameters coincide with DeMarco's parameter $\tau^2$ after an appropriate interpretation.
As a global consequence, we construct explicit blow-up models for the
iterate compactification over arbitrary algebraically closed fields of
characteristic $0$ and describe its boundary intrinsically.
Over $\bC$, this also recovers DeMarco's compactification.
\end{abstract}

\maketitle

\section{Introduction}\label{sec:intro}

Let $K$ be an algebraically closed field that is complete with respect to a non-trivial and non-archimedean absolute value $|\cdot|$. Put $k:=K^\circ/K^{\circ\circ}$, the residue field of $K$, which is still algebraically closed. Here and below, for a valued field $F$, we denote by $F^\circ$ and $F^{\circ\circ}$ the ring of $F$-integers and the unique maximal ideal of $F^{\circ}$. The projective line $\bP^1$ is densely embedded into its Berkovich
analytification $\sP^1=\sP^1_K$, the topology of which is the weak one and which is obtained by gluing two copies of the Berkovich affine line $\sA^1=\sA^1_K$ along the analytification of $\mathbb{G}_m$ via the coordinate change $z\mapsto z^{-1}$ (under the convention that $\infty=1/0=[0:1]\in\bP^1$). The $\sA^1$ consists of multiplicative seminorms on $K[z]$ extending $|\cdot|$; each $a\in K$ gives the evaluation seminorm $|h|_a=|h(a)|$. Similarly, the set
\begin{gather*}
\PGL(2,K^\circ)\backslash\PGL(2,K)\cong \Aff(2,K)\cdot K^\circ
=\{\text{closed disks in }K\text{ of diameter in }|K^\times|\}
\subset 2^{\bP^1}
\end{gather*}
is embedded into the Berkovich upper half space $\sH^1:=\sP^1\setminus\bP^1$ in $\sP^1$
as the dense subset $\sH^1_{\mathrm{II}}$ of type II points; each closed disk
$B=\ell(K^\circ)=\{z\in K:|z-a|\le r\}$, where $a\in K$,
$r\in|K^\times|$, and $\ell\in\Aff(2,K)$, gives the supremum norm
$|h|_B=\sup_{b\in B}|h|_b$, $h\in K[z]$, which is denoted by
$\xi_{a,r}$. For example, the Gauss (or canonical) point $\xi_g$ in $\sP^1$ is
the point $\xi_{0,1}\in\sH^1_{\mathrm{II}}$. Each closed disk of radius $r\in\bR_{>0}\setminus|K^\times|$ gives a type III point in $\sH^1$ in the same manner, and the totality of them is denoted by 
$\sH^1_{\mathrm{III}}$.
The Berkovich hyperbolic space $(\sH^1,\rho)$ is an $\bR$-tree, where
$\rho$ is the hyperbolic metric on $\sH^1$ so, e.g., that
$\rho(\xi_{a,s},\xi_{a,r})=|\log(s/r)|$, $a\in K$ and
$r,s\in\bR_{>0}$,
and the Gromov boundary of $\sH^1$ is canonically identified with $\bP^1$.
For simplicity, $\sH^1$ also denotes $(\sH^1,\rho)$.
Each direction $\Vv=\widevec{\xi\xi'}$ in the tangent (or direction) space
$T_\xi\sP^1$ of $\sP^1$ at a point $\xi\in\sP^1$ may be identified with
the connected component $U(\Vv)=U(\widevec{\xi\xi'})$ of $\sP^1\setminus\{\xi\}$ 
containing $\xi'\in\sP^1\setminus\{\xi\}$, or equivalently with a germ of
the half open interval $(\xi,\xi']$ between $\xi$ and $\xi'$.
For more details on $\sP^1$,
we refer to \cite{BR10,BenedettoBook,FR09,Jonsson15}.

\subsection{Intrinsic reductions of a rational function}
Let $\phi\in K(z)$ be a rational function on $\bP^1$ of degree $d\ge 1$. The action of $\phi$ on $\bP^1$ extends uniquely to an analytic action on $\sP^1$, so that both $\bP^1$ and $\sH^1$ are invariant under $\phi$. This analytic action of $\phi$ on $\sP^1$
remains continuous, open, discrete, surjective, 
and $d$ to $1$ taking into account the local degree $\deg_\xi(\phi)\in\{1,2,\ldots,d\}$ of $\phi$ at each $\xi\in\sP^1$, and along each sufficiently short segment issuing from a point of $\sH^1$ in a direction $\Vv$, it scales the $\rho$-length by the directional local degree $m_{\Vv}(\phi)\ge1$ (see Fact \ref{th:argument} for more details). For each $\xi\in\sP^1$,
the tangent map $\phi_{*,\xi}:T_\xi\sP^1\to T_{\phi(\xi)}\sP^1$ is
induced by $\phi$ at $\xi$, mapping
a direction $\Vv=\widevec{\xi\xi'}$ to the direction $\widevec{\phi(\xi)\phi(\xi')}$
so long as the point $\xi'\in U(\Vv)$ is close enough to $\xi$, and is surjective. When $\xi$ is of type II, this tangent map has degree $\deg_\xi(\phi)$ as a map between the natural projective $k$-lines (see Fact \ref{th:algred} and \eqref{eq:directional} below for more details).

In \cite[Definition 1.1]{okudepth}, both the {\em intrinsic reduction}
\begin{gather*}
 \tilde{\phi}_\xi\begin{cases}
		    :=\phi_{*,\xi} & \text{if }\phi(\xi)=\xi,\\
		    :\equiv\widevec{\xi(\phi(\xi))} & \text{otherwise}
		   \end{cases}:T_\xi\sP^1\to T_\xi\sP^1
\end{gather*}
of $\phi$ at each point $\xi\in\sP^1$ and the intrinsic depth
\begin{gather*}
 \dep_{\Vv}\tilde{\phi}_\xi:=(\phi^*\delta_{\xi})\bigl(U(\Vv)\bigr)\in\{0,\ldots,\deg\phi\}
\end{gather*}
of this intrinsic reduction $\tilde{\phi}_\xi$ at each direction $\Vv\in T_\xi\sP^1$
are introduced, where $\delta_\xi$ is the Dirac measure on $\sP^1$
at $\xi$ and where we set the pullback
\begin{gather*}
 \phi^*\delta_\xi=\sum_{\xi'\in\phi^{-1}(\xi)}(\deg_{\xi'}(\phi))\delta_{\xi'}\quad\text{on }\sP^1
\end{gather*}
of $\delta_\xi$ under $\phi$. In \cite{Okugeometric}, 
before those notions mentioned above, when $d>1$,
the hyperbolic resultant function
\begin{gather*}
 \hypRes_\phi(\xi):=
\frac{\rho(\xi,\xi_g)}{2}
+\frac{\rho(\xi,\phi(\xi)\wedge_{\xi_g}\xi)-\int_{\sP^1}\rho(\xi_g,\xi\wedge_{\xi_g}\cdot)(\phi^*\delta_{\xi_g})(\cdot)}{d-1}\in\bR
\end{gather*}
for $\phi$ on $(\sH^1,\rho)$ giving a global and explicit hyperbolic geometric formula for Rumely's function $\ordRes_\phi$ introduced in \cite{Rumely13}
is introduced, where $\xi\wedge_{\xi''}\xi'$ 
denotes the ``triple'' point in $\sP^1$
among $\xi,\xi',\xi''\in\sP^1$, i.e., the unique point in
the intersection of all the (possibly trivial) closed intervals $[\xi,\xi']$, $[\xi,\xi'']$, and $[\xi',\xi'']$ in $\sP^1$.
When $d>1$, for every $\xi\in\sH^1$ and every $\Vv\in T_\xi\sP^1$,
these three intrinsic notions are
tied by the following reduction-theoretic slope formula
\begin{gather}
 \rd_{\Vv}\hypRes_\phi=\frac{1}{d-1}
\biggl(-\dep_{\Vv}\tilde{\phi}_{\xi}
+\frac{1}{2}\cdot\begin{cases}
		  d-1 & \text{if }\tilde{\phi}_\xi(\Vv)=\Vv\\
		  d+1 & \text{otherwise}
		 \end{cases}\biggr)\label{eq:slope}
\end{gather}
\cite[(1.4)]{okudepth}
(where $\rd_{\widevec{\xi\xi'}}
=\frac{\rd}{\rd(\rho|(\xi,\xi'])}\bigl|_\xi$
is the directional derivative operator). Indeed, $\hypRes_\phi$ is convex, (metrically) proper, and piecewise affine, so in particular
it attains its minimum at some point in $\sH^1$.

Our aims in this paper are twofold: to establish the stationarity of the intrinsic semistable reduction loci in $\sP^1$ under iteration for quadratic rational functions, and to develop its moduli-theoretic applications, including intrinsic parameters for reduction data and algebraic arcs, their comparison with DeMarco's parameter over $\bC$, and the intrinsic stratification of $\overline{\rM}_2$.

\subsection{Main result: the stationarity of intrinsic semistable reduction loci for iterated quadratic rational functions}
We say a rational function $\phi\in K(z)$ of degree $d>1$
has an intrinsic {\em semistable (resp.\ stable)}
reduction at a point $\xi\in\sH^1$ if for every direction $\Vv\in T_\xi\sP^1$,
\begin{gather}
\begin{cases}
 \dep_{\Vv}\tilde{\phi}_{\xi}\le\dfrac{d+1}{2}
 \,\bigl(\text{resp.}\,\le\dfrac{d}{2}\bigr)
 \text{ and more strictly}\vspace*{5pt}\\
 \dep_{\Vv}\tilde{\phi}_{\xi}<\dfrac{d}{2}
 \,\bigl(\text{resp.}\,<\dfrac{d-1}{2}\bigr)
 \text{ if also }\tilde{\phi}_\xi(\Vv)=\Vv.
\end{cases}\label{eq:stable}
\end{gather}
This notion of intrinsic (semi)stable reduction on $\sH^1$ extends the corresponding
potential GIT (semi)stable reduction at points of $\sH^1_{\mathrm{II}}$
(see \cite[\S1.1]{okudepth}, which extends \cite[Theorem C]{Rumely17});
for the GIT
(semi)stable reduction, see Fact \ref{th:git}(\textit{\ref{item:git}}) below. By the slope formula \eqref{eq:slope}, the
intrinsic semistable reduction locus for $\phi$ in $\sH^1$
coincides with the minimum locus 
\begin{gather*}
 \MHRL(\phi):=\arg\hypRes_\phi 
\end{gather*}
of $\hypRes_\phi$ (so of $\ordRes_\phi$) in $\sH^1$.
Moreover, $\phi$ has an intrinsic {\em stable}
reduction at a point $\xi\in\sH^1$ if and only if $\hypRes_\phi$
attains its minimum {\em uniquely} at $\xi$. When $2\mid d$, this is always
the case; moreover, $\xi\in\sH^1_{\mathrm{II}}$. For details on the minimum locus of
$\ordRes_\phi$ (so of $\hypRes_\phi$), see \cite[\S4.11]{SilvermanDynamics}.

In \cite[Theorem 4]{NO24}, we established stationarity
of the minimum loci of $\hypRes_{\phi^j}$ for the $j$th iterations $\phi^j$ of $\phi$
($\phi^0(z):=z$ by convention) when
$\phi$ is any non-archimedean polynomial $P\in K[z]$ of degree $>1$; more precisely, $\hypRes_{P^j}$ attains its minimum
uniquely at the same point in $\sH^1$ for every integer $j\ge (\deg P)-1$. 

Our principal result is a similar stationarity result for the iterated
quadratic rational functions, which is as precise as our former result for iterated polynomials and provides a foundation for the moduli-theoretic applications developed in the sequel.

\begin{mainth}\label{th:quad}
Let $\phi\in K(z)$ be a quadratic rational function
on $\bP^1$, and let us denote by $\xi_\phi\in\sH^1$
the unique minimum point of $\hypRes_\phi$.

When the intrinsic reduction $\tilde{\phi}_{\xi_\phi}$ of $\phi$ at $\xi_\phi$ is
either constant or acyclic (in the monoid of all selfmaps of $T_{\xi_\phi}\sP^1$), 
we have $\MHRL(\phi^j)\equiv\{\xi_\phi\}$ for every integer $j\ge 1$.

Alternatively, when $\Char k\ne2$ and $\tilde{\phi}_{\xi_\phi}$ has finite order $p$, so in particular $\xi_\phi$ is not in the Berkovich (indeed maximal) ramification locus $\sR(\phi):=\{\xi\in\sP^1:\deg_\xi(\phi)>1\}$ of $\phi$, we have $p>1$ and, 
letting $\xi_0^\phi\in\sH^1\setminus\{\xi_\phi\}$ be the nearest point retraction of $\xi_\phi$ to $\sR(\phi)$, $U=U^\phi$ the Berkovich Fatou component of $\phi$ containing $\xi_\phi$, and $\xi_1^\phi$ the unique point of $\partial U$ in the half open interval $(\xi_\phi,\xi_0^\phi]$, the point $\xi_1^\phi$ is a periodic point of $\phi$ having the exact period
$q=q_\phi$ divisible by $p$, and we have
\begin{gather}
\MHRL(\phi^j)
 \begin{cases}
  =\{\xi_\phi\} & \text{if }1\le j<p,\\
  \subset U\bigl(\widevec{\xi_\phi\xi_1^\phi}\bigr)\cap
  U\bigl(\widevec{\xi_1^\phi\xi_\phi}\bigr)& \text{if } p\le j<q,\\
  =\{\xi_1^\phi\}& \text{if }j\ge q.
\end{cases}
\label{eq:stationarity}
\end{gather}
\end{mainth}

In equal characteristic greater than $2$, explicit semisimple and unipotent examples (see Subsection \ref{sec:residue} for the classification) show that the transitional range $p\le j<q_\phi$ in Theorem \ref{th:quad} can consist of several successive stationary blocks. On each block, $\MHRL(\phi^j)$ is a singleton independent of $j$, and the corresponding points occur successively along $(\xi_\phi,\xi_1^\phi)$ before the locus becomes stationary at $\xi_1^\phi$ for $j\ge q_\phi$. Analogous staircase behavior occurs in the unipotent order $2$ case in equal characteristic $2$, where the ramification locus has wild branches. We omit the detailed computations here.

We now turn to the moduli-theoretic applications of Theorem \ref{th:quad}. 
The rigidity and moduli-theoretic consequences of the constant-or-acyclic case of Theorem \ref{th:quad} are collected in Appendix A.
For the first moduli-theoretic application below, we work in the order $p$ case of Theorem \ref{th:quad}, we retain the notations 
$\xi_\phi$, $U^\phi$, $\xi_0=\xi_0^\phi$, $\xi_1=\xi_1^\phi$, and $q=q_\phi$ for $\phi$, and recall that the tangent space of $\sP^1$ at each point of $\sH^1_{\mathrm{II}}$ has a natural $k$-projective-line structure; see Fact \ref{th:algred} for more details. One of the key observations is
that $\widevec{\xi_1\xi_\phi}$ is a multiple fixed point of $\widetilde{(\phi^q)}_{\xi_1}$ 
(see Proposition \ref{th:multiple}).

\subsection{Intrinsic reduction data and the intrinsic prehole parameter}
\label{sec:residue}

First, we distinguish the two possible types of intrinsic reduction of $\phi$
at $\xi_\phi$ when $\tilde{\phi}_{\xi_\phi}$ has order $p>1$, and introduce the
first of our two new intrinsic parameters. In this order $p$ case,
the intrinsic reduction $\tilde{\phi}_{\xi_\phi}$ may be regarded as an element of
$\PGL(2,k)$, i.e., a projective transformation of $\bP^1_k$,
of order $p$, so has either exactly two fixed directions or a unique fixed
direction at $\xi_\phi$ since $k=\overline{k}$. 
We call those two cases the {\em semisimple} (elliptic) and
{\em unipotent} (parabolic) cases, respectively. In the semisimple case,
$\tilde{\phi}_{\xi_\phi}$ is $\PGL(2,k)$-conjugate to
$\zeta\mapsto\omega\zeta$ for a primitive $p$th root of unity
$\omega\in k^\times$; in the unipotent case, 
$\tilde{\phi}_{\xi_\phi}$ is $\PGL(2,k)$-conjugate
to $\zeta\mapsto\zeta+1$ (and $\Char k=p$).
In either case, the {\em intrinsic prehole parameter}
$\ph(\phi)\in\bP^1(k)$ of $\phi$ is defined as
\begin{gather} 
 \ph(\phi):=
 \begin{cases}
  \infty & \text{if }\xi_1\ne\xi_0,\\
  2+c_\phi+(c_\phi)^{-1}\in k & \text{if }\xi_1=\xi_0;
 \end{cases}
 \label{eq:rdef}
\end{gather}
no matter whether $\xi_1=\xi_0$,
as shown in case (\textit{\ref{item:finite}}) of Section \ref{sec:stationary}, we have
\begin{itemize}
 \item $\tilde{\phi}_{\xi_1}\equiv\overrightarrow{\xi_1\xi_\phi}=:\Vw_\infty=\Vw_\infty^\phi
$, 
 \item $\tilde{\phi}_{\xi_1}$ has the intrinsic depths $>0$
(indeed $=1$) in exactly two directions in $T_{\xi_1}\sP^1$, and
 \item one of such two directions is the $\Vw_\infty$ and the other is denoted by $\Vw_0=\Vw_0^\phi$,
\end{itemize}
and when $\xi_1=\xi_0$, as shown in subcase (\textit{\ref{item:boundary}}) of Section \ref{sec:stationary}, furthermore $\widetilde{(\phi^q)}_{\xi_1}(\Vw_0)=\Vw_\infty$ and
the cross ratio
\begin{gather*}
 c_\phi:=[\Vw_0,\Vw_\infty;\Vw_1,\Vw_2]\in k^\times
\end{gather*}
of the ordered 4 directions, where we write the pullback
$D_\phi:=\bigl(\widetilde{(\phi^q)}_{\xi_1}\bigr)^*(\Vw_0)$
of the divisor $\Vw_0$ on $T_{\xi_1}\sP^1$ under the intrinsic reduction
$\widetilde{(\phi^q)}_{\xi_1}=(\phi^q)_{*,\xi_1}$ as
\begin{gather}
 D_\phi=\Vw_1+\Vw_2,\quad\Vw_1=\Vw_1^\phi,\ \Vw_2=\Vw_2^\phi\in T_{\xi_1}\sP^1\text{ possibly }\Vw_1=\Vw_2.\label{eq:pullbackdiv} 
\end{gather}

Next, we introduce intrinsic reduction data together with the notion of projective equivalence between them and its simultaneous version; these notions play a key role in Theorem \ref{th:residuemoduli} below.
\begin{definition}
For a rational function $\chi\in K(z)$ of degree $\ge 1$,
{\em intrinsic reduction data of $\chi$ at a point $\xi\in\sH^1$} consist of the intrinsic
reduction $\tilde{\chi}_\xi$ together with the intrinsic depth function
$\Vv\mapsto\dep_{\Vv}\tilde{\chi}_\xi$ on $T_{\xi}\sP^1$.
For rational functions $\chi_1,\chi_2\in K(z)$ and points $\xi,\eta\in\sH^1_{\mathrm{II}}$, their intrinsic reduction data at $\xi$ and $\eta$, respectively, are
{\em projectively equivalent via} 
a projective isomorphism $\alpha:T_\xi\sP^1\to T_\eta\sP^1$ if
\begin{gather*}
 \alpha\circ\widetilde{(\chi_1)}_\xi
=\widetilde{(\chi_2)}_\eta\circ\alpha
 \quad\text{and}\quad
 \dep_{\alpha(\cdot)}\widetilde{(\chi_2)}_\eta=\dep_{\,\cdot}\widetilde{(\chi_1)}_\xi
\end{gather*}
on $T_\xi\sP^1$.
Two collections $C_1,C_2$ at $\xi_1,\xi_2\in\sH^1_{\mathrm{II}}$ respectively
of intrinsic reduction data of (various) rational functions in $K(z)$ of degree $>1$ 
labeled by the same (not necessarily finite) set are
{\em simultaneously projectively equivalent} if there is a projective
isomorphism $\alpha:T_{\xi_1}\sP^1\to T_{\xi_2}\sP^1$ such that, for every label, the corresponding data of $C_1,C_2$ are projectively equivalent via $\alpha$.
\end{definition}

The following is our first moduli-theoretic application of Theorem \ref{th:quad}.

\begin{mainth}[Classification by the intrinsic prehole parameter]\label{th:residuemoduli}
Let $\phi,\psi\in K(z)$ be quadratic rational functions in the order $p>1$
case of Theorem \ref{th:quad} such that
$q_\phi=q_\psi=:q$. Then the following are equivalent:
\begin{enumerate}[(i)]
\item\label{item:r} $\ph(\phi)=\ph(\psi)$;
\item\label{item:pdata} for some $j\ge q$, the intrinsic reduction data of $\phi^j$ and $\psi^j$ at
$\xi_1^\phi$ and $\xi_1^\psi$, respectively, are projectively equivalent;
\item\label{item:tail} 
the collections of intrinsic reduction data of $\phi^j$ and
$\psi^j$ at $\xi_1^\phi$ and $\xi_1^\psi$, respectively,
labeled by $j=1,2,\ldots$, are simultaneously projectively equivalent.
\end{enumerate}

Moreover, every projective isomorphism appearing in the (simultaneous) projective equivalence assertion (\textit{\ref{item:pdata}}) or (\textit{\ref{item:tail}}) carries the ordered pair
$(\Vw_0^\phi,\Vw_\infty^\phi)$ to $(\Vw_0^\psi,\Vw_\infty^\psi)$.
\end{mainth}

The construction and classification above apply to arbitrary quadratic rational
functions over $K$ in the order $p>1$ case and are entirely intrinsic: they use
only the Berkovich dynamics, the reduction-theoretic slope formula
\eqref{eq:slope}, and the intrinsic reduction data, without requiring that the
rational function arise from an algebraic family.

\subsection{Intrinsic antispecialization parameter and moduli realization}
\label{sec:moduli}
 
We recall in Fact \ref{th:git} below the standard construction
of the dynamical moduli schemes $\rM_d$ and their GIT compactifications
$\overline{\rM}_d$ over $\bZ$, together with the $j$th iterate maps
$\Phi_{d,j}:\overline{\rM}_d\dashrightarrow\overline{\rM}_{d^j}$
and the corresponding graph closures.
In the quadratic $d=2$ case,
$\rM_2\cong\bA^2_{\bZ}$,
$\overline{\rM}_2\cong\bP^2_{\bZ}$
via the Milnor--Silverman fixed-point multiplier coordinates.
Over a field $k$, a boundary point of $\rM_{d,k}$ is a point of
$\overline{\rM}_{d,k}(k)\setminus\rM_{d,k}(k)$.

We now focus on Puiseux lifts of degenerating algebraic arcs in quadratic
dynamical moduli over the specific valued field below.  The resulting multiplier branches and their
valuation data permit us explicit coefficient and specialization calculations.
In this subsection, for any algebraically closed field $k$ of characteristic $0$, we take
\begin{gather*}
 \bK:=\widehat{\overline{k((t))}},
\end{gather*}
i.e., the completion of the algebraic closure of the Laurent series field
$k((t))$ with respect to the extension of the $t$-adic absolute value.
As usual, we set
\begin{gather*}
 F:=k((t))\quad\text{and}\quad R:=F^\circ=k[[t]],
\end{gather*}
and denote by
\begin{gather*}
 \eta_R:=\Spec F
\end{gather*}
the generic point of $\Spec R$.

The following notion also plays a key role in Theorems \ref{th:arc} and \ref{th:demarco} below.

\begin{definition}[Order $p$ boundary points]\label{th:orderboundary}
A boundary point $x\in\overline{\rM}_{2,k}(k)$ of $\rM_{2,k}$ is called an
{\em order $p>1$ boundary point} if the fixed-point multiplier polynomial
$P_x$ in \eqref{eq:multiplierpoly} satisfies
\begin{gather*}
 \Res_T(P_x,\Phi_p)=0,
\end{gather*}
where $\Res_T$ denotes the resultant with respect to $T$
and $\Phi_p(T)$ is the $p$th cyclotomic polynomial. Equivalently
(see Fact \ref{th:git}(\textit{\ref{item:git}})),
$x=[1:\omega+\omega^{-1}:0]$ in the Milnor--Silverman fixed-point multiplier coordinates for some primitive
$p$th root of unity $\omega\in k^\times$, and $P_x(T)$ is a nonzero scalar multiple of
\begin{gather*}
 (T-\omega)(T-\omega^{-1}).
\end{gather*}
\end{definition}

By Theorem \ref{th:intrinsicstratification}(\textit{ii}) in Appendix
\ref{sec:djrappendix}, a boundary point $x\in\overline{\rM}_{2,k}(k)$ is of order $p$ if and only if
$f^\flat\in\PGL(2,k)$ has order $p$, where $f=H_f f^\flat$ is a stable
representative of $x$ (see Remark \ref{th:algdepth} for the notation).

Fix an order $p>1$ boundary point $x\in\overline{\rM}_{2,k}(k)$ and choose a primitive $p$th root of
unity $\omega\in k^\times$ such that $P_x(\omega)=0$.
The notation $\spe_R$ means taking the value at the unique closed point
of $\Spec R$ of the unique extension of a rational section of a
proper $R$-scheme. 
We call an algebraic arc $\gamma\in\overline{\rM}_{2,k}(R)$ an {\em antispecialization arc degenerating to} $x$ if
\begin{gather*}
 \gamma|_{\eta_R}\in\rM_{2,k}(F)\quad\text{and}\quad \spe_R(\gamma|_{\eta_R})=x.
\end{gather*}
We define the {\em intrinsic antispecialization parameter}
$\as(\gamma)\in\bP^1(k)$ as follows. Letting $\gamma=[s_1:s_2:s_3]$ 
be the Milnor--Silverman fixed-point multiplier coordinates over $R$, so that
$s_1\in R^\times$ and $0\ne s_3\in\mathfrak{m}_R$, we associate to $\gamma$
the fixed-point multiplier polynomial $P_\gamma(T)\in R[T]$ by evaluating the
expression over $\bZ$ in \eqref{eq:multiplierpoly} at $(s_1,s_2,s_3)$, and set
\begin{gather}
 \as(\gamma):=\spe_R\biggl(
 -\frac{p^2\Res_T(P_\gamma,\Phi_p)}
 {(2-\omega-\omega^{-1})\Phi_p'(\omega)\Phi_p'(\omega^{-1})
 s_1^{\deg\Phi_p-1}s_3}\biggr)
 \in\bP^1(k).
 \label{eq:tdef}
\end{gather}
Since $\Res_T(P_\gamma,\Phi_p)$ is homogeneous of degree $\deg\Phi_p$
in the coefficients of $P_\gamma$ and these coefficients are linear in
$(s_1,s_2,s_3)$, the numerator in \eqref{eq:tdef} is homogeneous of degree
$\deg\Phi_p$ in these coordinates. The denominator has the same degree
$(\deg\Phi_p-1)+1=\deg\Phi_p$, and the fraction is symmetric in $\omega$ and
$\omega^{-1}$, so $\as(\gamma)$ is independent of the choices of
$\omega$ and $(s_1,s_2,s_3)$.
In Subsection \ref{sec:proofintrinsic},
we express it in terms of the fixed-point multipliers of $\gamma$.
Finally, for an antispecialization arc $\gamma$ degenerating to $x$,
a {\em Puiseux lift} of $\gamma$ is a quadratic rational function
$\phi\in\bK(z)$ such that, for a finite extension $F'/F$ contained in
$\overline F$, $\phi\in F'(z)$ and $[\phi]\in\rM_{2,k}(F')$ is the base change
of $\gamma|_{\eta_R}$ to $F'$.

The following is our second moduli-theoretic application of Theorem \ref{th:quad}.

\begin{mainth}[Moduli realization of the intrinsic antispecialization parameter]\label{th:arc}
Let $k,\bK,R$ be as above. Let $x\in\overline{\rM}_{2,k}(k)$ be an order $p>1$ boundary point
of $\rM_{2,k}$. Then the following statements hold.
\begin{enumerate}[(i)]
\item\label{item:intrinsic}
For every antispecialization arc $\gamma$ 
degenerating to $x$, every Puiseux lift $\phi\in\bK(z)$ of $\gamma$  
is in the order $p$ case of Theorem \ref{th:quad} and satisfies $q_\phi=p$ and
\begin{gather}
 \ph(\phi)=\as(\gamma)
 \quad\text{in }\bP^1(k).
 \label{eq:rt}
\end{gather}

\item\label{item:limits}
For any antispecialization arcs $\gamma,\gamma'\in\overline{\rM}_{2,k}(R)$
degenerating to $x$ and Puiseux lifts $\phi,\psi$ satisfying
$\xi_1^\phi=\xi_1^\psi=\xi_g$ after $\PGL(2,\bK)$-conjugations
of them if necessary, for every $j\ge p$,
\begin{multline}
 \spe_R\bigl(\Phi_{2,j,k}\circ(\gamma|_{\eta_R})\bigr)
 =\spe_R\bigl(\Phi_{2,j,k}\circ(\gamma'|_{\eta_R})\bigr)\\
\Leftrightarrow
 \widehat{\phi^j}\text{ and }\widehat{\psi^j}\text{ are }\PGL(2,k)\text{-conjugate}
 \Leftrightarrow \as(\gamma)=\as(\gamma').
\label{eq:tclass}
\end{multline} 

\item\label{item:graph}
For every antispecialization arc $\gamma$ degenerating to $x$ and every
Puiseux lift $\phi$ of $\gamma$, after a $\PGL(2,\bK)$-conjugation of
$\phi$ if necessary so that $\xi_1^\phi=\xi_g$, for every $j\ge p$ we have
\begin{gather}
 \widehat{\phi^j}\in\bP(V_{2^j,k})^{\mathrm{s}}(k)
 \quad\text{and}\quad
 \Bigl[\widehat{\phi^j}\Bigr]
 =\spe_R\bigl(\Phi_{2,j,k}\circ(\gamma|_{\eta_R})\bigr)
 \quad\text{in }\overline{\rM}_{2^j,k}(k).
 \label{eq:modlim}
\end{gather}
There is a bijective morphism
\begin{gather}
 \Theta_x:\bP^1_k\to
 \varprojlim_{n\ge p}
 \im_{\Sch}\Bigl(
 \pr_{p,n}:
 \pr_1^{-1}(x)
 \to
 \prod_{j=p}^{n}\overline{\rM}_{2^j,k}
 \Bigr)
 \label{eq:Theta}
\end{gather}
where the inverse limit is taken with respect to the natural truncation
morphisms from the scheme indexed by $m$ to that indexed by $n$, $m\ge n\ge p$, 
and these truncation morphisms are indeed bijective.
For every $\gamma$ and every $\phi$ normalized as above,
\begin{gather*}
 \Theta_x\bigl(\as(\gamma)\bigr)
 =\Bigl(\bigl[\widehat{\phi^j}\bigr]\Bigr)_{j\ge p}
 =\Bigl(\spe_R\bigl(\Phi_{2,j,k}\circ(\gamma|_{\eta_R})\bigr)\Bigr)_{j\ge p}.
\end{gather*}
Here $\pr_1:\mathcal{G}_{2,n,k}\to\overline{\rM}_{2,k}$ and
$\pr_{p,n}:\mathcal{G}_{2,n,k}\to
\prod_{j=p}^n\overline{\rM}_{2^j,k}$
are the projections onto the first factor and onto the factors indexed by
$p,\ldots,n$, respectively.
See Fact \ref{th:git} for details on these notions.

\item\label{item:findet}
Let $\gamma,\gamma'$ be antispecialization arcs degenerating to $x$. If
$\spe_R(\Phi_{2,p,k}\circ(\gamma|_{\eta_R}))
 =\spe_R(\Phi_{2,p,k}\circ(\gamma'|_{\eta_R}))$, then
$\spe_R(\Phi_{2,j,k}\circ(\gamma|_{\eta_R}))
 =\spe_R(\Phi_{2,j,k}\circ(\gamma'|_{\eta_R}))$ for every $j\ge 2$.
Correspondingly, by Theorem \ref{th:residuemoduli}, for every Puiseux lift
$\phi\in\bK(z)$ of $\gamma$, the intrinsic reduction data of $\phi^p$ at $\xi_1$
determine the intrinsic reduction data of $\phi^j$ at $\xi_1$ for every
$j\ge 1$, up to simultaneous projective equivalence.
\end{enumerate}
\end{mainth}

\begin{remark}
By Theorem \ref{th:arc} and the regularity argument in
Theorem \ref{th:acyclicmoduli}, the indeterminacy locus of $\Phi_{2,n,k}$
consists precisely of order $p$ boundary points with $2\le p\le n$;
see Section \ref{sec:proofexcfib} for the proof. Over $\bC$, this recovers
\cite[Theorem~5.1]{DeMarco07}.
By Theorem \ref{th:arc}(\textit{\ref{item:findet}}), the iterate coordinates over
each such boundary point are determined by finitely many of them.
\end{remark}

\subsection{Exceptional fibers and their realization over $k$}
\label{sec:excfib}

We retain the setting of Subsection \ref{sec:moduli}, with $k$ an
algebraically closed field of characteristic $0$.

\begin{definition}[Blow-up models of iterate graphs]\label{th:blowupmodels}
For every $n\ge2$, let $\widetilde{\mathcal{G}}_{2,n,k}\to\overline{\rM}_{2,k}$ be the blow-up of
the ideal sheaf equal to $\cO_{\overline{\rM}_{2,k}}$ away from the
order $p$ boundary points with $2\le p\le n$ and having, at each such boundary point
$x=[1:\omega+\omega^{-1}:0]\in\overline{\rM}_{2,k}(k)$
of $\rM_{2,k}$, the local ideal
\begin{gather*}
 \begin{cases}
 (u,v)&\text{if }p=2,\\
 (u^2,v)&\text{if }p>2,
 \end{cases}
 \quad\text{where}\quad
 u=s_2/s_1-(\omega+\omega^{-1}),\quad v=s_3/s_1.
\end{gather*}
There is a unique finite birational bijective morphism
\[
 \nu_n:\widetilde{\mathcal{G}}_{2,n,k}\to\mathcal{G}_{2,n,k}
\]
extending the iterate-graph morphism on $\rM_{2,k}$.
It is the normalization morphism of $\mathcal{G}_{2,n,k}$, and
$\pr_1\circ \nu_n$ is the natural blow-up morphism. In particular, $\nu_n$ is
a homeomorphism; the existence of $\nu_n$ and its stated properties are proved in
Section \ref{sec:proofexcfib}.
\end{definition}

For every order $p>1$ boundary point $x\in\overline{\rM}_{2,k}(k)$ of $\rM_{2,k}$ and every $n\ge p$, put
\begin{gather*}
 E_{p,n}(x):=(\pr_1\circ \nu_n)^{-1}(x)\subset \widetilde{\mathcal{G}}_{2,n,k},
\end{gather*}
where $\pr_1:\mathcal{G}_{2,n,k}\to\overline{\rM}_{2,k}$ is
the first projection as in Theorem \ref{th:arc}(\textit{\ref{item:graph}}).
We call an algebraic arc $\widetilde{\gamma}\in \widetilde{\mathcal{G}}_{2,n,k}(R)$
an {\em antispecialization arc degenerating to}
$e\in E_{p,n}(x)(k)$ if
\begin{gather*}
 (\pr_1\circ \nu_n\circ\widetilde{\gamma})|_{\eta_R}
 \in\rM_{2,k}(F)
 \quad\text{and}\quad
 \widetilde{\gamma}(0)=e.
\end{gather*}
Then $\gamma:=\pr_1\circ \nu_n\circ\widetilde{\gamma}$
is an antispecialization arc degenerating to $x$,
and $\widetilde{\gamma}$ is a lift of $\gamma$.

Our next moduli-theoretic application is the following, which extends the
exceptional-fiber realization of \cite[Theorem~7.1]{DeMarco07} to $k$
and identifies its parameter intrinsically.

\begin{mainth}[Exceptional-fiber realization]\label{th:demarco}
Let $x\in\overline{\rM}_{2,k}(k)$ be an order $p>1$ boundary point
of $\rM_{2,k}$, and pick any $n\ge p$. Then the intrinsic antispecialization parameter $\as$
defined in \eqref{eq:tdef} induces an isomorphism of projective curves
$\as_{p,n}:E_{p,n}(x)\overset{\sim}{\to}\bP^1_k$
such that
\begin{gather*}
 \pr_{p,n}\circ \nu_n|_{E_{p,n}(x)}
 =\bigl(\pi_n\circ\Theta_x\bigr)\circ\as_{p,n},
\end{gather*}
where $\pi_n$ denotes the natural projection from the inverse limit
in \eqref{eq:Theta} onto the scheme indexed by $n$.
Moreover,
for every $e\in E_{p,n}(x)(k)$, every antispecialization arc
$\gamma\in\overline{\rM}_{2,k}(R)$
degenerating to $x$ which lifts to
an antispecialization arc in $\widetilde{\mathcal{G}}_{2,n,k}$ degenerating to $e$,
and every Puiseux lift
$\phi\in\bK(z)$ of $\gamma$, we have
\begin{gather*}
 \as_{p,n}(e)
 =\as(\gamma)
 =\ph(\phi).
\end{gather*}
\end{mainth}

When $k=\bC$, $\widetilde{\mathcal{G}}_{2,n,k}$ and $\nu_n$
are identified with DeMarco's blow-up space and regular homeomorphism
in \cite[\S7, Theorem~7.1]{DeMarco07}.
The isomorphism $E_{p,n}(x)\cong\bP^1_k$ can also be obtained from
\cite[Theorem~7.1 and its proof]{DeMarco07} by the Lefschetz principle,
using the blow-up models and the corrected parameter formulas there, all defined
over a cyclotomic field (see Remark~\ref{rem:demarco}
for the seemingly required corrections to \cite{DeMarco07}).
We emphasize that our proofs of Theorem \ref{th:demarco} and also of
Theorem \ref{th:compactification} in the next subsection, based on
the Milnor--Silverman fixed-point multiplier polynomial, work
directly on the dynamical moduli space and entirely avoid the
complex-analytic arguments and sophisticated asymptotic computations carried out
by DeMarco on lifts to the space of quadratic rational functions,
which go back at least to Epstein's seminal work
\cite[\S\S3--4, Proposition~2]{Epstein00}.
When $k=\bC$, so that $\bK=\widehat{\overline{\bC((t))}}$ is the
Levi-Civita field, these values also equal $\tau^2(e)$, where $\tau^2$
denotes DeMarco's complex-dynamical parameter (equivalently,
the class $\pm\tau$; see Section \ref{sec:proofexcfib}).

\subsection{Iterate compactification over $k$}
\label{sec:compactification}

We retain the assumption that $k$ is algebraically closed of characteristic
$0$. With the truncation morphisms $\rho_{n,m,k}$ from Fact \ref{th:git},
we regard the inverse system
\begin{gather*}
 \widehat{\rM}_{2,k}:=(\mathcal{G}_{2,n,k})_{n\ge2}
\end{gather*}
as a pro-scheme and call it the {\em iterate compactification} of $\rM_{2,k}$.
Its set of $k$-points is
\begin{gather*}
 \widehat{\rM}_{2,k}(k)
 :=\varprojlim_{n\ge2}\mathcal{G}_{2,n,k}(k)
 \subset\prod_{j\ge1}\overline{\rM}_{2^j,k}(k),
\end{gather*}
and $\pr_1$ denotes projection onto the first coordinate.
The following is a global consequence of Theorems \ref{th:arc} and
\ref{th:demarco} and the regularity result in Appendix \ref{sec:acyclicappendix}.

\begin{mainth}[Iterate compactification over $k$]\label{th:compactification}
The following statements hold.
\begin{enumerate}[(i)]
\item
For every $n\ge m\ge2$, there is a unique morphism
\begin{gather*}
 \widetilde{\rho}_{n,m,k}:
 \widetilde{\mathcal{G}}_{2,n,k}\to\widetilde{\mathcal{G}}_{2,m,k}
\end{gather*}
over $\overline{\rM}_{2,k}$. These morphisms form an inverse system and satisfy
\begin{gather*}
 \nu_m\circ\widetilde{\rho}_{n,m,k}
 =\rho_{n,m,k}\circ\nu_n.
\end{gather*}
Thus the normalization morphisms $\nu_n$ define a morphism of pro-schemes
from this inverse system to $\widehat{\rM}_{2,k}$, inducing a bijection
between their sets of $k$-points.

\item
For every order $p>1$ boundary point $x$ and every $n\ge m\ge p$,
$\widetilde{\rho}_{n,m,k}$ restricts to an isomorphism
$E_{p,n}(x)\to E_{p,m}(x)$ satisfying
\begin{gather*}
 \as_{p,m}\circ\widetilde{\rho}_{n,m,k}|_{E_{p,n}(x)}
 =\as_{p,n}.
\end{gather*}
These compatible parameters identify
$\pr_1^{-1}(x)\subset\widehat{\rM}_{2,k}(k)$ with $\bP^1(k)$.
The fiber over every other point of $\overline{\rM}_{2,k}(k)$ is a singleton.

\item
Every point $(y_j)_{j\ge1}\in\widehat{\rM}_{2,k}(k)$ is uniquely determined
by finitely many coordinates: by $(y_1,y_p)$ if $y_1$ is an order $p>1$
boundary point, and by $y_1$ otherwise.
No projection onto a finite level $\mathcal{G}_{2,n,k}(k)$ is injective,
and no finite sequence of blow-ups of $\overline{\rM}_{2,k}$ resolves all
the rational iterate maps $\Phi_{2,j,k}$, $j\ge2$, simultaneously.
\end{enumerate}
\end{mainth}

We emphasize that our proof derives Theorem \ref{th:compactification}
directly from the intrinsic classification and parameters established
in this paper. This can also be obtained from
\cite[Theorems~1.3--1.4,~5.4, and~7.1 and their proofs]{DeMarco07}
by applying the Lefschetz principle to the finite-level models over
cyclotomic fields and then passing to the inverse systems.
For the latter derivation, the identification of the exceptional-fiber
parameters with $\as_{p,n}$ additionally uses the computations in
\eqref{eq:excp2} and \eqref{eq:excp}.

The rest of the paper is organized as follows. In Section \ref{sec:background},
we recall background on non-archimedean dynamics including the non-archimedean
argument principle and the composition formula for intrinsic depths. Together
with the slope formula \eqref{eq:slope}, these formulas underlie the
stationarity and classification results. In Section \ref{sec:stationary},
we prove Theorem \ref{th:quad}. In Section \ref{sec:residuemoduli}, we prove Theorem
\ref{th:residuemoduli} using the intrinsic pullback divisor $D_\phi$ and computations
of intrinsic depths. In Section \ref{sec:arc}, we recall the GIT construction
of dynamical moduli and the iterate-graph closures, and then prove Theorem
\ref{th:arc}, that is, the coincidence of the two intrinsic parameters,
the moduli realization, and finite determination of iterate specializations.
In Section \ref{sec:proofexcfib}, we construct the blow-up spaces and the morphisms
$\nu_n$ over $k$, prove the exceptional-fiber realization in Theorem
\ref{th:demarco}, and identify the parameter $\as_{p,n}$ of this realization with DeMarco's parameter over
$\bC$. We then prove Theorem \ref{th:compactification}, describing the iterate
compactification through the compatible blow-up models and establishing finite
determination and the impossibility of resolving all iterates by finitely many
blow-ups. In Appendix \ref{sec:acyclicappendix}, we record rigidity results in the
constant-or-acyclic case in Theorem \ref{th:quad} and prove regularity of iteration at the
corresponding boundary points. In Appendix
\ref{sec:djrappendix}, we establish the intrinsic stratification in multiplier
coordinates, which simplifies the statements and proofs of our results,
and recover the quadratic crucial-set stratification of
Doyle--Jacobs--Rumely. In Appendix \ref{sec:normalformsappendix}, we recall the
Type II/III local normal forms from \cite{Jonsson15} used in Section \ref{sec:residuemoduli} and
derive the non-archimedean argument principle mentioned above from these local normal forms, an intrinsic definition of surplus local degrees.

\section{Background}
\label{sec:background}
Let $\phi\in K(z)$ be a rational function on $\bP^1$
of degree $d\ge 1$.

\begin{fact}\label{th:algred}
The reduction $\tilde{\phi}:=\widehat{\phi}^{\flat}\in k(\zeta)$ of $\phi$ modulo $K^{\circ\circ}$ in terms of the notation recalled in Remark \ref{th:algdepth} below 
is identified with $\tilde{\phi}_{\xi_g}$ 
under the canonical identification $\bP^1(k)\ni\hat{a}\leftrightarrow\widevec{\xi_ga}\in T_{\xi_g}\sP^1$, where $\hat{a}$ is the reduction of $a\in\bP^1(K)$ modulo $K^{\circ\circ}$.
In particular, when $\tilde{\phi}_{\xi_g}$ is non-constant, 
it is regarded as a rational function on $\bP^1(k)$ 
of degree $\deg\tilde{\phi}$, which indeed equals $\deg_{\xi_g}(\phi)$.
More generally, 
for every $\xi\in\sH^1_{\mathrm{II}}$, 
each $M\in\PGL(2,K)$ such that $M(\xi)=\xi_g$ 
identifies $T_\xi\sP^1$ with $\bP^1(k)$; for any two such $M_1,M_2$,
$\widetilde{M_1\circ M_2^{-1}}\in\PGL(2,k)$, that is, $T_\xi\sP^1$ has a natural $k$-projective-line structure. For $\xi,\eta\in\sH^1_{\mathrm{II}}$, 
a bijection $\alpha:T_\xi\sP^1\to T_\eta\sP^1$ is called a {\em projective isomorphism} 
if it is an element of $\PGL(2,k)$ with respect to the natural projective $k$-structures
of $T_\xi\sP^1$ and $T_\eta\sP^1$. Every such $\alpha$ is $g_{*,\xi}$ for some
$g\in\PGL(2,K)$ with $g(\xi)=\eta$.
\end{fact}

\begin{fact}\label{th:argument}
We repeatedly use the non-archimedean argument principle
(\cite[Proposition 3.10]{Faber13topologyI}):
 for every $\xi'\in\sP^1$ and every direction  $\Vv\in T_\xi\sP^1$ at every point $\xi\in\sH^1$,
 \begin{gather}
  \bigl((\phi^*\delta_{\xi'})\bigr)(U(\Vv))
 =s_{\Vv}(\phi)+
 \begin{cases}
  m_{\Vv}(\phi) & \text{if }U(\phi_*\Vv)\ni\xi',\\
  0 & \text{otherwise}
 \end{cases}\label{eq:FKN}
 \end{gather}
 in terms of the directional local degree $m_{\Vv}(\phi)\in\{1,2,\ldots,\deg_\xi(\phi)\}$ and
 the surplus local degree $s_{\Vv}(\phi)\in\{0,1,\ldots,d-\deg_\xi(\phi)\}$ of $\phi$ at $\Vv$; more precisely, they respectively satisfy
 \begin{gather}
 \sum_{\Vv\in T_\xi\sP^1:\,\phi_*\Vv=\Vw}m_{\Vv}(\phi)=\deg_\xi(\phi)
\label{eq:directional}
 \end{gather}
for every $\Vw\in T_{\phi(\xi)}\sP^1$
and
 \begin{gather}
 \sum_{\Vv\in T_\xi\sP^1}s_{\Vv}(\phi)=d-\deg_\xi(\phi),\label{eq:surplus}
 \end{gather}
 and moreover, for every $\Vv=\widevec{\xi\xi'}\in T_\xi\sP^1$,
if $\xi'\in U(\Vv)$ is sufficiently close to $\xi$, then the map $\phi:(\xi,\xi')\to(\phi(\xi),\phi(\xi'))$ is affine with respect to the $\rho$-length parameters, with slope equal to the directional local degree $m_{\Vv}(\phi)$.  For such $\xi'$, we also have
\begin{gather}
 m_{\Vv}(\phi)=1\quad\text{if and only if}\quad
\deg_{\,\cdot\,}(\phi)\equiv 1\text{ on }(\xi,\xi').\label{eq:isomet}
\end{gather}
For further details on $m_{\Vv}(\phi)$, $s_{\Vv}(\phi)$, and the
non-archimedean argument principle \eqref{eq:FKN}, see Appendix
\ref{sec:normalformsappendix}.
\end{fact}

We also use the following composition formula for intrinsic depths,
which is a simple application of the non-archimedean argument principle \eqref{eq:FKN}:
for any rational functions $\phi,\psi\in K(z)$ of degree $\ge 1$, 
every $\xi\in\sH^1_{\mathrm{II}}$ such that $\psi(\xi)=\xi$,
and for every $\Vv\in T_\xi\sP^1$, 
\begin{gather}
 \dep_{\Vv}\widetilde{(\phi\circ\psi)}_\xi=(\deg\phi)\dep_{\Vv}\tilde{\psi}_\xi
 +m_{\Vv}(\psi)\dep_{\psi_{*}\Vv}\tilde{\phi}_\xi;
 \label{eq:depthcomp}
\end{gather}
for, we first have
$(\phi\circ\psi)^*\delta_\xi=\psi^*(\phi^*\delta_\xi)
=\sum_{x\in\phi^{-1}(\xi)}\deg_x(\phi)\,\psi^*\delta_x$, and in turn
\begin{align*}
 &\dep_{\Vv}\widetilde{(\phi\circ\psi)}_\xi
 =\sum_{x\in\phi^{-1}(\xi)}\deg_x(\phi)\cdot(\psi^*\delta_x)(U(\Vv))
 =\sum_{x\in\phi^{-1}(\xi)}\deg_x(\phi)
 \bigl(s_{\Vv}(\psi)+m_{\Vv}(\psi)
 1_{U(\psi_{*}\Vv)}(x)\bigr)\\
 =&(\deg\phi)s_{\Vv}(\psi)
 +m_{\Vv}(\psi)(\phi^*\delta_\xi)
 \bigl(U(\psi_{*}\Vv)\bigr)
 =(\deg\phi)\dep_{\Vv}\tilde{\psi}_\xi
 +m_{\Vv}(\psi)\dep_{\psi_{*}\Vv}\tilde{\phi}_\xi
\end{align*}
by the definition of intrinsic depth and \eqref{eq:FKN}. See also \cite[\S2]{DeMarco05}.

\begin{remark}\label{th:algdepth}
When $\xi=\xi_g$, the intrinsic depth $\dep_{\widevec{\xi_ga}}\tilde{\phi}_{\xi_g}$ can also be computed as the depth at $\hat{a}\in\bP^1(k)$ 
of the {\em coefficient} reduction $\hat{\phi}\in\bP^{2d+1}(k)$ of $\phi$ modulo $K^{\circ\circ}$
(see the depth formulas \cite[Lemma 3.17]{Faber13topologyI} and \cite[Proof of Corollary 2.11]{KN23}).
For a coefficient point $f=[F:G]\in\bP^{2d+1}(k)$, where $f(\zeta)=G(1,\zeta)/F(1,\zeta)$ in the affine coordinate $\zeta=Y/X$, putting
\begin{gather*}
 H_f:=\gcd(F,G)\quad\text{and}\quad
 f^\flat:=\frac{G/H_f}{F/H_f},
\end{gather*}
we write $f=H_f f^\flat$. As mentioned in Fact \ref{th:algred},
$\hat{\phi}^\flat=:\tilde{\phi}\leftrightarrow\tilde{\phi}_{\xi_g}$,
and when $\phi(\xi_g)=\xi_g$, 
\begin{gather*}
 \deg_{\hat{a}}\tilde{\phi}=m_{\widevec{\xi_ga}}(\phi)\quad\text{and}\quad
 \ord_{\hat{a}}[H_{\hat{\phi}}=0]=s_{\widevec{\xi_ga}}(\phi).
\end{gather*}
We use this coefficient point notion only for clarifying the exposition.
\end{remark}

\begin{fact}\label{th:FJ}
Suppose in addition that $d>1$. The Berkovich Julia set $\sJ(\phi)$ of $\phi$ is by definition the support of the canonical equilibrium measure $\nu_\phi$
 for $\phi$ on $\sP^1$, the defining property of which is
 the weak convergence
 \begin{gather}
 \nu_\phi=\lim_{j\to\infty}\frac{(\phi^j)^*\delta_\xi}{d^j}\label{eq:canonical}
 \end{gather}
 on $\sP^1$ for any $\xi\in\sH^1$.
 Unless $\sJ(\phi)$ is a singleton, $\sJ(\phi)$
 coincides with the set of points $\xi\in\sP^1$ such that, for every open neighborhood of $\xi$ in $\sP^1$,
 the union of its iterated images under $\phi$ covers
 all of $\sP^1$ except possibly a countable subset of $\bP^1$.
 The Berkovich Fatou set of $\phi$ is the complement in $\sP^1$ of $\sJ(\phi)$, each connected component of which
 is called a Berkovich Fatou component of $\phi$.
 Both the compact subset $\sJ(\phi)$ and the open subset
 $\sF(\phi)$ in $\sP^1$ are totally invariant under $\phi$,
 and we have $\sF(\phi^j)=\sF(\phi)$ (and $\sJ(\phi^j)=\sJ(\phi)$) for any integer $j\ge 1$.

A point $\xi\in\sH^1$ fixed by $\phi$ is said to be repelling if $\deg_{\xi}(\phi)>1$,
and then $\xi\in\sJ(\phi)\cap\sH^1_{\mathrm{II}}$.
A Berkovich Fatou component $U$ of $\phi$ is said to be cyclic
if it is mapped onto itself by $\phi^p$ for some integer $p\ge 1$, and then
$U$ is either an immediate attractive basin, i.e., there is a (super)attracting fixed point $a\in U\cap\bP^1$ of $\phi^p$ and $\lim_{j\to\infty}(\phi^p|U)^j=a$
or a singular domain, i.e., the restriction $\phi^p:U\to U$ is bijective (\cite{Juan03}, see also \cite[Theorem 9.14]{BenedettoBook} in the general setting here).
\end{fact}

\section{Proof of Theorem \ref{th:quad}}\label{sec:stationary}

Let $\phi\in K(z)$ be a quadratic rational function
on $\bP^1=\bP^1(K)$, i.e., $\deg\phi=2$,
and, for notational simplicity, denote by
$\xi:=\xi_\phi\in\sH^1$ the unique minimum point
of $\hypRes_\phi$. Set
\begin{gather*}
 A:=A_{\tilde{\phi}_\xi}
 :=\bigl\{\Vv\in T_\xi\sP^1:\dep_{\Vv}\tilde{\phi}_\xi>0\bigr\}.
\end{gather*}
For an automorphism of a set, an orbit of it is acyclic if and only if this orbit is infinite.

First, if the intrinsic reduction $\tilde{\phi}_\xi$
of $\phi$ at $\xi$ is (non-constant and) $\deg\phi(=2)$ to $1$ or equivalently
$\phi^{-j}(\xi)=\{\xi\}$ for every $j\ge 1$, then
$\dep_{\Vv}\widetilde{(\phi^j)}_\xi=0$
for every $j\ge 1$ and every $\Vv\in T_{\xi}\sP^1$, so
$\MHRL(\phi^j)\equiv\{\xi\}$ for every $j\ge 1$ by the slope formula \eqref{eq:slope}.

Next, if $\tilde{\phi}_\xi$ is constant (i.e., $\equiv\widevec{\xi\phi(\xi)}$) or equivalently $\phi(\xi)\neq\xi$, then by the slope formula \eqref{eq:slope}
and the non-archimedean argument principle \eqref{eq:FKN}, 
we have not only
\begin{gather*}
\dep_{\Vv}\tilde{\phi}_\xi=
 \begin{cases}
  0+m_{\Vv}(\phi) &\text{if }\phi_{*,\xi}\Vv=\widevec{\phi(\xi)\xi} \\
  s_{\Vv}(\phi)+0 & \text{otherwise}
 \end{cases}=1\quad\text{for each }\Vv\in A,
\end{gather*}
so in particular $\#A=2$, but also $A\subset (T_\xi\sP^1)\setminus\{\widevec{\xi\phi(\xi)}\}$, so in particular
\begin{gather*}
\{\xi\}\sqcup\bigsqcup_{\Vv\in A}U(\Vv)\subset U(\widevec{\phi(\xi)\xi}).
\end{gather*}
Then by the non-archimedean argument principle \eqref{eq:FKN}
applied to $\phi$ at $\xi$, for each $\Vv\in A$ and every $j\ge 2$,
we have the inclusion
\begin{gather*}
 \phi^{-j}\Bigl(\{\xi\}\sqcup\bigsqcup_{\Vv\in A}U(\Vv)\Bigr)
\subset\bigsqcup_{\Vv\in A}U(\Vv)\subset U(\widevec{\phi(\xi)\xi})
\end{gather*}
by induction and compute
\begin{align*}
 \dep_{\Vv}\widetilde{(\phi^j)}_\xi&=\int_{\sP^1}\bigl(s_{\Vv}(\phi)+m_{\Vv}(\phi)\cdot 1_{U(\phi_{*,\xi}\Vv)}(\cdot)\bigr)\bigl((\phi^{j-1})^*\delta_\xi\bigr)(\cdot)\\
&=\begin{cases}
  0+m_{\Vv}(\phi)\cdot((\phi^{j-1})^*\delta_\xi)\bigl(U(\widevec{\phi(\xi)\xi})\bigr) &\text{if }\phi_{*,\xi}\Vv=\widevec{\phi(\xi)\xi}, \\
  s_{\Vv}(\phi)\cdot\bigl((\phi^{j-1})^*\delta_\xi\bigr)(\sP^1) +0& \text{otherwise}
 \end{cases}\\
&=
1\cdot\deg(\phi^{j-1})=\frac{\deg(\phi^j)}{2}
\Bigl(<\frac{\deg(\phi^j)+1}{2}\Bigr),
\end{align*}
so that $\phi^j(\xi)\neq\xi$, 
that for every $\Vv\in A$, $\widetilde{(\phi^j)}_{\xi}\equiv\widevec{\xi\phi^j(\xi)}\neq\Vv$, and that $\dep_{\Vu}\widetilde{(\phi^j)}_\xi=0$ for every $\Vu\in(T_\xi\sP^1)\setminus A$.
Hence we still have $\MHRL(\phi^j)=\{\xi\}$ for every $j\ge 2$ 
by the slope formula \eqref{eq:slope} in this case.

From now on, suppose that $\tilde{\phi}_\xi$ is bijective
or equivalently $\deg_\xi(\phi)=(\phi^*\delta_\xi)(\{\xi\})=1$ (and then by \eqref{eq:directional}, 
\begin{gather*}
 m_{\,\cdot\,}(\phi)\equiv\deg_\xi(\phi)(=1) 
\end{gather*}
on $T_\xi\sP^1$). Then by the non-archimedean argument principle \eqref{eq:FKN} and \eqref{eq:surplus},
we have $\#A=(\deg\phi)-(\phi^*\delta_\xi)(\{\xi\})=1$;
writing $A=\{\Vv_1\}$, we have $s_{\Vv_1}(\phi)=\dep_{\Vv_1}\tilde{\phi}_\xi=1$,
and by the slope formula \eqref{eq:slope}, we also have
\begin{gather}
 \Vv_1\not\in F=\Fix_{\tilde{\phi}_\xi}:=\bigl\{\Vv\in T_\xi\sP^1:\tilde{\phi}_\xi\Vv=\Vv\bigr\}.\label{eq:notfixed}
\end{gather}
We also note that $\#F\in\{1,2\}$ and that some $\Vv\in(T_\xi\sP^1)\setminus F$ is acyclic under $\tilde{\phi}_\xi$ if and only if $\tilde{\phi}_\xi$ is acyclic or equivalently of infinite order (in the monoid of all selfmaps of $T_\xi\sP^1$), recalling Fact \ref{th:algred}.

There are two subcases.

\begin{inparaenum}[(i)]
 \item \label{item:infinite}
 If $\Vv_1$ is acyclic under $\tilde{\phi}_\xi$,
then for every $j\ge 2$, we have $\widetilde{(\phi^j)}_\xi\Vv_1\neq\Vv_1$ (and still $((\phi^j)^*\delta_\xi)(\{\xi\})=1$), 
and by the non-archimedean argument principle \eqref{eq:FKN}
(or the composition formula \eqref{eq:depthcomp} applied to $\phi^j=\phi^{j-1}\circ\phi$
at $\xi$), we compute
 \begin{align}
 \dep_{\Vv_1}\widetilde{(\phi^j)}_\xi
&=\deg(\phi^{j-1})\cdot\bigl(\dep_{\Vv_1}\tilde{\phi}_\xi\bigr)
+\bigl(\deg_\xi(\phi)\bigr)\cdot\dep_{\phi_{*,\xi}\Vv_1}\widetilde{(\phi^{j-1})}_\xi\notag\\
&=\deg(\phi^{j-1})\cdot 1+1\cdot 0=\frac{\deg(\phi^j)}{2}\Bigl(<\frac{\deg(\phi^j)+1}{2}\Bigr),\tag{*}
 \end{align}
 which in turn yields, for every $\Vv\in(T_{\xi}\sP^1)\setminus\{\Vv_1\}$,
 \begin{gather*}
 \dep_{\Vv}\widetilde{(\phi^j)}_\xi\le
\deg(\phi^j)-\dep_{\Vv_1}\widetilde{(\phi^j)}_\xi-\bigl((\phi^j)^*\delta_{\xi}\bigr)(\{\xi\})=\frac{\deg(\phi^j)-2}{2}\Bigl(<\frac{\deg(\phi^j)-1}{2}\Bigr);
 \end{gather*}
in the line (*) above, unless $\dep_{\phi_{*,\xi}\Vv_1}\widetilde{(\phi^{j-1})}_\xi=0$,
by the non-archimedean argument principle \eqref{eq:FKN} applied to $\phi^{j-1}$ and $\phi_{*,\xi}\Vv_1$ and the above uniqueness of $\Vv_1\in A$,
we must have $(\phi^\ell)_{*,\xi}(\phi_{*,\xi}\Vv_1)=\Vv_1$ for some $\ell\in\{1,\ldots,j-1\}$, so that $\widetilde{(\phi^{\ell+1})}_\xi\Vv_1=(\phi^{\ell+1})_{*,\xi}\Vv_1=(\phi_{*,\xi})^{\ell+1}\Vv_1=\Vv_1$, which contradicts the acyclicity assumption of $\Vv_1$.

Hence for every $j\ge 2$, 
we still have $\MHRL(\phi^j)=\{\xi\}$
by the slope formula \eqref{eq:slope} in this case.

 \item\label{item:finite} Alternatively, assume that $\Char k\ne2$ and that $\Vv_1$ is cyclic under $\tilde{\phi}_\xi$ having the (minimal) period $p$ (and then $p>1$ by \eqref{eq:notfixed}). Then $\xi\in\sF(\phi)$ (see \cite[Theorem 2.1]{Kiwi14}, or \cite[Theorem 8.7]{BenedettoBook} in the present generality).

Denoting by $U$ the Berkovich Fatou component of $\phi$ containing $\xi(=\phi(\xi))$,
 the restriction $\phi:U\to\phi(U)=U$ is bijective (see Fact \ref{th:FJ}),
 and then $\#\partial U<\infty$, while the restriction
$\phi:\partial U\to\partial U\subset\sJ(\phi)$
 is a (possibly non-cyclic) permutation (see \cite[Theorem 9.7]{BenedettoBook} in the general setting here). Moreover, the Berkovich (maximal) ramification locus
\begin{gather*}
 \sR(\phi):=\{\xi\in\sP^1:\deg_\xi(\phi)>1\}
 \bigl(=\{\xi\in\sP^1:\deg_\xi(\phi)=2=\deg\phi\}\text{ in this quadratic case}\bigr)
\end{gather*}
is a closed and connected non-singleton set in $\sP^1$ (see \cite[Theorem 8.2]{Faber13topologyI} for more details), and 
since $\phi$ is tame under the assumption $\Char k\ne2$, 
we have $U\cap\sR(\phi)=\emptyset$; by the non-archimedean argument principle \eqref{eq:FKN}, \eqref{eq:directional}, and the slope formula \eqref{eq:slope} applied to $\phi$ at $\xi$,
we also have $\sR(\phi)\subset U(\Vv_1)$.  Let
$\xi_0\in\sR(\phi)\setminus\{\xi\}$ be the nearest point retraction of $\xi$
in $\sH^1$ to $\sR(\phi)$ (so that $\xi_0\neq\xi$ and that
$\sR(\phi)\subset\sP^1\setminus U(\widevec{\xi_0\xi})$)
as in the statement in Theorem \ref{th:quad}.
Then there is a unique $\xi_1\in\partial U$ in the half open interval $(\xi,\xi_0]$.
Since $\partial U$ is finite and $\phi:\partial U\to\partial U$ is a permutation,
$\xi_1$ is a periodic point of $\phi$, 
and then the exact period $q=q_\phi$ of $\xi_1$ under $\phi$ is divisible by $p$
since we have the isometric restriction $\phi^q:(\xi,\xi_1]\to(\xi,\xi_1]$ 
so in particular $(\phi^q)_*\Vv_1=\Vv_1$. 

By an argument similar to that in the previous case (\ref{item:infinite}), 
we first have $\MHRL(\phi^j)=\{\xi\}$ for every $j\in\{1,\ldots,p-1\}$
(no matter whether $\xi\in\sF(\phi)$). We next note that for every $j\ge p$, 
\begin{gather}
 \MHRL(\phi^j)\subset U(\widevec{\xi\xi_1});
 \label{eq:transout}
\end{gather}
indeed, for every $j\ge p$, using the composition formula \eqref{eq:depthcomp},
we recursively compute
\begin{align*} 
 \dep_{\Vv_1}\widetilde{(\phi^j)}_\xi
 &=(\deg\phi)\dep_{\Vv_1}\widetilde{(\phi^{j-1})}_\xi
 +\deg_\xi(\phi)\dep_{(\phi^{j-1})_{*,\xi}\Vv_1}\tilde{\phi}_\xi\\
 &=(\deg\phi)^j\cdot\dep_{\Vv_1}\widetilde{(\Id)}_\xi+
\sum_{\ell=0}^{j-1}(\deg\phi)^{j-1-\ell}\cdot 1\cdot\dep_{(\phi^\ell)_{*,\xi}\Vv_1}\tilde{\phi}_\xi\\
 &\begin{cases}
 \displaystyle=\frac{\deg(\phi^j)}{2} &\text{if }j=p,\\
 \displaystyle\ge\frac{\deg(\phi^j)}{2}+(\deg\phi)^{j-1-p} &\text{if }j\ge p+1, 
\end{cases}
\end{align*}
and then the slope formula \eqref{eq:slope}
applied to $\phi^j$ and $\widevec{\xi\xi_1}$ 
(together with the convexity of $\hypRes_{\phi^j}$)
 yields \eqref{eq:transout}.

There are still two subcases.

 \begin{inparaenum}[(\ref{item:finite}-a)]
 \item\label{item:nonend}
If $\xi_1\in(\xi,\xi_0)$, then
$\widetilde{(\phi^q)}_{\xi_1}=(\phi^q)_{*,\xi_1}$ is also bijective or equivalently
$((\phi^q)^*\delta_{\xi_1})(\{\xi_1\})=1$
(and then by \eqref{eq:directional},
\begin{gather*}
 m_{\,\cdot\,}(\phi)\equiv\deg_{\xi_1}(\phi)=1
\end{gather*}
on $T_{\xi_1}\sP^1$),
and then for every $\ell\in\{1,\ldots,q\}$,
\begin{gather}
\begin{aligned}
 & \phi_{*,\phi^{\ell-1}(\xi_1)}\text{ is bijective and }
 \phi_{*,\phi^{\ell-1}(\xi_1)}\widevec{\phi^{\ell-1}(\xi_1)\xi}=\widevec{\phi^\ell(\xi_1)\xi},\text{ so in particular}\\
&\widetilde{(\phi^q)}_{\xi_1}\widevec{\xi_1\xi}=\widevec{\xi_1\xi}
\quad\text{and }\widetilde{(\phi^j)}_{\xi_1}\equiv\widevec{\xi_1\xi}=\widevec{\xi_1\phi^j(\xi_1)}
\text{ for every }j\text{ not divisible by }q.
\end{aligned}
\label{eq:tangentbij}
\end{gather}
Noting also that $\sR(\phi)\subset U(\widevec{\xi_1\xi_0})$ in this case, by \eqref{eq:tangentbij} and the non-archimedean argument principle \eqref{eq:FKN}, we have
\begin{gather}
\begin{aligned}
 & \bigl\{\Vw\in T_{\xi_1}\sP^1:\dep_{\Vw}\tilde{\phi}_{\xi_1}>0\text{ (and then}=1)\bigr\}
 =\bigl\{\widevec{\xi_1\xi},\widevec{\xi_1\xi_0}\bigr\}\quad\text{and indeed}\\
 &\dep_{\widevec{\xi_1\xi_0}}(\tilde{\phi}_{\xi_1})=s_{\widevec{\xi_1\xi_0}}(\phi)=1,
\end{aligned}
\label{eq:holes}
\end{gather}
and recalling that $\xi_1\in\sJ(\phi)$, the direction $\widevec{\xi_1\xi_0}$ is acyclic under $\widetilde{(\phi^q)}_{\xi_1}=(\phi^q)_{*,\xi_1}$ (see \cite[Theorem 2.1]{Kiwi14}, or \cite[Theorem 8.7]{BenedettoBook} in the present generality).

Once this acyclicity assertion on $\widevec{\xi_1\xi_0}$ is at our disposal, we conclude with the following argument similar to that in the case (\ref{item:infinite}); for later use,
we show slightly more than what is needed in this proof.
First, for every $j\ge 1$, $\widetilde{(\phi^j)}_{\xi_1}\widevec{\xi_1\xi_0}\neq\widevec{\xi_1\xi_0}$, and by \eqref{eq:tangentbij}, \eqref{eq:holes}, and the non-archimedean argument principle \eqref{eq:FKN}, we compute
 \begin{align}
 \dep_{\widevec{\xi_1\xi_0}}\widetilde{(\phi^j)}_{\xi_1}
&=\int_{\sP^1}\bigl(s_{\widevec{\xi_1\xi_0}}(\phi)+m_{\widevec{\xi_1\xi_0}}(\phi)\cdot 1_{U(\phi_{*,\xi_1}\widevec{\xi_1\xi_0})}(\cdot)\bigr)\bigl((\phi^{j-1})^*\delta_{\xi_1}\bigr)(\cdot)\notag\\
&=\bigl(\dep_{\widevec{\xi_1\xi_0}}\tilde{\phi}_{\xi_1}\bigr)\cdot\deg(\phi^{j-1})
 +\bigl(\deg_{\xi_1}(\phi)\bigr)\cdot
\bigl((\phi^{j-1})^*\delta_{\xi_1}\bigr)(U(\phi_{*,\xi_1}\widevec{\xi_1\xi_0}))\notag\\
 &=1\cdot\deg(\phi^{j-1}) +1\cdot 0=\frac{\deg(\phi^j)}{2};\tag{**}
 \end{align}
in the line (**), unless
$\bigl((\phi^{j-1})^*\delta_{\xi_1}\bigr)(U(\phi_{*,\xi_1}\widevec{\xi_1\xi_0}))=0$,
by \eqref{eq:tangentbij} and the non-archimedean argument principle \eqref{eq:FKN} applied to $\phi_{*,\xi_1}\widevec{\xi_1\xi_0}$,
 we must have $j\ge q$ and $(\phi^\ell)_{*,\phi(\xi_1)}(\phi_{*,\xi_1}\widevec{\xi_1\xi_0})=\widevec{\xi_1\xi_0}$ for some $\ell\in\{q-1,\ldots,j-1\}$, so that ($q|(\ell+1)$ and)
  $((\phi^q)_{*,\xi_1})^{(\ell+1)/q}\widevec{\xi_1\xi_0}=\widevec{\xi_1\xi_0}$, which contradicts the above acyclicity assertion on $\widevec{\xi_1\xi_0}$.

Next, pick any $1\le j<q$, and recall that
$\widetilde{(\phi^j)}_{\xi_1}\equiv\widevec{\xi_1\xi}$.
Then similarly to the above computation 
for the direction $\widevec{\xi_1\xi_0}$,
we inductively compute 
\begin{align*}
 \dep_{\widevec{\xi_1\xi}}\widetilde{(\phi^j)}_{\xi_1}
&=\int_{\sP^1}\bigl(s_{\widevec{\xi_1\xi}}(\phi)
+m_{\widevec{\xi_1\xi}}(\phi)\cdot 1_{U(\widevec{\phi(\xi_1)\xi})}(\cdot)\bigr)\bigl((\phi^{j-1})^*\delta_{\xi_1}\bigr)(\cdot)\notag\\
&=0\cdot\deg(\phi^{j-1})
 +m_{\widevec{\xi_1\xi}}(\phi)\cdot\bigl((\phi^{j-1})^*\delta_{\xi_1}\bigr)(U(\widevec{\phi(\xi_1)\xi}))\notag\\
 &\biggl(=0+\int_{\sP^1}\bigl(s_{\widevec{\phi(\xi_1)\xi}}(\phi)+m_{\widevec{\phi(\xi_1)\xi}}(\phi)\cdot 1_{U(\widevec{\phi^2(\xi_1)\xi})}(\cdot)\bigr)\bigl((\phi^{j-2})^*\delta_{\xi_1}\bigr)(\cdot)\quad(\text{if }j>1)\biggr)\notag\\
&=0+1\cdot\sum_{\ell=1}^{j-1}\deg(\phi^{\ell-1})
+m_{\widevec{\phi^{j-1}(\xi_1)\xi}}(\phi)\cdot
\delta_{\xi_1}(U(\widevec{\phi^j(\xi_1)\xi}))\notag\\
&=0+\frac{2^{j-1}-1}{2-1}+1=\frac{\deg(\phi^j)}{2}-\bigl((\phi^j)^*\delta_{\xi_1}\bigr)(\{\xi_1\}),
\end{align*}
where for every $\ell\in\{1,\ldots,q-1\}$, 
by \eqref{eq:tangentbij}, we have
$\sR(\phi)\subset U(\widevec{\phi^\ell(\xi_1)\xi})=U(\phi_*^\ell\widevec{\xi_1\xi})$,
and then by the non-archimedean argument principle \eqref{eq:FKN}, we have
$s_{\phi_*^\ell\widevec{\xi_1\xi}}(\phi)>0$, and indeed $=1$. 
Hence also by the above computation of $\dep_{\widevec{\xi_1\xi_0}}\widetilde{(\phi^j)}_{\xi_1}$, we have
\begin{multline*}
 \sum_{\Vw\in(T_{\xi_1}\sP^1)\setminus\{\widevec{\xi_1\xi_0},\widevec{\xi_1\xi}\}}
 \dep_{\Vw}\widetilde{(\phi^j)}_{\xi_1}\\
=\deg(\phi^j)-
\dep_{\widevec{\xi_1\xi_0}}\widetilde{(\phi^j)}_{\xi_1}
-\Bigl(\dep_{\widevec{\xi_1\xi}}\widetilde{(\phi^j)}_{\xi_1}
+\bigl((\phi^j)^*\delta_{\xi_1}\bigr)(\{\xi_1\})\Bigr)
=0,
\end{multline*}
and moreover, we have the opposite inclusion
$\MHRL(\phi^j)\subset U(\widevec{\xi_1\xi})$
to \eqref{eq:transout} by the slope formula \eqref{eq:slope} applied to $\phi^j$ 
and $\widevec{\xi_1\xi}$ (and the convexity of $\hypRes_{\phi^j}$).

Finally, pick any $j\ge q$, and
recall \eqref{eq:tangentbij}, \eqref{eq:holes}. 
If $q\mid j$ or equivalently $\bigl((\phi^j)^*\delta_{\xi_1}\bigr)(\{\xi_1\})=1$,
then by the above computation of $\dep_{\widevec{\xi_1\xi_0}}\widetilde{(\phi^j)}_{\xi_1}$,
we have
\begin{gather*}
 \sup_{\Vw\in(T_{\xi_1}\sP^1)\setminus\{\widevec{\xi_1\xi_0}\}}
 \dep_{\Vw}\widetilde{(\phi^j)}_{\xi_1}
 \le
\deg(\phi^j)-\bigl((\phi^j)^*\delta_{\xi_1}\bigr)(\{\xi_1\})
-\dep_{\widevec{\xi_1\xi_0}}\widetilde{(\phi^j)}_{\xi_1}
 =\frac{\deg(\phi^j)}{2}-1.
\end{gather*}
Alternatively if $q\nmid j$
or equivalently $\widetilde{(\phi^j)}_{\xi_1}\equiv\widevec{\xi_1\xi}$, then 
by the acyclicity assertion on $\widevec{\xi_1\xi_0}$ 
and the non-archimedean argument principle \eqref{eq:FKN}, 
we have $(\phi^q)_{*,\xi_1}^{-\lfloor j/q \rfloor}\widevec{\xi_1\xi_0}
\in T_{\xi_1}\sP^1\setminus\{\widevec{\xi_1\xi_0},\widevec{\xi_1\xi}\}$
and $\dep_{(\phi^q)_{*,\xi_1}^{-\lfloor j/q \rfloor}\widevec{\xi_1\xi_0}}\widetilde{(\phi^j)}_{\xi_1}\ge s_{(\phi^q)_{*,\xi_1}^{-\lfloor j/q \rfloor}\widevec{\xi_1\xi_0}}(\phi^j)\ge s_{\widevec{\xi_1\xi_0}}(\phi)=1$,
which together with the above computation of $\dep_{\widevec{\xi_1\xi_0}}\widetilde{(\phi^j)}_{\xi_1}$ yield
\begin{gather*}
\dep_{\widevec{\xi_1\xi}}\widetilde{(\phi^j)}_{\xi_1}
\le\deg(\phi^j)
-\dep_{\widevec{\xi_1\xi_0}}\widetilde{(\phi^j)}_{\xi_1}
-\dep_{(\phi^q)_{*,\xi_1}^{-\lfloor j/q \rfloor}\widevec{\xi_1\xi_0}}\widetilde{(\phi^j)}_{\xi_1} 
\le\frac{\deg(\phi^j)}{2}-1\quad\text{and}\\
\sup_{\Vw\in(T_{\xi_1}\sP^1)\setminus\{\widevec{\xi_1\xi_0}\}}\dep_{\Vw}\widetilde{(\phi^j)}_{\xi_1}
\le\deg(\phi^j)-\dep_{\widevec{\xi_1\xi_0}}\widetilde{(\phi^j)}_{\xi_1}
\displaystyle=\frac{\deg(\phi^j)}{2}.
\end{gather*} 
Hence for every $j\ge q$, the minimum locus of $\hypRes_{\phi^j}$ equals $\{\xi_1\}$ 
by the slope formula \eqref{eq:slope} applied to $\phi^j$ at $\xi_1$.

\item\label{item:boundary}
Suppose finally that $\xi_1=\xi_0$. Since $\xi_1\in\sR(\phi)$, we have $\deg_{\xi_1}(\phi)=2$; moreover, $\xi_1$ has period $q$ under $\phi$ and $\deg_{\,\cdot\,}(\phi)\equiv1$ on $\partial U\setminus\{\xi_1\}$. Hence
 $\widetilde{(\phi^q)}_{\xi_1}$ is $2$ to $1$
 or equivalently $((\phi^q)^*\delta_{\xi_1})(\{\xi_1\})=\deg\phi=2$,
 and then $s_{\,\cdot\,}(\phi)\equiv 0$ on $T_{\xi_1}\sP^1$ by \eqref{eq:surplus};
more precisely,
 \begin{gather}
\begin{cases}
  \phi_{*,\xi_1}:T_{\xi_1}\sP^1\to T_{\phi(\xi_1)}\sP^1
 \text{ is }2\text{ to }1\text{ and }\phi_{*,\xi_1}\widevec{\xi_1\xi}=\widevec{\phi(\xi_1)\xi}\bigl(=\widevec{\phi(\xi_1)\xi_1}\bigr)\quad\text{and}\\
 \text{for every }\ell\in\{2,\ldots,q\},
 \phi_{*,\phi^{\ell-1}(\xi_1)}\text{ is bijective and }
 \phi_{*,\phi^{\ell-1}(\xi_1)}\widevec{\phi^{\ell-1}(\xi_1)\xi}=\widevec{\phi^\ell(\xi_1)\xi},
\end{cases}\label{eq:tangent21}
\end{gather}
so in particular $\deg_{\phi(\xi_1)}(\phi^{q-1})=1$ and,
by \eqref{eq:directional},
there is a unique direction $\Vw_1\in(T_{\xi_1}\sP^1)\setminus\{\widevec{\xi_1\xi}\}$ such that
$\phi_{*,\xi_1}\Vw_1=\widevec{\phi(\xi_1)\xi}$.
Then by the non-archimedean argument principle \eqref{eq:FKN},
we have
\begin{gather*}
 \bigl\{\Vw\in T_{\xi_1}\sP^1:\dep_{\Vw}\tilde{\phi}_{\xi_1}>0\text{ (and then }\dep_{\Vw}\tilde{\phi}_{\xi_1}=m_{\Vw}(\phi)=1)\bigr\}
=\bigl\{\widevec{\xi_1\xi},\Vw_1\bigr\}.
\end{gather*}
For every $j\ge 0$, set
\begin{gather*}
 A_j:=\dep_{\Vw_1}\widetilde{(\phi^j)}_{\xi_1},\quad
 B_j:=
 \sum_{\Vw\in(T_{\xi_1}\sP^1)\setminus\{\widevec{\xi_1\xi},\Vw_1\}}
 \dep_{\Vw}\widetilde{(\phi^j)}_{\xi_1},\quad\text{and}\quad
 C_j:=\dep_{\widevec{\xi_1\xi}}\widetilde{(\phi^j)}_{\xi_1}
\end{gather*}
(so that $A_0=B_0=C_0=0$). Then by the non-archimedean argument principle
\eqref{eq:FKN} and \eqref{eq:tangent21}, for every $j\ge1$, we have
the equalities
\begin{gather}
\begin{aligned}
 A_j
 =&\int_{\sP^1}\bigl(s_{\Vw_1}(\phi)+m_{\Vw_1}(\phi)\cdot
 1_{U(\widevec{\phi(\xi_1)\xi})}\bigr)(\phi^{j-1})^*\delta_{\xi_1}
 =0+1\cdot\bigl((\phi^{j-1})^*\delta_{\xi_1}\bigr)
 \bigl(U(\widevec{\phi(\xi_1)\xi})\bigr)\\
 =&\int_{\sP^1}\bigl(s_{\widevec{\xi_1\xi}}(\phi)
 +m_{\widevec{\xi_1\xi}}(\phi)\cdot
 1_{U(\widevec{\phi(\xi_1)\xi})}\bigr)(\phi^{j-1})^*\delta_{\xi_1}
 =C_j,
\end{aligned}\label{eq:A}\\
\begin{aligned}
 B_j
 &=\sum_{\Vw\in(T_{\xi_1}\sP^1)\setminus\{\widevec{\xi_1\xi},\Vw_1\}}
 \int_{\sP^1}\bigl(s_{\Vw}(\phi)+m_{\Vw}(\phi)\cdot
 1_{U(\phi_{*,\xi_1}\Vw)}\bigr)(\phi^{j-1})^*\delta_{\xi_1},\\
 &\hspace*{250pt}\text{which is }\equiv 0\text{ if }j<q,\quad\text{and}
\end{aligned}\label{eq:Borig}\\
\begin{aligned}
  &A_j+B_j+C_j
 =\deg(\phi^j)-\bigl((\phi^j)^*\delta_{\xi_1}\bigr)(\{\xi_1\})
 =\deg(\phi^j)-
 \begin{cases}
  2^{j/q} & \text{if }q\mid j,\\
  0 & \text{if }q\nmid j,
 \end{cases}\\
&\hspace*{250pt}\text{so }A_j=C_j=\frac{\deg(\phi^j)}{2}\text{ if }j<q.
\end{aligned}  
\label{eq:C}
\end{gather}
In particular,
for every $p\le j<q$, recalling that
$\widetilde{(\phi^j)}_{\xi_1}\equiv\widevec{\xi_1\xi}$,
we have the opposite inclusion
$\MHRL(\phi^j)\subset U(\widevec{\xi_1\xi})$
to \eqref{eq:transout} by the slope formula \eqref{eq:slope} applied to $\phi^j$ 
and $\widevec{\xi_1\xi}$ (and the convexity of $\hypRes_{\phi^j}$).

Finally, pick any $j\ge q$. Then we also have the recursive equality
\begin{gather}
\begin{aligned}
  B_j
&=\sum_{\Vw\in(T_{\xi_1}\sP^1)\setminus\{\widevec{\xi_1\xi},\Vw_1\}}
\int_{\sP^1}\bigl(s_{\Vw}(\phi)+m_{\Vw}(\phi)\cdot 1_{U(\phi_{*,\xi_1}\Vw)}\bigr)(\phi^{j-1})^*\delta_{\xi_1}\quad(\text{seen in \eqref{eq:Borig}})\\
&=0+\sum_{\Vu\in(T_{\phi(\xi_1)}\sP^1)\setminus\{\widevec{\phi(\xi_1)\xi}\}}\int_{\sP^1}
\bigl(\deg_{\xi_1}(\phi)\cdot 1_{U(\Vu)}\bigr)(\phi^{j-1})^*\delta_{\xi_1}\quad(\text{also using \eqref{eq:directional}})\\
&=\bigl(\deg_{\xi_1}(\phi)\bigr)\cdot\sum_{\Vu\in(T_{\phi(\xi_1)}\sP^1)\setminus\{\widevec{\phi(\xi_1)\xi}\}}\int_{\sP^1}\Bigl((\phi^{q-1})_*1_{U(\Vu)}\Bigr)(\phi^{j-q})^*\delta_{\xi_1}\\
&=\bigl(\deg_{\xi_1}(\phi)\bigr)\cdot\sum_{\Vw\in(T_{\xi_1}\sP^1)\setminus\{\widevec{\xi_1\xi}\}}\int_{\sP^1}1_{U(\Vw)}(\phi^{j-q})^*\delta_{\xi_1}
=2\cdot(A_{j-q}+B_{j-q}),
\end{aligned}\label{eq:B}
\end{gather}
so that from \eqref{eq:B}, \eqref{eq:A}, and \eqref{eq:C}, we have
\begin{gather*}
 B_j-B_{j-q}=2A_{j-q}+B_{j-q}
=2^{j-q}-
\begin{cases}
  2^{(j-q)/q} & \text{if }q\mid j,\\
 0 & \text{if } q\nmid j,
\end{cases}
\end{gather*}
and then also using \eqref{eq:Borig}, 
\begin{gather}
 B_j
\le B_{j-\lfloor\frac{j}{q}\rfloor q}
+\sum_{\ell=1}^{\lfloor\frac{j}{q}\rfloor}2^{j-\ell q}
=0+2^j\cdot\frac{2^{-q}\bigl(1-2^{-q\lfloor\frac{j}{q}\rfloor}\bigr)}{1-2^{-q}}
=\frac{2^j-2^{j-\lfloor\frac{j}{q}\rfloor q}}{2^q-1}<\frac{\deg(\phi^j)-1}{2}\quad\text{and}
\label{eq:outward} \\
B_j\begin{cases}
    =0 &\text{if }j=q,\\
    \ge 2 &\text{if }j>q.
   \end{cases}\notag
\end{gather}
Hence by \eqref{eq:A} and \eqref{eq:C}, we also have
\begin{gather*}
 A_j=C_j
\begin{cases}
 \le\displaystyle\frac{1}{2}\biggl(2^j-B_j-
 \begin{cases}
  2^{j/q} & \text{if }q|j\\
 0 & \text{if } q\nmid j
 \end{cases}\biggr)<\displaystyle\frac{\deg(\phi^j)-1}{2} &\text{if }j>q,\\
 =\displaystyle\frac{\deg(\phi^q)-2}{2} &\text{if }j=q.
\end{cases}
\end{gather*}
Hence for every $j\ge q$,
the minimum locus of $\hypRes_{\phi^j}$ equals $\{\xi_1\}$
by the slope formula \eqref{eq:slope} applied to $\phi^j$ at $\xi_1$.
\end{inparaenum}
\end{inparaenum}
Now the proof of Theorem \ref{th:quad} is complete. \qed
  
\begin{remark}[The quadratic polynomial case]
For every quadratic polynomial $\phi=P$, by \cite[Theorem 4]{NO24},
the minimum locus of $\hypRes_{P^j}$ equals $\{\xi_P\}$ for every integer
$j\ge (\deg P)-1=1$. Hence the finite-order case in Theorem \ref{th:quad}
never occurs for $P$. Indeed, $\deg_{\xi_P}P=2$ by
\cite[Lemma 4.1]{NO24}, so the case $(P^*\delta_{\xi_P})(\{\xi_P\})=1$
in the proof of Theorem \ref{th:quad} cannot occur.
\end{remark}

\section{Proof of Theorem \ref{th:residuemoduli}}\label{sec:residuemoduli}

For a quadratic rational function $\phi\in K(z)$ in $\Char k\neq 2$ and
the order $p>1$ case of Theorem \ref{th:quad}, 
we retain the notations $\xi_\phi$, $U$, $\xi_0$, $\xi_1$, $q$ for $\phi$
in Theorem \ref{th:quad} and the notations $\Vw_0,\Vw_\infty$, 
and also $\Vw_1,\Vw_2$ when $\xi_1=\xi_0$, 
prepared in Section \ref{sec:residue}. 

As seen in (\ref{item:finite}) in Section \ref{sec:stationary},
the Berkovich (indeed maximal) ramification locus
$\sR(\phi)$ of $\phi$ is in $\sP^1\setminus U(\widevec{\xi_0\xi_\phi})$.
By the computations in (\textit{\ref{item:nonend}}) and
(\textit{\ref{item:boundary}}) in Section \ref{sec:stationary},
for every $1\le r<q$, the intrinsic reduction data of $\phi^r$ at $\xi_1$ is
\begin{equation}
 \begin{aligned}
  &\widetilde{(\phi^r)}_{\xi_1}\equiv\Vw_\infty,\quad\text{or more generally, }
\widetilde{(\phi^j)}_{\xi_1}\equiv\Vw_\infty\text{ for every }j\ge 1\text{ such that }q\nmid j,\\
  &\dep_{\Vw_0}\widetilde{(\phi^r)}_{\xi_1}=
  \dep_{\Vw_\infty}\widetilde{(\phi^r)}_{\xi_1}=\frac{\deg(\phi^r)}{2},\quad\text{and}\\
  &\dep_{\Vw}\widetilde{(\phi^r)}_{\xi_1}=0
  \quad\text{for every }\Vw\notin\{\Vw_0,\Vw_\infty\},
 \end{aligned}
 \label{eq:short}
\end{equation}
and there are two possibilities for that of $\phi^q$ at $\xi_1$;
if $\xi_1\ne\xi_0$, then
\begin{equation}
 \begin{aligned}
  &\widetilde{(\phi^q)}_{\xi_1}\text{ is bijective},\\
&\widetilde{(\phi^q)}_{\xi_1}(\Vw_\infty)=\Vw_\infty\text{ but }
  \Vw_0\text{ is acyclic under }\widetilde{(\phi^q)}_{\xi_1},\\
  &\dep_{\Vw_0}\widetilde{(\phi^q)}_{\xi_1}
=\frac{\deg(\phi^q)}{2}
>\frac{\deg(\phi^q)-2}{2}
=\dep_{\Vw_\infty}\widetilde{(\phi^q)}_{\xi_1}>0,
\quad\text{and}\\
  &\dep_{\Vw}\widetilde{(\phi^q)}_{\xi_1}=0
  \quad\text{for every }\Vw\notin\{\Vw_0,\Vw_\infty\}.
 \end{aligned}
 \label{eq:qinf}
\end{equation}
Indeed, in (\textit{\ref{item:nonend}}),
$\Vw_0$ and $\Vw_\infty$ correspond respectively to
$\widevec{\xi_1\xi_0}$ and $\widevec{\xi_1\xi}$ there, and \eqref{eq:qinf} is seen there.
If $\xi_1=\xi_0$, then
\begin{equation}
\begin{aligned}
 &\widetilde{(\phi^q)}_{\xi_1}\text{ is $2$ to $1$},\\
&\widetilde{(\phi^q)}_{\xi_1}(\Vw_\infty)
 =\Vw_\infty=\widetilde{(\phi^q)}_{\xi_1}(\Vw_0)
,\text{ so $\Vw_0$ is not acyclic under }\widetilde{(\phi^q)}_{\xi_1},\\
 &\dep_{\Vw_0}\widetilde{(\phi^q)}_{\xi_1}
=\dep_{\Vw_\infty}\widetilde{(\phi^q)}_{\xi_1}=\frac{\deg(\phi^q)-2}{2},\quad\text{and}\\
 &\dep_{\Vw}\widetilde{(\phi^q)}_{\xi_1}=0
 \quad\text{for every }\Vw\notin\{\Vw_0,\Vw_\infty\}.
\end{aligned}
\label{eq:qfin}
\end{equation}
Indeed, in (\textit{\ref{item:boundary}}),
$\Vw_\infty$ and $\Vw_0$ correspond respectively to
$\widevec{\xi_1\xi}$ and $\Vw_1$ there, and \eqref{eq:qfin} is seen there.
Moreover, if $\xi_1\neq\xi_0$, then for every $j>q$, it is also seen 
in (\textit{\ref{item:nonend}}) that
\begin{equation}
\begin{aligned}
 &\dep_{\Vw_\infty}\widetilde{(\phi^j)}_{\xi_1}
 <\frac{\deg(\phi^j)}{2}=\dep_{\Vw_0}\widetilde{(\phi^j)}_{\xi_1},\\
 &\sup_{\Vw\in T_{\xi_1}\sP^1\setminus\{\Vw_0\}}
 \dep_{\Vw}\widetilde{(\phi^j)}_{\xi_1}
 \le\frac{\deg(\phi^j)}{2}-
\begin{cases}
 1 & \text{if }q\mid j,\\
 0 & \text{if }q\nmid j.
\end{cases} 
\end{aligned}
\label{eq:nonenddepth}
\end{equation}

To classify the intrinsic reduction data of $\phi^q$ at $\xi_1$, we use the natural
$k$-projective-line structure on $T_{\xi_1}\sP^1$.  Set
\begin{gather*}
 \bT_\phi
 :=\bigl\{
 h\in\operatorname{Aut}_k(T_{\xi_1}\sP^1):
 h(\Vw_0)=\Vw_0,\ 
 h(\Vw_\infty)=\Vw_\infty
 \bigr\}
 \cong\mathbb{G}_m.
\end{gather*}
Recall the definition \eqref{eq:rdef} of the intrinsic prehole parameter
$\ph(\phi)$ of $\phi$; 
if $\xi_1=\xi_0$, then 
\begin{gather*}
 D_\phi:=\bigl(\widetilde{(\phi^q)}_{\xi_1}\bigr)^*(\Vw_0)=\Vw_1+\Vw_2 
\end{gather*}
(see \eqref{eq:pullbackdiv}) as the divisor on $T_{\xi_1}\sP^1$
and 
\begin{gather}
 \ph(\phi):=2+c_\phi+(c_\phi)^{-1},
 \label{eq:cross}
\end{gather}
where the (symmetrized) cross ratio is
\begin{equation*}
 \begin{aligned}
 c_\phi+(c_\phi)^{-1}
 &=[\Vw_0,\Vw_\infty;\Vw_1,\Vw_2]
  +[\Vw_0,\Vw_\infty;\Vw_2,\Vw_1]\\
 &=\frac{(\Vw_1-\Vw_0)(\Vw_2-\Vw_\infty)}
 {(\Vw_1-\Vw_\infty)(\Vw_2-\Vw_0)}
 +\frac{(\Vw_2-\Vw_0)(\Vw_1-\Vw_\infty)}
 {(\Vw_2-\Vw_\infty)(\Vw_1-\Vw_0)}.
 \end{aligned}
\end{equation*}
If $\xi_1\neq\xi_0$, $\ph(\phi):=\infty$.

The following observation plays a key role in the proof of Theorem \ref{th:residuemoduli}.

\begin{proposition}\label{th:multiple}
Suppose that $\phi$ is in the order $p>1$ case of Theorem \ref{th:quad}.
Then the fixed point (direction) $\Vw_\infty$ of $\widetilde{(\phi^q)}_{\xi_1}$ 
is multiple or equivalently
the multiplier of $\widetilde{(\phi^q)}_{\xi_1}$ at $\Vw_\infty$
equals $1$.
\end{proposition}

\begin{proof}
It is shown by a computation of non-archimedean analytic function on an annulus.
For a simpler argument, we use Corollary \ref{th:fixednormal} and Proposition
\ref{th:typeIInormal} in Appendix \ref{sec:normalformsappendix} and 
the notations therein.

We note that, for every $\xi\in[\xi_\phi,\xi_1)$,
$\bigl((\phi^q)^*\delta_\xi\bigr)(\{\xi\})=1$, and
if $\xi\in(\xi_\phi,\xi_1)$, then
the directions $\widevec{\xi\xi_\phi},\widevec{\xi\xi_1}$
are fixed by $\widetilde{(\phi^q)}_{\xi}$,
and $\widetilde{(\phi^q)}_{\xi_\phi}
=\widetilde{(\phi^p)}_{\xi_\phi}^{q/p}=\Id$
since $\widetilde{\phi}_{\xi_\phi}$ has order $p$ and $p\mid q$. 
Replacing the projective coordinate of $\bP^1$ if necessary, 
suppose here that
\begin{gather*}
 \xi_1=\xi_g\quad\text{and}\quad \xi_\phi\in(0,\xi_g).
\end{gather*}

First, for every $\xi\in[\xi_\phi,\xi_g)\cap\sH^1_{\mathrm{II}}$, 
choose $M_\xi\in\Aff(2,K)$ fixing both $0$ and $\infty$ such that
$M_\xi(\xi)=\xi_g$. By Corollary \ref{th:fixednormal}, denoting the $a$ and
$h$ there by $a_\xi$ and $h_\xi$, we have the normal form
\begin{gather*}
 (M_\xi\circ\phi^q\circ M_\xi^{-1})(z)
 =a_\xi z\bigl(1+h_\xi(z)\bigr)\quad\text{in }K(z),
 \quad |a_\xi|=1,\text{ and }|h_\xi|_{\xi_g}<1,
\end{gather*}
so that
\begin{gather*}
 \phi^q(z)=a_\xi z\bigl(1+(h_\xi\circ M_\xi)(z)\bigr)\quad\text{in }K(z),
 \quad |a_\xi|=1,\text{ and }|h_\xi\circ M_\xi|_{\xi}=|h_\xi|_{\xi_g}<1.
\end{gather*}
Similarly, for every $\xi\in[\xi_\phi,\xi_g)\cap\sH^1_{\mathrm{III}}$, 
by Corollary \ref{th:fixednormal}, setting $M_\xi=\Id$, we have the normal form
\begin{gather*}
 \phi^q(z)=a_\xi z\bigl(1+(h_\xi\circ M_\xi)(z)\bigr)\quad\text{in }K(z),
 \quad |a_\xi|=1,\text{ and }|h_\xi\circ M_\xi|_{\xi}<1.
\end{gather*}
We observe that $\widehat{a_{\,\cdot}}\equiv 1$ on $[\xi_\phi,\xi_g)$;
indeed, for each $\xi\in[\xi_\phi,\xi_g)$, 
by the continuity of $|h_\xi\circ M_\xi|_{\cdot}$ on $\sP^1$,
we have 
\begin{gather*}
 \left|\frac{a_\cdot}{a_\xi}-1\right|
 =\left|\frac{1+h_\xi\circ M_\xi}{1+h_\cdot\circ M_\cdot}-1\right|_\cdot
 =\left|\frac{h_\xi\circ M_\xi-h_\cdot\circ M_\cdot}{1+h_\cdot\circ M_\cdot}\right|_\cdot<1,
 \quad\text{so that}\quad
 \widehat{\left(\frac{a_\cdot}{a_\xi}\right)}=1\quad\text{near }\xi,
\end{gather*}
that is, $\widehat{a_{\,\cdot}}$ is locally constant 
on the interval $[\xi_\phi,\xi_g)$ in $\sH^1$, and 
$\widehat{a_{\xi_\phi}}=1$ by the following computation
\begin{multline*}
1\cdot\zeta=\bigl((M_{\xi_\phi})_*\circ\Id\circ(M_{\xi_\phi})_*^{-1}\bigr)(\zeta)
= 
\bigl((M_{\xi_\phi})_*\circ\widetilde{(\phi^q)}_{\xi_\phi}\circ(M_{\xi_\phi})_*^{-1}\bigr)(\zeta)\\
= 
 \widetilde{(M_{\xi_\phi}\circ\phi^q\circ M_{\xi_\phi}^{-1})}(\zeta)
=
\widehat{a_{\xi_\phi}}\cdot\zeta
\quad\text{in } k(\zeta).
\end{multline*}
Next, noting \eqref{eq:qinf} and \eqref{eq:qfin},
by Proposition \ref{th:typeIInormal},
we have the normal form
\begin{gather*}
 \phi^q(z)
 =a z\frac{\prod_{i=1}^{\ell-1}(z-\alpha_i)}
 {\prod_{j=1}^{k}(z-\beta_j)}\bigl(1+h(z)\bigr)\quad\text{in }K(z)
 \quad\text{and}\quad \max\{\ell,k\}=\deg_{\xi_g}(\phi^q),
\end{gather*}
of $\phi^q$ at $\xi_g$, where
$\ell\ge1$, $k\ge0$, $|a|=|\alpha_i|=|\beta_j|=1$,
$\widehat{\alpha_i}\ne\widehat{\beta_j}$
for every $i,j$, and $|h|_{\xi_g}<1$, so that
\begin{gather*}
 \widetilde{(\phi^q)}(\zeta)
 =\widehat a\,\zeta
 \frac{\prod_{i=1}^{\ell-1}(\zeta-\widehat{\alpha_i})}
 {\prod_{j=1}^{k}(\zeta-\widehat{\beta_j})}.
\end{gather*}
Consequently, 
by the continuity of $|h|_{\cdot}$ on $\sP^1$,
choosing $\xi\in[\xi_{\phi},\xi_g)\cap\sH^1_{\mathrm{II}}$ close enough to $\xi_g$, we have
$|h\circ M_\xi^{-1}|_{\xi_g}<1$, and then also by the above observation
on $\widehat{a_\cdot}$,
\begin{gather*}
 1\cdot\zeta=\widehat{a_\xi}\zeta
 =\widetilde{(M_\xi\circ\phi^q\circ M_\xi^{-1})}(\zeta)
 =\widehat{\biggl(a\cdot(-1)^{\ell-1-k}
 \frac{\prod_{i=1}^{\ell-1}\alpha_i}{\prod_{j=1}^{k}\beta_j}\biggr)}\zeta
 =\Bigl(\frac{\rd\widetilde{(\phi^q)}}{\rd\zeta}(\zeta)
 \Bigm|_{\zeta=0}\Bigr)\zeta\quad\text{in } k(\zeta),
\end{gather*}
which completes the proof.
\end{proof}

\begin{remark}
 In the latter half of the above proof, replacing further the projective coordinate of $\bP^1$ so in addition that $\Vw_0=\infty\in\bP^1(k)$ if necessary, we have even
\begin{gather*}
 (k,\ell)=
 \begin{cases}
  (1,1) & \text{if }\xi_1\ne\xi_0,\\
  (2,1) & \text{if }\xi_1=\xi_0
 \end{cases}
\end{gather*}
noting \eqref{eq:qinf} and \eqref{eq:qfin} and 
replacing also the constants and functions in the above normal forms if necessary.
\end{remark}

The following proposition will also be used in the proof of Theorem
\ref{th:arc}. 

\begin{proposition}[Recovery from the pullback divisor]\label{th:divisor}
Suppose that $\ph(\phi)\in k$.
Then the pullback divisor $D_\phi$ recovers
$\widetilde{(\phi^q)}_{\xi_1}$ equivariantly under
$\bT_\phi$-conjugation of it;
more precisely, choosing $M\in\PGL(2,K)$ such that
\begin{gather*}
 M(\xi_1)=\xi_g,\quad
 M_*(\Vw_0)=\widevec{\xi_g0},\quad\text{and that}\quad
 M_*(\Vw_\infty)=\widevec{\xi_g\infty},
\end{gather*}
letting $[X:Y]$ be the homogeneous coordinates on
$T_{\xi_g}\sP^1\cong\bP^1(k)$ with $\zeta=Y/X$, and writing
\begin{gather*}
 M_*D_\phi=(AX^2+BXY+CY^2=0)\quad\text{and}\quad AC\ne0,\quad\text{on }\bP^1(k),
\end{gather*}
the intrinsic prehole parameter of $\phi$ is computed as
\begin{gather}
 \ph(\phi)=\frac{B^2}{AC}
 \label{eq:rdiv}
\end{gather}
and moreover
\begin{gather}
 \widetilde{(M\circ\phi^q\circ M^{-1})}(\zeta)
 =\frac{A+B\zeta+C\zeta^2}{C\zeta}\quad\text{in }k(\zeta).
 \label{eq:return}
\end{gather}
\end{proposition}

\begin{proof}
Fix $M$ and homogeneous coordinates $[X:Y]$ on $T_{\xi_g}\sP^1$
as in the statement and write $M_*D_\phi=(A+B\zeta+C\zeta^2=0)=x_1+x_2$, $AC\neq 0$.
Then we compute
\begin{gather*}
 \ph(\phi)
 =2+\frac{x_1}{x_2}+\frac{x_2}{x_1}
=\frac{(x_1+x_2)^2}{x_1x_2}=\frac{(B/C)^2}{A/C}
 =\frac{B^2}{AC}
\end{gather*}
and, moreover, recalling \eqref{eq:qfin}, there is $c\in k^\times$ such that
\begin{gather*}
 \widetilde{(M\circ\phi^q\circ M^{-1})}(\zeta)
 =\frac{A+B\zeta+C\zeta^2}{c\zeta}\quad\text{in }k(\zeta).
\end{gather*}
Hence, by Proposition \ref{th:multiple}, we compute the multiplier of $\widetilde{(M\circ\phi^q\circ M^{-1})}$
at the fixed point $\infty$ as
\begin{gather*}
 1=\frac{\rd}{\rd\zeta}
 \biggl(\frac{1}{\widetilde{(M\circ\phi^q\circ M^{-1})}(1/\zeta)}\biggr)\biggm|_{\zeta=0}
 =\frac{\rd}{\rd\zeta}\biggl(\frac{c\zeta}{A\zeta^2+B\zeta+C}\biggr)\biggm|_{\zeta=0} 
 =\frac{cC}{C^2}=\frac{c}{C}, 
\end{gather*}
so that $c=C$. The equivariance under $\bT_\phi$-conjugation is clear from \eqref{eq:return}.
\end{proof}

Now let $\phi,\psi\in K(z)$ be the quadratic rational functions in the
order $p>1$ case of Theorem \ref{th:quad} with $q_\phi=q_\psi=q$. 
Pick a projective isomorphism
$\beta:T_{\xi_1^\phi}\sP^1\to T_{\xi_1^\psi}\sP^1$
carrying the ordered pair
$(\Vw_0^\phi,\Vw_\infty^\phi)$ to
$(\Vw_0^\psi,\Vw_\infty^\psi)$,
and fix $M\in\PGL(2,K)$ such that
\begin{gather*}
M(\xi_1^\phi)=\xi_g,\quad M_*(\Vw_0^\phi)=\widevec{\xi_g0},
\quad\text{and that}\quad M_*(\Vw_\infty^\phi)=\widevec{\xi_g\infty}.
\end{gather*}

Let us see (\textit{\ref{item:r}})$\Leftrightarrow$(\textit{\ref{item:pdata}}) for $j=q$.
Suppose first that $\ph(\phi),\ph(\psi)\in k$. Writing
\begin{gather*}
 M_*D_\phi=\hat{a}_1+\hat{a}_2\quad\text{and}\quad M_*(\beta^*D_\psi)=\hat{b}_1+\hat{b}_2\quad\text{on }\bP^1(k),
\end{gather*}
we have
\begin{align*}
 \ph(\phi)=\ph(\psi)
 &\Leftrightarrow
 \frac{\hat{a}_1}{\hat{a}_2}+\frac{\hat{a}_2}{\hat{a}_1}
 =\frac{\hat{b}_1}{\hat{b}_2}+\frac{\hat{b}_2}{\hat{b}_1}
 \Leftrightarrow
 (\hat{a}_1\hat{b}_2-\hat{b}_1\hat{a}_2)(\hat{a}_1\hat{b}_1-\hat{a}_2\hat{b}_2)=0\\
 &\Leftrightarrow h_*D_\phi=\beta^*D_\psi
 \quad\text{for some }h\in\bT_\phi\\
 &\Leftrightarrow
 (\beta\circ h)\circ\widetilde{(\phi^q)}_{\xi_1^\phi}
 =\widetilde{(\psi^q)}_{\xi_1^\psi}\circ(\beta\circ h)
 \quad\text{for some }h\in\bT_\phi\\
 &\Leftrightarrow
 \text{the intrinsic reduction data of $\phi^q$ and $\psi^q$ are projectively equivalent}.
\end{align*}
The first equivalence is by the projective invariance of cross ratios.
The ``$\Rightarrow$'' implication in the fourth equivalence is by Proposition
\ref{th:divisor}. The ``$\Rightarrow$'' implication in the final equivalence is by
the computation of the intrinsic depths in \eqref{eq:qfin}; conversely, 
if the intrinsic reduction data of $\phi^q$ and $\psi^q$ are 
projectively equivalent via $\alpha:T_{\xi_1^\phi}\sP^1\to T_{\xi_1^\psi}\sP^1$, then
by the computation of the intrinsic depths in \eqref{eq:qfin},
we have $\alpha(\{\Vw_0^\phi,\Vw_\infty^\phi\})
=\{\Vw_0^\psi,\Vw_\infty^\psi\}$, and in turn, 
by the mapping property of the intrinsic reductions in
\eqref{eq:qfin}, $\alpha$ even carries the ordered pair $(\Vw_0^\phi,\Vw_\infty^\phi)$ to $(\Vw_0^\psi,\Vw_\infty^\psi)$. In particular, $\beta^{-1}\circ\alpha\in\bT_\phi$.

Suppose next that
$\ph(\phi)=\ph(\psi)=\infty$. 
By the mapping property of the intrinsic reductions in \eqref{eq:qinf} and Proposition \ref{th:multiple}, both
$\widetilde{(\phi^q)}_{\xi_1^\phi}$ and
$\beta^{-1}\circ\widetilde{(\psi^q)}_{\xi_1^\psi}\circ\beta$
are non-trivial unipotent projective transformations of $T_{\xi_1^\phi}\sP^1$
fixing $\Vw_\infty^\phi$, which are respectively written as
$\zeta\mapsto\zeta+\hat{a}$ and $\zeta\mapsto\zeta+\hat{b}$ with respect to the coordinate induced by $M_*$, where $\hat{a},\hat{b}\in k^\times$. Hence
setting
$h:=M_*^{-1}\circ(\zeta\mapsto(\hat{b}/\hat{a})\zeta)\circ M_*\in\bT_\phi$, 
we have $(\beta\circ h)\circ\widetilde{(\phi^q)}_{\xi_1^\phi}
 =\widetilde{(\psi^q)}_{\xi_1^\psi}\circ(\beta\circ h)$, and then
by the computation of the intrinsic depths in \eqref{eq:qinf}, the intrinsic reduction data of $\phi^q$ and $\psi^q$ are projectively equivalent via $\beta\circ h$.
Conversely, if the intrinsic reduction data of $\phi^q$ and $\psi^q$ are 
projectively equivalent via $\alpha:T_{\xi_1^\phi}\sP^1\to T_{\xi_1^\psi}\sP^1$, then
by the computation of the depths in \eqref{eq:qinf},
$\alpha$ immediately carries the ordered pair $(\Vw_0^\phi,\Vw_\infty^\phi)$ to $(\Vw_0^\psi,\Vw_\infty^\psi)$. 

Finally, in the other cases than the above, neither (\ref{item:r}) nor (\ref{item:pdata}) for $j=q$ is the case. Hence the proof of (\textit{\ref{item:r}})$\Leftrightarrow$(\textit{\ref{item:pdata}}) for $j=q$ and the final assertion for $j=q$ are complete.

We next show (\textit{\ref{item:pdata}}) for $j=q$ $\Leftrightarrow$(\textit{\ref{item:tail}}); the ``$\Leftarrow$'' implication is clear.
Suppose that the intrinsic reduction data of $\phi^q$ and $\psi^q$ are projectively equivalent via a projective isomorphism
$\alpha:T_{\xi_1^\phi}\sP^1\to T_{\xi_1^\psi}\sP^1$.
First, for every $1\le r<q$, by the mapping property and the computation of the intrinsic depths
of the intrinsic reductions of $\phi^r$ and $\psi^r$ in \eqref{eq:short}, 
the same $\alpha$ is a projective equivalence of the
intrinsic reduction data of $\phi^r$ and $\psi^r$.

Let us see by induction on $\ell\ge1$ that the same $\alpha$ is a
projective equivalence of the intrinsic reduction data of $\phi^{\ell q}$ and
$\psi^{\ell q}$. The case of $\ell=1$ is clear by the assumption. 
Pick $\ell\ge1$, and assume the assertion holds for $\ell$. 
Then for every $\Vw\in T_{\xi_1^\phi}\sP^1$, 
we compute not only
\begin{align*}
 \dep_{\alpha(\Vw)}\widetilde{(\psi^{(\ell+1)q})}_{\xi_1^\psi}
 &=\deg(\psi^q)\dep_{\alpha(\Vw)}\widetilde{(\psi^{\ell q})}_{\xi_1^\psi}
 +m_{\alpha(\Vw)}(\psi^{\ell q})
 \dep_{(\psi^{\ell q})_{*,\xi_1^\psi}(\alpha(\Vw))}
 \widetilde{(\psi^q)}_{\xi_1^\psi}\\
 &=\deg(\phi^q)\dep_{\Vw}\widetilde{(\phi^{\ell q})}_{\xi_1^\phi}
 +m_{\Vw}(\phi^{\ell q})
 \dep_{(\phi^{\ell q})_{*,\xi_1^\phi}(\Vw)}
 \widetilde{(\phi^q)}_{\xi_1^\phi}\\
&=\dep_{\Vw}\widetilde{(\phi^{(\ell+1)q})}_{\xi_1^\phi}
\end{align*}
using the composition formula \eqref{eq:depthcomp} 
applied to $\psi^{(\ell+1)q}=\psi^q\circ\psi^{\ell q}$ and
$\phi^{(\ell+1)q}=\phi^q\circ\phi^{\ell q}$ 
and the chain rule for the tangent maps
(and recalling Fact \ref{th:argument}), but also
\begin{align*} 
\alpha\circ\widetilde{(\phi^{(\ell+1)q})}_{\xi_1^\phi}
 &=\alpha\circ(\phi^q)_{*,\xi_1^\phi}
 \circ(\phi^{\ell q})_{*,\xi_1^\phi}
=\alpha\circ\widetilde{(\phi^q)}_{\xi_1^\phi}
 \circ\widetilde{(\phi^{\ell q})}_{\xi_1^\phi}
 =\widetilde{(\psi^q)}_{\xi_1^\psi}\circ\alpha
 \circ\widetilde{(\phi^{\ell q})}_{\xi_1^\phi}\\
 &=\widetilde{(\psi^q)}_{\xi_1^\psi} \circ\widetilde{(\psi^{\ell q})}_{\xi_1^\psi}\circ\alpha=(\psi^q)_{*,\xi_1^\psi}
 \circ(\psi^{\ell q})_{*,\xi_1^\psi}\circ\alpha
=\widetilde{(\psi^{(\ell+1)q})}_{\xi_1^\psi}\circ\alpha.
\end{align*}
using again the chain rule for the tangent maps. 
Hence the assertion for $\ell+1$ also holds. We also note that 
the same $\alpha$ is a
projective equivalence of the intrinsic reduction data of $\phi^{0q}=\Id$ and
$\psi^{0q}=\Id$.

Finally, once this implication
(\textit{\ref{item:pdata}}) for $j=q$ 
$\Rightarrow$(\textit{\ref{item:tail}}) 
for every $j\ge 0$ divisible by $q$ is at our disposal,
for every $j=\ell q+r$, $\ell\ge1$ and $1\le r<q$, 
we have
$\dep_{\alpha(\Vw)}\widetilde{(\psi^j)}_{\xi_1^\psi}
 =\dep_{\Vw}\widetilde{(\phi^j)}_{\xi_1^\phi}$
for every $\Vw\in T_{\xi_1^\phi}\sP^1$
by a computation similar to that in the above induction 
using the composition formula \eqref{eq:depthcomp} 
applied to $\psi^j=\psi^r\circ\psi^{\ell q}$ and $\phi^j=\phi^r\circ\phi^{\ell q}$
and the chain rule for the tangent maps,
and
\begin{align*}
 \alpha\circ\widetilde{(\phi^j)}_{\xi_1^\phi}
 \equiv\alpha(\Vw_\infty^\phi)=\Vw_\infty^\psi
 \equiv\widetilde{(\psi^j)}_{\xi_1^\psi}\circ\alpha
\end{align*}
by the mapping property of the intrinsic reductions when $q\nmid j$ in \eqref{eq:short}. Hence (\textit{\ref{item:pdata}}) for $j=q$ $\Leftrightarrow$ (\textit{\ref{item:tail}}).

Finally let us see (\textit{\ref{item:pdata}}) for some $j>q$
$\Rightarrow$(\textit{\ref{item:r}}).
Pick any such $j>q$ and suppose that the intrinsic reduction data of $\phi^j$
and $\psi^j$ are projectively equivalent via
$\alpha:T_{\xi_1^\phi}\sP^1\to T_{\xi_1^\psi}\sP^1$.
As seen in (\ref{item:nonend}) and (\ref{item:boundary}), 
\begin{gather*}
 \sup_{\Vw\in T_{\xi_1}\sP^1}
 \dep_{\Vw}\widetilde{(\phi^j)}_{\xi_1}
 \le\frac{\deg(\phi^j)}{2}\quad\text{and the equality holds iff } \ph(\phi)=\infty,
\end{gather*}
and the above supremum is preserved under the projective equivalence $\alpha$. Hence we are done when
$\ph(\phi)=\ph(\psi)=\infty$.
Suppose  that $\ph(\phi)\in k$, 
and we claim that
\begin{gather*}
 \left\{\Vw\in T_{\xi_1}\sP^1:
 \dep_{\Vw}\widetilde{(\phi^j)}_{\xi_1}
 =
 \max_{\Vw'\in T_{\xi_1}\sP^1}
 \dep_{\Vw'}\widetilde{(\phi^j)}_{\xi_1}\right\}
=\{\Vw_0,\Vw_\infty\}\quad\text{and}\\
 \biggl\{\Vw\in T_{\xi_1}\sP^1\setminus\{\Vw_0,\Vw_\infty\}:
 \dep_{\Vw}\widetilde{(\phi^j)}_{\xi_1}
 = \max_{\Vw'\in T_{\xi_1}\sP^1\setminus\{\Vw_0,\Vw_\infty\}}
 \dep_{\Vw'}\widetilde{(\phi^j)}_{\xi_1}\biggr\}=\supp(D_\phi),
\end{gather*}
so in particular $\{\Vw_0,\Vw_\infty\}$ and $\supp(D_\phi)$ are determined by 
the data, and the same holds for $\psi$; noting that
$(\phi^q)_{*,\xi_1}\bigl(\{\Vw_0,\Vw_\infty\}\bigr)
 \subset\{\Vw_0,\Vw_\infty\}$, set
\begin{gather*}
 N_{\Vw}=N_{\Vw}^{(j)}:= 
 \min\bigl\{\ell\ge 0: \widetilde{(\phi^{q\ell})}_{\xi_1}\Vw\in\{\Vw_0,\Vw_\infty\}\bigr\}
\in\{0,1,\ldots\}\sqcup\{+\infty\}
\quad\text{for every }\Vw\in T_{\xi_1}\sP^1, 
\end{gather*}
so that $N_{\Vw}=N_{\widetilde{(\phi^q)}_{\xi_1}\Vw}+1$ if $N_{\Vw}>0$
and that $\{\Vw:N_{\Vw}=0\}=\{\Vw_0,\Vw_\infty\}$ and 
$\{\Vw:N_{\Vw}=1\}=\supp(D_\phi)$ by \eqref{eq:qfin}, 
under the usual convention $\min\emptyset=+\infty$ and $+\infty+1=+\infty$.
Pick $\Vw\in T_{\xi_1}\sP^1$. By \eqref{eq:short},
if $N_{\Vw}>\lfloor j/q\rfloor$, 
$\dep_{\Vw}\widetilde{(\phi^j)}_{\xi_1}=0$. 
If $N_{\Vw}\le\lfloor j/q\rfloor$, then
by the composition formula
\eqref{eq:depthcomp}, \eqref{eq:qfin}, the choice of $N_{\Vw}$,
and the non-archimedean argument principle \eqref{eq:FKN},
we compute
\begin{multline}
 \dep_{\Vw}\widetilde{(\phi^j)}_{\xi_1}
 =\deg(\phi^{j-N_{\Vw}q})\dep_{\Vw}\widetilde{(\phi^{N_{\Vw}q})}_{\xi_1}+
m_{\Vw}(\phi^{N_{\Vw}q})
 \dep_{\widetilde{(\phi^{N_{\Vw}q})}_{\xi_1}\Vw}
 \widetilde{(\phi^{j-N_{\Vw}q})}_{\xi_1}\\
 =0+m_{\Vw}(\phi^{N_{\Vw}q})
\dep_{\Vu}\widetilde{(\phi^{j-N_{\Vw}q})}_{\xi_1},\label{eq:descend}
\end{multline}
where $\Vu$ is either $\Vw_0$ or $\Vw_\infty$, and
no matter which of $\Vw_0$ or $\Vw_\infty$ the $\Vu$ is, 
for every $1\le n\le \lfloor j/q\rfloor$,
we similarly and recursively compute
\begin{align*}
 \dep_{\Vu}\widetilde{(\phi^{j-(n-1)q})}_{\xi_1}
 &=
 \deg(\phi^{j-nq})\dep_{\Vu}\widetilde{(\phi^q)}_{\xi_1}
 +
 m_{\Vu}(\phi^q)\dep_{\widetilde{(\phi^q)}_{\xi_1}\Vu}
 \widetilde{(\phi^{j-nq})}_{\xi_1}\\
 &=
 \deg(\phi^{j-nq})\frac{\deg(\phi^q)-2}{2}
 +1\cdot\dep_{\Vw_\infty}\widetilde{(\phi^{j-nq})}_{\xi_1}\\
&=\frac{\deg(\phi^q)-2}{2}\sum_{\ell=n}^{\lfloor j/q \rfloor}
\deg(\phi^{j-\ell q})+\dep_{\Vu}\widetilde{(\phi^{j-\lfloor j/q \rfloor q})}_{\xi_1}
\end{align*}
so that
\begin{multline*}
\dep_{\Vu}\widetilde{\bigl(\phi^{j-(n-1)q}\bigr)}_{\xi_1}
-(\deg\phi)\dep_{\Vu}\widetilde{\bigl(\phi^{j-nq}\bigr)}_{\xi_1}\\
 \biggl(= \deg(\phi^{j-nq})\frac{\deg(\phi^q)-2}{2}
 -\dep_{\Vw_\infty}\widetilde{\bigl(\phi^{j-nq}\bigr)}_{\xi_1}\biggr)\\
 =\frac{\deg(\phi^q)-2}{2}\biggl(\deg(\phi^{j-nq})-\sum_{\ell={n+1}}^{\lfloor j/q \rfloor}
\deg(\phi^{j-\ell q})\biggr)
-\dep_{\Vw_\infty}\widetilde{(\phi^{j-\lfloor j/q \rfloor q})}_{\xi_1}\\
\ge\deg(\phi^{j-\lfloor j/q \rfloor q})\biggl(\frac{2^q-2}{2}
\cdot\frac{(2^q-2)2^{(\lfloor j/q \rfloor-n)q}+1}{2^q-1}-\frac{1}{2}\biggr)
\ge\deg(\phi^{j-\lfloor j/q \rfloor q})\cdot
\frac{2^q-3}{2}>0,
\end{multline*}
recalling \eqref{eq:short} and $q>1$. Consequently, if $0<N_{\Vw}\le\lfloor j/q\rfloor$, we compute
\begin{align*} 
 \dep_{\Vw}\widetilde{(\phi^j)}_{\xi_1} 
&=m_{\Vw}(\phi^q)\cdot m_{\widetilde{(\phi^q)}_{\xi_1}\Vw}(\phi^{N_{\widetilde{(\phi^q)}_{\xi_1}\Vw}q})
 \cdot\dep_{\Vu}\widetilde{(\phi^{j-N_{\Vw}q})}_{\xi_1}\\
&\le m_{\widetilde{(\phi^q)}_{\xi_1}\Vw}(\phi^{N_{\widetilde{(\phi^q)}_{\xi_1}\Vw}q})\cdot
(\deg\phi)\dep_{\Vu}\widetilde{(\phi^{j-N_{\Vw}q})}_{\xi_1}\\
&<m_{\widetilde{(\phi^q)}_{\xi_1}\Vw}\bigl(\phi^{N_{\widetilde{(\phi^q)}_{\xi_1}\Vw}q}\bigr)\cdot
\dep_{\Vu}\widetilde{\bigl(\phi^{j-N_{\widetilde{(\phi^q)}_{\xi_1}\Vw}q}\bigr)}_{\xi_1} 
=\dep_{\widetilde{(\phi^q)}_{\xi_1}\Vw}\widetilde{(\phi^j)}_{\xi_1}.  
\end{align*}  
Hence the claim holds.
Since $\deg D_\phi=\deg\phi=2$, the divisor $D_\phi$ itself is determined by
$\supp D_\phi=\{\Vw_1,\Vw_2\}$, and by \eqref{eq:qfin}, we have 
$\widetilde{(\phi^j)}_{\xi_1}(\{\Vw_0,\Vw_\infty\})=\{\Vw_\infty\}$ if $q\mid j$;
if $q\nmid j$ or equivalently $\widetilde{(\phi^j)}_{\xi_1}\equiv\Vw_\infty$, then
$\widetilde{(\phi^j)}_{\xi_1}(\{\Vw_0,\Vw_\infty\})=\{\Vw_\infty\}$ is clear.
Hence the ordered pair $(\Vw_0,\Vw_\infty)$ is still determined by
the intrinsic reduction data of $\phi^j$, and the same holds for $\psi$.
Consequently, the projective equivalence $\alpha$ carries
$(\Vw_0^\phi,\Vw_\infty^\phi)$ and $D_\phi$ to
$(\Vw_0^\psi,\Vw_\infty^\psi)$ and $D_\psi$, respectively.
In particular, by Proposition \ref{th:divisor}, we have
$\ph(\phi)=\ph(\psi)$.
Now the proof of (\textit{\ref{item:pdata}}) for some $j\ge q$
$\Rightarrow$(\textit{\ref{item:r}}) is complete.

It remains only to prove the final assertion under
\begin{gather*}
 (\text{\ref{item:pdata}) for some }j>q\quad\text{and}\quad\ph(\phi)=\ph(\psi)=\infty.
\end{gather*}
Pick any such integer $j>q$ that the intrinsic reduction data of
$\phi^j$ and $\psi^j$ are projectively equivalent via
$\alpha:T_{\xi_1^\phi}\sP^1\to T_{\xi_1^\psi}\sP^1$. 
If $q\nmid j$ or equivalently both 
$\widetilde{(\phi^j)}_{\xi_1}$ and $\widetilde{(\psi^j)}_{\xi_1}$
are $\equiv\Vw_\infty$, we first have $\alpha(\Vw_\infty^\phi)=\Vw_\infty^\psi$.
Recall \eqref{eq:nonenddepth} and \eqref{eq:qinf}, and $q>1$. Then we have
\begin{gather*}
 \biggl\{\Vw\in(T_{\xi_1}\sP^1)\setminus\{\Vw_\infty\}:
 (\deg(\phi^j)\ge)\dep_{\Vw}\widetilde{(\phi^j)}_{\xi_1}
 \ge\frac{\deg(\phi^j)}{2}(=\dep_{\Vw_0}\widetilde{(\phi^j)}_{\xi_1})\biggr\}=\{\Vw_0\}
\end{gather*}
noting that
$\dep_{\Vw_\infty}\widetilde{(\phi^j)}_{\xi_1}
>\deg(\phi^{j-q})\cdot 0+m_{\Vw_\infty}(\phi^q)\dep_{\Vw_\infty}(\widetilde{(\phi^{j-q})}_{\xi_1})\ge 0$ 
by the composition formula \eqref{eq:depthcomp} applied to 
$\phi^j=\phi^{j-q}\circ\phi^q$. Hence $\alpha$
carries the ordered pair $(\Vw_0^\phi,\Vw_\infty^\phi)$ to 
$(\Vw_0^\psi,\Vw_\infty^\psi)$.

Suppose now that $q\mid j$ or equivalently $\lfloor j/q\rfloor=j/q$. We claim that
\begin{gather*}
 \dep_{\Vw_0}\widetilde{(\phi^j)}_{\xi_1}
 >\dep_{\Vw_\infty}\widetilde{(\phi^j)}_{\xi_1}
 >\sup_{\Vw\in(T_{\xi_1}\sP^1)\setminus\{\Vw_0,\Vw_\infty\}}
\dep_{\Vw}\widetilde{(\phi^j)}_{\xi_1},
\end{gather*}
and then $\alpha$ would still carry the ordered pair $(\Vw_0^\phi,\Vw_\infty^\phi)$ to
$(\Vw_0^\psi,\Vw_\infty^\psi)$;
for, pick $\Vw\in(T_{\xi_1}\sP^1)\setminus\{\Vw_0,\Vw_\infty\}$.
Setting
\begin{gather*} 
  N_{\Vw}=N_{\Vw}^{(j)}:= \min\bigl\{\ell\ge 0:
 \widetilde{(\phi^{\ell q})}_{\xi_1}(\Vw)
 \in\{\Vw_0,\Vw_\infty\}\bigr\}>0
\end{gather*} 
under the convention $\min_\emptyset:=+\infty$,
by \eqref{eq:short},
if $N_{\Vw}>\lfloor j/q\rfloor$, then
$\dep_{\Vw}\widetilde{(\phi^j)}_{\xi_1}=0$. 
Alternatively if $N_{\Vw}\le\lfloor j/q\rfloor$, then by \eqref{eq:qinf},
we have $\widetilde{(\phi^{N_{\Vw}q})}_{\xi_1}(\Vw)=\Vw_0$ and,
by \eqref{eq:short} and a computation similar to \eqref{eq:descend},
\begin{align*}
 \dep_{\Vw}\widetilde{(\phi^j)}_{\xi_1}
 =
 0+1\cdot\dep_{\Vw_0}
 \widetilde{(\phi^{j-N_{\Vw}q})}_{\xi_1}
 =
 \begin{cases}
 0  &\text{if }N_{\Vw}=\lfloor j/q\rfloor\\
 \displaystyle\frac{\deg(\phi^{j-N_{\Vw}q})}{2} &\text{if }N_{\Vw}<\lfloor j/q\rfloor
 \end{cases}
\le\frac{\deg(\phi^{j-q})}{2}.
\end{align*}
On the other hand, by \eqref{eq:nonenddepth} and 
the composition formula \eqref{eq:depthcomp}
applied to $\phi^j=\phi^{j-q}\circ\phi^q$ and \eqref{eq:qinf}, 
we have
\begin{align*}
 \dep_{\Vw_0}\widetilde{(\phi^j)}_{\xi_1}>\dep_{\Vw_\infty}\widetilde{(\phi^j)}_{\xi_1}
 &=\deg(\phi^{j-q})\dep_{\Vw_\infty}\widetilde{(\phi^q)}_{\xi_1}
 +m_{\Vw_\infty}(\phi^q)
 \dep_{\widetilde{(\phi^q)}_{\xi_1}\Vw_\infty}\widetilde{(\phi^{j-q})}_{\xi_1}\\
 &\ge\deg(\phi^{j-q})\frac{\deg(\phi^q)-2}{2}
 >\frac{\deg(\phi^{j-q})}{2}.
\end{align*}
Hence the claim holds. Now the proof of Theorem \ref{th:residuemoduli} is complete. \qed

\section{Proof of Theorem \ref{th:arc}}\label{sec:arc}
For the standard constructions on dynamical moduli recalled below, see
\cite{Silverman98}, \cite[Chapter~4]{SilvermanDynamics}, \cite[\S2]{DeMarco05},
and \cite[Chapter~1]{MFK94}.

\begin{fact}\label{th:git}
We denote an arbitrary field by $k$. Fix any integer $d\ge2$.

\begin{inparaenum}[(i)]
\item\label{item:git} {\em GIT construction and quotients.}
Put $W_{\bZ}:=\bZ^2$ and
\begin{gather*}
 V_{d,\bZ}:=W_{\bZ}\otimes_{\bZ}
 \Sym^d(W_{\bZ}^{\vee})
 \quad\text{and}\quad
 \bP(V_{d,\bZ})\cong\bP^{2d+1}_{\bZ},
\end{gather*}
and let $A_{d,\bZ}:=\Sym_{\bZ}
(V_{d,\bZ}^{\vee})$.  The conjugation action of
$\SL_{2,\bZ}$ on degree-$d$ rational selfmaps of $\bP^1$
induces a linear action on $V_{d,\bZ}$.
The resultant $\Res$ is a homogeneous
$\SL_{2,\bZ}$-invariant, and
\begin{gather*}
 \Rat_d:=D_+(\Res)
 \subset\bP(V_{d,\bZ}).
\end{gather*}
Its geometric quotient
\begin{gather*}
 \rM_d:=\Rat_d/\SL_{2,\bZ}
\end{gather*}
is the moduli scheme of degree-$d$ rational selfmaps of $\bP^1$.  By
Silverman's GIT construction
\cite[Theorems~1.3 and~2.1]{Silverman98}, the projective GIT compactification
of $\rM_d$ is
\begin{gather*}
 \overline{\rM}_d:=
 \operatorname{Proj}\bigl(A_{d,\bZ}^{\SL_{2,\bZ}}\bigr).
\end{gather*}
For every
algebraically closed field $\Omega$, the map
$\SL_2(\Omega)\to\PGL(2,\Omega)$ is surjective, so the
$\Omega$-points of $\rM_d$ are the $\PGL(2,\Omega)$-conjugacy classes; see
\cite[Theorem~3.2]{Silverman98}.

Over $k$, put $V_{d,k}:=V_{d,\bZ}\otimes_{\bZ}k$ and denote by
\begin{gather*}
 \Rat_{d,k}:=(\Rat_d)_k,
 \quad
 \overline{\rM}_{d,k}:=(\overline{\rM}_d)_k,
 \quad\text{and}\quad
 \rM_{d,k}:=(\rM_d)_k
\end{gather*}
the base changes to $k$ of $\Rat_d$, $\overline{\rM}_d$, and
$\rM_d$, respectively.  Then $\Rat_{d,k}
=D_+(\Res_k)\subset\bP(V_{d,k})$, and
$\overline{\rM}_{d,k}$ is projective over $k$, with $\rM_{d,k}$ as an open
dense subscheme.  Let $\bP(V_{d,k})^{\mathrm{ss}}$ and
$\bP(V_{d,k})^{\mathrm{s}}$ denote the semistable and stable loci,
respectively, with respect to the induced $\SL_{2,k}$-linearization
of $\cO_{\bP(V_{d,k})}(1)$.  By the relative GIT construction over
$\bZ$ in
\cite[\S2, especially Theorem~2.1 and its proof]{Silverman98}, which invokes Seshadri's
geometric reductivity theorem over arbitrary bases, the quotient morphism
$q_{d,k}:\bP(V_{d,k})^{\mathrm{ss}}\to\overline{\rM}_{d,k}$
is categorical, and its restriction to $\bP(V_{d,k})^{\mathrm{s}}$ is geometric.
Thus, for every algebraically closed extension $\Omega/k$, two stable
$\Omega$-points have the same image under $q_{d,k}$ if and only if they belong to
the same $\PGL(2,\Omega)$-orbit.
Moreover, the quotient morphism on the semistable locus is compatible with
specialization: if $R$ is a valuation ring with fraction field $F$ and residue
field $k$, and a coefficient point $f\in\bP(V_{d,F})(F)$ extends to an
$R$-point with semistable closed fiber $\widehat{f}$, then
\begin{gather*}
 q_{d,k}(\widehat{f})=\spe_R\bigl(q_{d,F}(f)\bigr).
\end{gather*}
We call a point of
$\overline{\rM}_{d,k}(k)\setminus\rM_{d,k}(k)$ a
{\em boundary point of $\rM_{d,k}$}.

Over $\bZ$, the elementary symmetric functions of the $d+1$
fixed-point multipliers, counted with multiplicity, define global regular
functions
\begin{gather*}
 \sigma_1,\ldots,\sigma_{d+1}
 \in\Gamma(\rM_d,\cO_{\rM_d})
\end{gather*}
by \cite[Proposition~4.2, p.~56]{Silverman98}.
In the quadratic $d=2$ case, the Milnor--Silverman fixed-point multiplier
coordinates identify $\rM_2\cong\bA^2_{\bZ}$ and
$\overline{\rM}_2\cong\bP^2_{\bZ}$, with
$\rM_2\cong\{s_3\ne0\}$ and affine coordinates $(s_1/s_3,s_2/s_3)$
equal to $(\sigma_1,\sigma_2)$
(\cite[Theorem~1.2, Corollary~1.4, and Theorem~1.5]{Silverman98}).
Since $\sigma_3=\sigma_1-2$
(\cite[Remark before Lemma~5.3, p.~63]{Silverman98}), the fixed-point
multiplier polynomial over $\bZ$ is
\begin{gather*}
 T^3-\sigma_1T^2+\sigma_2T-(\sigma_1-2)
 \in\bZ[\sigma_1,\sigma_2][T].
\end{gather*}
Homogenizing this polynomial in the Milnor--Silverman coordinates over
$\bZ$ and then changing base to $k$, we associate to each
$y=[s_1:s_2:s_3]\in\overline{\rM}_{2,k}(k)$ the fixed-point multiplier polynomial
\begin{gather}
 P_y(T):=s_3T^3-s_1T^2+s_2T-(s_1-2s_3).
 \label{eq:multiplierpoly}
\end{gather}
Changing the homogeneous representative of $y$ multiplies $P_y$ by a
nonzero element of $k$. When $y\in\rM_{2,k}(k)$, the polynomial
$P_y(T)/s_3$ is monic and has the three fixed-point multipliers as its roots,
counted with multiplicity.
For a boundary point $y=[s_1:s_2:0]$ and $\omega\in k^\times$, if
$P_y(\omega)=0$, then $s_2\omega=s_1(\omega^2+1)$, which forces $s_1\ne0$
and $s_2=s_1(\omega+\omega^{-1})$, so that
\begin{gather*}
 P_y(T)=-s_1T^2+s_2T-s_1
 =-s_1(T-\omega)(T-\omega^{-1})
 \quad\text{and}\quad
 y=[1:\omega+\omega^{-1}:0].
\end{gather*}

\item\label{item:graphs} {\em Iteration and graph closures.}
For every $j\ge2$, iteration defines a morphism
$\rM_d\to\rM_{d^j}$ and hence a rational map
\begin{gather*}
 \Phi_{d,j}:\overline{\rM}_d
 \dashrightarrow\overline{\rM}_{d^j}\quad\text{and}\quad
 \Phi_{d,j}([\psi])=[\psi^j].
\end{gather*}
Over $k$, denote by
\begin{gather*}
 \Phi_{d,j,k}:\overline{\rM}_{d,k}\dashrightarrow\overline{\rM}_{d^j,k}
\end{gather*}
the base change of $\Phi_{d,j}$; thus
$\Phi_{d,j,k}([\psi])=[\psi^j]$ for every $[\psi]\in\rM_{d,k}$.

For every $n\ge2$, put
\begin{gather*}
 \mathcal{G}_{d,n}
 :=
 \overline{
 (\Id,\Phi_{d,2},\ldots,\Phi_{d,n})(\rM_d)
 }^{\,\Sch}
 \subset
 \overline{\rM}_d\times_{\bZ}
 \prod_{j=2}^n\overline{\rM}_{d^j},
\end{gather*}
where $\overline{(\,\cdot\,)}^{\,\Sch}$ denotes
scheme-theoretic closure.  Over $k$, put
\begin{gather*}
 \mathcal{G}_{d,n,k}
 :=
 \overline{
 (\Id,\Phi_{d,2,k},\ldots,\Phi_{d,n,k})(\rM_{d,k})
 }^{\,\Sch}
 \subset
 \overline{\rM}_{d,k}\times_k\prod_{j=2}^n\overline{\rM}_{d^j,k}.
\end{gather*}
When $\Char k=0$, we have
$\mathcal{G}_{d,n,k}=(\mathcal{G}_{d,n})_k$.
For every field $k$, the scheme $\mathcal{G}_{d,n,k}$ is integral and
projective over $k$.
The schemes $(\mathcal{G}_{d,n,k})_{n\ge2}$ form an inverse system.
More precisely, for $2\le m\le n$, coordinate truncation of tails induces
\begin{gather*}
 \rho_{n,m,k}:\mathcal{G}_{d,n,k}\to\mathcal{G}_{d,m,k},\quad
 (x_1,x_2,\ldots,x_m,\ldots,x_n)\mapsto(x_1,x_2,\ldots,x_m).
\end{gather*}

\item\label{item:arcs} {\em Algebraic arcs.}
Put $R:=k[[t]]$, $F:=k((t))$, and $\eta_R:=\Spec F$.  Let
$\gamma\in\overline{\rM}_{d,k}(R)$ be an algebraic arc whose generic fiber lies
in $\rM_{d,k}$.  For every $j\ge2$, the rational section
\begin{gather*}
 \Phi_{d,j,k}\circ(\gamma|_{\eta_R}):\eta_R\to\overline{\rM}_{d^j,k}
\end{gather*}
has a unique extension to an $R$-point by properness of
$\overline{\rM}_{d^j,k}$; denote the value of this extension at the unique
closed point of $\Spec R$ by
\begin{gather*}
 \spe_R\bigl(\Phi_{d,j,k}\circ(\gamma|_{\eta_R})\bigr)
 \in\overline{\rM}_{d^j,k}(k).
\end{gather*}
For every $n\ge2$, these unique extensions together with $\gamma$ define a
morphism
\begin{gather*}
 \Spec R\to
 \overline{\rM}_{d,k}\times_k\prod_{j=2}^n\overline{\rM}_{d^j,k}
\end{gather*}
which factors through $\mathcal{G}_{d,n,k}$.
\end{inparaenum}
\end{fact}

\begin{remark}
When $k=\bC$, the schemes $\mathcal{G}_{2,n,\bC}$ appearing in Subsection
\ref{sec:moduli} are the finite-level graph closures whose inverse
limit is the iterate compactification considered by Kiwi and Nie;
see \cite[\S1]{KiwiNie26}.
\end{remark}

\subsection{Proof of Theorem \ref{th:arc}(\textit{\ref{item:intrinsic}})}\label{sec:proofintrinsic}

Let $k$, $\bK$, and $R$ be as in Theorem \ref{th:arc}. Pick an order $p>1$
boundary point $x\in\overline{\rM}_{2,k}(k)$ and choose a primitive $p$th
root of unity $\omega\in k^\times$ such that $P_x(\omega)=0$.

Pick any antispecialization arc $\gamma\in\overline{\rM}_{2,k}(R)$ degenerating to $x$. 
We first express the intrinsic antispecialization parameter
$\as(\gamma)$ in terms of the fixed-point multipliers and prove
the first assertion in (\textit{\ref{item:intrinsic}}). After a finite ramified extension $R\subset R'=k[[u]]$ such that $F'=k((u))$, put
\begin{gather*}
 \gamma_{R'}:=\gamma\times_R R'
 \quad\text{and}\quad
 \eta_{R'}:=\Spec F'.
\end{gather*}
By the factorization of $P_x$ in Definition \ref{th:orderboundary}, after
enlarging $R'$ if necessary, we can choose two roots $\mu_1$ and $\mu_2$ in
$R'$, counted with multiplicity, of the fixed-point multiplier polynomial
\eqref{eq:multiplierpoly} of $\gamma_{R'}|_{\eta_{R'}}$, labeled so that
\begin{gather*}
 \widehat{\mu_1}=\omega\quad\text{and}\quad
 \widehat{\mu_2}=\omega^{-1}.
\end{gather*}
After enlarging $R'$ if necessary, let $\mu_3\in F'$ be the remaining root
of the fixed-point multiplier polynomial of the generic fiber
$\gamma_{R'}|_{\eta_{R'}}\in\rM_{2,k}(F')$. Noting that
$\mu_1+\mu_2+\mu_3=\mu_1\mu_2\mu_3+2$ by \eqref{eq:multiplierpoly}, we have
\begin{gather*}
 \widehat{2-\mu_1-\mu_2}
 =2-\omega-\omega^{-1}
 =-\frac{(\omega-1)^2}{\omega}\ne0,
 \quad\text{and on the other hand},\quad
 \widehat{\mu_1^p-1}=\omega^p-1=0.
\end{gather*}
Hence setting
\begin{align*}
 \varepsilon&:=1-\mu_1\mu_2
 \in\mathfrak{m}_{R'}\setminus\{0\}
 \subset\bK^{\circ\circ}\setminus\{0\} \quad\text{and}\\
 \delta&:=\mu_1^p-1
 \in\mathfrak{m}_{R'}\subset\bK^{\circ\circ},
\end{align*}
we have
\begin{gather*}
  \mu_3=\frac{2-\mu_1-\mu_2}{1-\mu_1\mu_2}
 \quad\text{and}\quad
 |\mu_3|=\frac{1+0}{|1-\mu_1\mu_2|}>1.
\end{gather*}
By a direct computation, the rational function
\begin{gather}
 \phi_{F'}(z):=
 z\frac{(1-\mu_1)z+\mu_1(1-\mu_2)}
 {\mu_2(1-\mu_1)z+(1-\mu_2)}
 \in R'(z)\subset F'(z)\label{eq:MSnormal}
\end{gather}
is of degree $2$ and fixes $0,1,\infty$, with
\begin{gather*}
 \phi_{F'}'(z)
 =\frac{\mu_2(1-\mu_1)^2z^2+2(1-\mu_1)(1-\mu_2)z
 +\mu_1(1-\mu_2)^2}
 {\bigl(\mu_2(1-\mu_1)z+(1-\mu_2)\bigr)^2},
\end{gather*}
and the corresponding multipliers are computed as
\begin{gather*}
 \phi_{F'}'(0)=\mu_1,
 \quad
 \frac{\rd}{\rd z}
 \left(\frac{1}{\phi_{F'}(1/z)}\right)\biggm|_{z=0}
 =\mu_2,
 \quad\text{and},\\
 \phi_{F'}'(1)
 =1+\frac{(1-\mu_1)(1-\mu_2)}{1-\mu_1\mu_2}
 =\frac{2-\mu_1-\mu_2}{1-\mu_1\mu_2}
 =\mu_3,
\end{gather*}
so $\phi_{F'}$ is a Puiseux lift of $\gamma$
by the Milnor--Silverman fixed-point
multiplier isomorphism (see Fact \ref{th:git}(\textit{\ref{item:git}})).
By the normal form \eqref{eq:MSnormal}, we have 
$\widetilde{\phi_{F'}}_{\xi_g}(\zeta)=\omega\zeta$ in $k(\zeta)$, so that
$(\phi_{F'}^*\delta_{\xi_g})(\{\xi_g\})=1$ (see Fact \ref{th:algred}). 
We also note that the two simple poles of $\phi_{F'}$ are $\infty$ and
\begin{gather*}
 b:=\frac{1-\mu_2}{\mu_2(\mu_1-1)},
 \quad\text{and note that}\quad
 \widehat{b}=\frac{1-\omega^{-1}}{\omega^{-1}(\omega-1)}=1,
\end{gather*}
and that $(\phi_{F'})_*(\widevec{\xi_g1})=\widevec{\xi_g\omega}\neq\widevec{\xi_g1},\widevec{\xi_g\infty}$.
Hence by the non-archimedean argument principle \eqref{eq:FKN}, we compute
\begin{gather*}
 1=(\phi_{F'}^*\delta_\infty)\bigl(U(\widevec{\xi_g1})\bigr)
 =s_{\widevec{\xi_g1}}(\phi_{F'})+0
 =\dep_{\widevec{\xi_g1}}\widetilde{\phi_{F'}}_{\xi_g},
\end{gather*}
and in turn $\sum_{\Vv\in(T_{\xi_g}\sP^1)\setminus\{\widevec{\xi_g1}\}}\dep_{\Vv}\widetilde{\phi_{F'}}_{\xi_g}
=2-\dep_{\widevec{\xi_g1}}\widetilde{\phi_{F'}}_{\xi_g}-(\phi_{F'}^*\delta_{\xi_g})(\{\xi_g\})=0$. Hence by the slope formula \eqref{eq:slope} together with the convexity of
$\hypRes_{\phi_{F'}}$, we have $\xi_{\phi_{F'}}=\xi_g$.
Now for every Puiseux lift $\phi\in\bK(z)$ of $\gamma$, there is $M\in\PGL(2,\bK)$ such that
\begin{gather*}
 M\circ\phi\circ M^{-1}=\phi_{F'}
 \quad\text{and}\quad
 M(\xi_\phi)=\xi_g.
\end{gather*}
By the construction over $\bZ$ in Fact \ref{th:git}, the fixed-point
multiplier polynomial $P_\gamma(T)\in R[T]$ satisfies
$P_\gamma(T)=s_3\prod_{i=1}^3(T-\mu_i)$ after base change to $F'$, so
$\mu_1+\mu_2+\mu_3=s_1/s_3$.
Using $(1-\mu_1\mu_2)\mu_3=2-\mu_1-\mu_2$ and
$\widehat{2-\mu_1-\mu_2}=2-\omega-\omega^{-1}$, together with the compatibility
of specialization with the base change $R\subset R'$, we compute
\begin{gather}
\spe_{R'}\biggl(
 -\frac{(\mu_1^p-1)(\mu_2^p-1)}{1-\mu_1\mu_2}\biggr)
=\spe_{R'}\biggl(
 -\frac{(\mu_1^p-1)(\mu_2^p-1)}{1-\mu_1\mu_2}\biggr)\cdot 1^3
=\spe_{R'}\Biggl(-\frac{(\mu_1^p-1)(\mu_2^p-1)}{1-\mu_1\mu_2}\cdot\label{eq:asmultiplier}\\
 \frac{p\omega^{p-1}\cdot p\omega^{-(p-1)}}{ \prod_{\substack{\zeta^p=1 \\ \zeta\ne\omega}}(\mu_1-\zeta) \prod_{\substack{\zeta^p=1 \\ \zeta\ne\omega^{-1}}}(\mu_2-\zeta)}\cdot
 \frac{\prod_{\substack{\Phi_p(\zeta)=0 \\ \zeta\ne\omega}}(\mu_1-\zeta)
\prod_{\substack{\Phi_p(\zeta)=0 \\ \zeta\ne\omega^{-1}}}(\mu_2-\zeta)}
 {\Phi_p'(\omega)\Phi_p'(\omega^{-1})}\cdot
 \frac{\prod_{\Phi_p(\zeta)=0}(\mu_3-\zeta)}
 {\mu_3\cdot(\mu_1+\mu_2+\mu_3)^{\deg\Phi_p-1}}
 \Biggr)\notag\\
 =\spe_{R'}\left( 
 -\frac{\Phi_p(\mu_1)\Phi_p(\mu_2)\Phi_p(\mu_3)p^2}
 {(1-\mu_1\mu_2)\Phi_p'(\omega)\Phi_p'(\omega^{-1})}
\frac{1}
 {\mu_3(s_1/s_3)^{\deg\Phi_p-1}}
 \right)\notag\\
=\spe_{R'}\biggl(
 -\frac{p^2\Res_T(P_\gamma,\Phi_p)}
 {(1-\mu_1\mu_2)\mu_3\Phi_p'(\omega)\Phi_p'(\omega^{-1})s_1^{\deg\Phi_p-1}s_3}
 \biggr)\notag\\
=\spe_{R'}\biggl(
 -\frac{p^2\Res_T(P_\gamma,\Phi_p)}
 {(2-\omega-\omega^{-1})\Phi_p'(\omega)\Phi_p'(\omega^{-1})s_1^{\deg\Phi_p-1}s_3}
 \biggr)=\as(\gamma)\notag. 
\end{gather}

\begin{remark}
Here, for later use, we recalled the Milnor--Silverman fixed-point normal form 
\eqref{eq:MSnormal} (see \cite[(3-4) and Appendix~C, (C-1)]{Milnor93} and
\cite[Lemma~5.3(i)]{Silverman98})
as the explicit representative $\phi_{F'}$ above, rather than by invoking a general
algebro-geometric lifting argument.
Under the weighted scaling $(\delta,\varepsilon)\mapsto(c\delta,c^2\varepsilon)$,
the ratio $\delta^2/\varepsilon$ is invariant.  This ratio is the
algebraic counterpart of the intrinsic prehole parameter
$\ph$ in Section \ref{sec:residuemoduli}.
\end{remark}

We now prove $q_\phi=p$ and the second assertion \eqref{eq:rt} in
(\textit{\ref{item:intrinsic}}).
Let us see
\begin{gather*}
 \xi_0^{\phi_{F'}}=\xi_{1,|\varepsilon|^{1/2}}; 
\end{gather*}
indeed, after a further finite ramified base change of $R'$ if necessary, choose
$s\in R'$ with $s^2=\varepsilon$. Solving $\phi_{F'}'(z)=0$, we obtain
the two critical points $c_1,c_2$ of $\phi_{F'}$ as
\begin{gather*}
 \{c_1,c_2\}
 =\left\{\frac{(1-\mu_2)(1\pm s)}{\mu_2(\mu_1-1)}\right\}
 =\left\{1+\frac{s\bigl(s\pm(1-\mu_2)\bigr)}{\mu_2(\mu_1-1)}\right\},
 \quad\text{so that}\\
 |c_1-1|=|c_2-1|
 =\left|\frac{s\bigl(s\pm(1-\mu_2)\bigr)}{\mu_2(\mu_1-1)}\right|
 =|s|=|\varepsilon|^{1/2},\\
 |c_1-c_2|
 =\left|\frac{2s(1-\mu_2)}{\mu_2(\mu_1-1)}\right|
 =|s|=|\varepsilon|^{1/2},
\end{gather*}
since $|s|<1$, $|2|=|\mu_2|=1$, and $|\mu_1-1|=|1-\mu_2|=1$ by
\begin{gather*}
 \widehat{\mu_1-1}=\omega-1\ne0,
 \quad
 \widehat{1-\mu_2}=1-\omega^{-1}\ne0.
\end{gather*}
Under $\Char k=0$, 
the nearest point retraction $\xi_0^{\phi_{F'}}$ of $\xi_g$ to $\sR(\phi_{F'})=[c_1,c_2]$ is
$\xi_{1,|\varepsilon|^{1/2}}$.

Next, writing $\spe_{R'}:\bP^1(F')\to\bP^1(k)$ for closed-fiber specialization, let us see
\begin{gather}
 \as(\gamma)
 =\spe_{R'}\Bigl(\frac{\delta^2}{\varepsilon}\Bigr)
 \in\bP^1(k); \label{eq:dels}
\end{gather}
indeed, noting that $\mu_2=(1-\varepsilon)/\mu_1$ and
$\delta=\mu_1^p-1$, we have
\begin{gather*}
 \mu_2^p-1
 =\mu_1^{-p}\bigl((1-\varepsilon)^p-\mu_1^p\bigr)
 =-\mu_1^{-p}\bigl(\delta+p\varepsilon+\varepsilon^2h\bigr),\quad\text{so that}\\
 -\frac{(\mu_1^p-1)(\mu_2^p-1)}{1-\mu_1\mu_2}
 =\frac{\delta\bigl(\delta+p\varepsilon+\varepsilon^2h\bigr)}{\mu_1^p\varepsilon},
\end{gather*} 
for some $h\in R'$, and then also also by \eqref{eq:asmultiplier},
\begin{gather*}
\begin{aligned}
 \delta=0
 &\Rightarrow
 \as(\gamma)
 =\spe_{R'}\Bigl(
 \frac{\delta(\delta+p\varepsilon+\varepsilon^2h)}
 {\mu_1^p\varepsilon}\Bigr)
 =0
 =\spe_{R'}\!\Bigl(\frac{\delta^2}{\varepsilon}\Bigr),\\
 |\varepsilon|<|\delta|
 &\Rightarrow
 \widehat{\Bigl(
 1+\frac{p\varepsilon}{\delta}
 +\frac{\varepsilon^2h}{\delta}\Bigr)}=1
 \Rightarrow
 \as(\gamma)
 =\spe_{R'}\biggl(
 \frac{\delta^2}{\mu_1^p\varepsilon}
 \Bigl(1+\frac{p\varepsilon}{\delta}
 +\frac{\varepsilon^2h}{\delta}\Bigr)\biggr)
 =\spe_{R'}\Bigl(\frac{\delta^2}{\varepsilon}\Bigr),\quad\text{and}\\
 |\varepsilon|\ge|\delta|>0
 &\Rightarrow
 \left|
 \frac{\delta(\delta+p\varepsilon+\varepsilon^2h)}
 {\mu_1^p\varepsilon}
 \right|<1
 \quad\text{and}\quad
 \left|\frac{\delta^2}{\varepsilon}\right|<1\\
 &\Rightarrow
 \as(\gamma)
 =\spe_{R'}\Bigl(
 \frac{\delta(\delta+p\varepsilon+\varepsilon^2h)}
 {\mu_1^p\varepsilon}\Bigr)
 =0
 =\spe_{R'}\Bigl(\frac{\delta^2}{\varepsilon}\Bigr).
\end{aligned}
\end{gather*}

To see how the intrinsic reduction of
$\phi_{F'}^p$ changes along the segment from $\xi_g=\xi_{\phi_{F'}}$ toward
$\xi_0^{\phi_{F'}}=\xi_{1,|\varepsilon|^{1/2}}$, choose
$c\in R'\setminus\{0\}$ such that
\begin{gather*}
 \max\{|\delta|,|\varepsilon|^{1/2}\}\le|c|<1,\quad\text{so in particular that}
\quad 0<|\varepsilon|<|c|,
\end{gather*}
and put
\begin{gather*}
 N(z)=N_c(z):=\frac{1}{c}(z-1), \quad\text{so that}\quad N^{-1}(z)=1+cz.
\end{gather*}
Then, in the $N_*$-coordinate on $T_{\xi_{1,|c|}}\sP^1$, the intrinsic
reduction of $\phi_{F'}^p$ at $\xi_{1,|c|}$ is identified with
$\widetilde{N\circ\phi_{F'}^p\circ N^{-1}}(\zeta)$, which we compute as
\begin{align}
 \widetilde{N\circ\phi_{F'}^p\circ N^{-1}}(\zeta)
 &= \zeta+\widehat{\Bigl(\frac{\delta}{c}\Bigr)}
 +\widehat{\Bigl(\frac{\varepsilon}{c^2}\Bigr)}\zeta^{-1} \label{eq:semiaffineremaining}
\\
&=\begin{cases}
   \zeta &\text{if }|c|>\max\{|\delta|,|\varepsilon|^{1/2}\},\\
   \zeta +1 &\text{if } c=\delta\text{ and }|\delta|>|\varepsilon|^{1/2},\\
   \zeta +\widehat{\bigl(\delta/c\bigr)} +\zeta^{-1} &\text{if } c^2=\varepsilon\text{ and }|\delta|\le|\varepsilon|^{1/2}
  \end{cases}
 \quad\text{in }k(\zeta)\label{eq:interval}
\end{align}
after a finite ramified base change of $R'$ in the third case if necessary;
indeed, noting that $\delta=\mu_1^p-1$, we have
\begin{align*}
 N\circ\phi_{F'}^p\circ N^{-1}(z)
 &=
 \frac{1}{c}\bigl(\phi_{F'}^p(1+cz)-1\bigr)\notag\\
 &=
 \frac{1}{c}\biggl(
 \sum_{\ell=0}^{p-1}\mu_1^{p-1-\ell}
 \Bigl(\phi_{F'}^{\ell+1}(1+cz)-\mu_1\phi_{F'}^\ell(1+cz)\Bigr)
 +\mu_1^p\cdot(1+cz)-1\biggr)\notag\\
 &=
 \mu_1^p z+\frac{\delta}{c}
 +\frac{\mu_1^{p-1}}{c}\sum_{\ell=0}^{p-1}\mu_1^{-\ell}
 \bigl(\phi_{F'}^{\ell+1}(1+cz)-\mu_1\cdot\phi_{F'}^\ell(1+cz)\bigr)\quad\text{in }R'(z),
\end{align*}
and noting also the normal form \eqref{eq:MSnormal}, we have
\begin{gather}
 \phi_{F'}(1+cz)-\mu_1\cdot(1+cz)
 =
 \frac{
 \mu_1\frac{\varepsilon}{c}(1+cz)^2
 }{
 (1-\varepsilon)z+
 \frac{\varepsilon}{c}\frac{\mu_1}{1-\mu_1}
 }\quad\text{in }R'(z),\label{eq:difference}
\end{gather}
so that, on one hand,
\begin{gather*}
\widetilde{\Bigl(\frac{1}{c}\bigl(\phi_{F'}(1+cz)-\mu_1\cdot(1+cz)\bigr)\Bigr)}(\zeta)
=\omega\cdot\widehat{\Bigl(\frac{\varepsilon}{c^2}\Bigr)}\zeta^{-1}\quad\text{in }k(\zeta)\quad\text{and}\\
\widetilde{\bigl(\phi_{F'}(1+cz)\bigr)}(\zeta)-\omega= 
\widetilde{\bigl(\phi_{F'}(1+cz)-\mu_1\cdot(1+cz)\bigr)}(\zeta)
=\frac{\omega\cdot 0}{
 \zeta+0\cdot\frac{\omega}{1-\omega}
 }=0\quad\text{in }k(\zeta), 
\end{gather*}
the latter of which, since $\omega^j\ne1$ for $1\le j<p$, implies by the
composition formula \eqref{eq:depthcomp} and \eqref{eq:FKN} that, for every
$1\le\ell\le p-1$,
\begin{gather*}
 \widetilde{\phi_{F'}^\ell(1+cz)}
 =\bigl(\widetilde{\phi_{F'}}_{\xi_g}\bigr)^{\ell-1}(\omega)
 =\omega^\ell\quad\text{in }k(\zeta). 
\end{gather*}
Moreover, rewriting \eqref{eq:difference} as
\begin{gather*}
 \frac{1}{c}(\phi_{F'}(w)-\mu_1\cdot w)
 =
 \frac{
 \mu_1\frac{\varepsilon}{c}w^2
 }{
 (1-\varepsilon)(w-1)+ \varepsilon\frac{\mu_1}{1-\mu_1}
 }\quad\text{in }R'(w),
\end{gather*}
for every $1\le\ell\le p-1$, we inductively compute 
\begin{multline*}
 \widetilde{\Bigl(\frac{1}{c}\bigl(\phi_{F'}^{\ell+1}(1+cz)-\mu_1\cdot\phi_{F'}^\ell(1+cz)\bigr)\Bigr)}=\widetilde{
 \biggl(\frac{ 
 \mu_1\frac{\varepsilon}{c}(\phi_{F'}^\ell(1+cz))^2
 }{
 (1-\varepsilon)(\phi_{F'}^\ell(1+cz)-1)+ \varepsilon\frac{\mu_1}{1-\mu_1}
 }\biggr)}
=\frac{\omega^{2\ell+1} \cdot 0}{\omega^\ell-1}=0.
\end{multline*}
Hence \eqref{eq:semiaffineremaining} holds.

Let us conclude the proof of Theorem \ref{th:arc}(\ref{item:intrinsic}).
Suppose first that $|\delta|>|\varepsilon|^{1/2}$.
Then using the second case in \eqref{eq:interval},
we have $\phi_{F'}^p(\xi_{1,|\delta|})=\xi_{1,|\delta|}$ and
$\deg_{\xi_{1,|\delta|}}(\phi_{F'}^p)=1$ and, under $\Char k=0$, we also have
\begin{gather*}
 \phi_{F'}^p(1)=1\not\in U\Bigl(\widetilde{(\phi_{F'}^p)}_{\xi_{1,|\delta|}}(\widevec{\xi_{1,|\delta|}1})\Bigr),  
\end{gather*}
which together with the non-archimedean argument principle \eqref{eq:FKN} yields
\begin{gather*}
\dep_{\widevec{\xi_{1,|\delta|}1}}\bigl(\widetilde{(\phi_{F'}^p)}_{\xi_{1,|\delta|}}\bigr)
=\bigl((\phi_{F'}^p)^*\delta_1\bigr)\bigl(U(\widevec{\xi_{1,|\delta|}1})\bigr)
\ge\bigl((\phi_{F'}^p)^*\delta_1\bigr)(\{1\})\ge 1, 
\end{gather*}
and in turn $\xi_{1,|\delta|}\in \sJ(\phi_{F'})$ 
(see \cite[Theorem 2.1]{Kiwi14}, or \cite[Theorem 8.7]{BenedettoBook} in the present generality). Then also using the first case in \eqref{eq:interval}, 
by Theorem \ref{th:quad}, we have 
$\xi_1^{\phi_{F'}}=\xi_{1,|\delta|}\neq\xi_{1,|\varepsilon|^{1/2}}=\xi_0^{\phi_{F'}}$ 
and $q_{\phi_{F'}}=p$, which with the definition of $\ph$ and
\eqref{eq:dels} yields
\begin{gather*}
 \ph(\phi)=\ph(\phi_{F'})
 =\infty=\as(\gamma).
\end{gather*}
Suppose alternatively that $|\delta|\le|\varepsilon|^{1/2}$.
Then using the first and third cases in \eqref{eq:interval}, we have $q_{\phi_{F'}}=p$
and $\xi_1^{\phi_{F'}}=\xi_{1,|\varepsilon|^{1/2}}=\xi_0^{\phi_{F'}}$.
On the other hand, by the third case in \eqref{eq:interval},
where $c^2=\varepsilon$, in homogeneous coordinates $[X:Y]$
with $\zeta=Y/X$, the intrinsic pullback divisor $D_{\phi_{F'}}$ is computed as
\begin{gather*}
 N_*D_{\phi_{F'}}
 =\bigl(\zeta+\widehat{\delta/c}+\zeta^{-1}\bigr)^*(0)
 =\bigl(X^2+(\widehat{\delta/c})XY+Y^2=0\bigr), 
\end{gather*}
which together with Proposition \ref{th:divisor} and \eqref{eq:dels} yields
\begin{gather*}
 \ph(\phi)=\ph(\phi_{F'})
=\widehat{\Bigl(\frac{\delta/c}{1}\Bigr)^2}
=\widehat{\Bigl(\frac{\delta^2}{\varepsilon}\Bigr)}=\as(\gamma).
\end{gather*}
Now the proof of (\textit{\ref{item:intrinsic}}) is complete. \qed

\subsection{Proof of Theorem \ref{th:arc}(\textit{\ref{item:limits}})}
We retain the setting and notation of Subsection \ref{sec:proofintrinsic}
for each given $x$ and for each picked $\gamma$.
Let us see \eqref{eq:modlim} in (\textit{\ref{item:graph}}) first.
By (\ref{item:intrinsic}), we have $q_{\phi_{F'}}=p$,
and the description of $\MHRL(\phi_{F'})$ and Theorem \ref{th:quad},
for every $j\ge p$, 
\begin{gather*}
 \widehat{\Bigl(\bigl((N_c\circ M)\circ\phi\circ(N_c\circ M)^{-1}\bigr)^j\Bigr)}
 =\widehat{N_c\circ \phi_{F'}^j\circ N_c^{-1}}\in\bP(V_{2^j,k})^{\mathrm{s}}(k),
\end{gather*}
where $c$ is chosen as in the second or third case of \eqref{eq:interval}.
Then by the compatibility of the quotient morphism with specialization and the properness of $\overline{\rM}_{2^j,k}$ (see Fact \ref{th:git}),
\begin{align*}
 &\biggl[\widehat{\Bigl(\bigl((N_c\circ M)\circ\phi\circ(N_c\circ M)^{-1}\bigr)^j\Bigr)}\biggr]\\
 =&
 \spe_{R'}\Bigl(
 \bigl[N_c\circ\phi_{F'}^j\circ N_c^{-1}\bigr]
 \Bigr)
=
 \spe_{R'}\Bigl(
 \Phi_{2,j,k}\bigl(\gamma_{R'}|_{\eta_{R'}}\bigr)
 \Bigr)
=
 \spe_R\bigl(\Phi_{2,j,k}\circ(\gamma|_{\eta_R})\bigr)
\quad\text{in }\overline{\rM}_{2^j,k}(k),
\end{align*}
Hence \eqref{eq:modlim} in (\textit{\ref{item:graph}}) holds. Moreover,
let $\gamma'$ be another antispecialization arc degenerating to $x$, and let
$\psi$ be a Puiseux lift of $\gamma'$, so the assertion similar to that for $\gamma$
also holds for $\gamma'$. After $\PGL(2,\bK)$-conjugations of
$\phi$ and $\psi$ if necessary, for every $j\ge p$, we have
\begin{align*}
 &\spe_R\bigl(\Phi_{2,j,k}\circ(\gamma|_{\eta_R})\bigr)
 =\spe_R\bigl(\Phi_{2,j,k}\circ(\gamma'|_{\eta_R})\bigr)\\
 \Leftrightarrow\,&
 \widehat{\phi^j}\text{ and }\widehat{\psi^j}\text{ are }\PGL(2,k)\text{-conjugate}\\
 \Leftrightarrow\,&
 \text{the intrinsic reduction data of $\phi^j$ and $\psi^j$ at
 $\xi_1^\phi$ and $\xi_1^\psi$, respectively, are projectively equivalent}\\
 \Leftrightarrow\,&
 \ph(\phi)=\ph(\psi)
 \Leftrightarrow
 \as(\gamma)=\as(\gamma').
\end{align*}
The first equivalence follows from \eqref{eq:modlim} and Fact \ref{th:git}(\ref{item:git}).
The second equivalence is clear (see Fact \ref{th:algred} and Remark \ref{th:algdepth})
the third equivalence is by Theorem \ref{th:residuemoduli}, and the last equivalence
is by (\ref{item:intrinsic}). Hence \eqref{eq:tclass} holds. Now the proof of
(\textit{\ref{item:limits}}) is complete. \qed

\subsection{Proof of Theorem \ref{th:arc}(\textit{\ref{item:graph}})}\label{sec:proofgraph}
We retain the setting and notation of Subsection \ref{sec:proofintrinsic}
for the fixed $x$. We have already seen \eqref{eq:modlim} in (\textit{\ref{item:graph}}).

For every $\lambda\in\bP^1(k)$, 
choosing $\tau\in k$ such that $\tau^2=\lambda$ when $\lambda\in k$, 
set 
\begin{gather*}
 \bigl(\mu_{1,\lambda}(t),\mu_{2,\lambda}(t)\bigr)
 :=
 \begin{cases}
 \displaystyle\biggl(
 \omega\Bigl(1+\frac{\tau}{p}t\Bigr),
 \frac{1-t^2}{\omega(1+(\tau/p)t)}
 \biggr) & \text{if }\lambda\in k\\
 \displaystyle\biggl(
 \omega\Bigl(1+\frac{1}{p}t\Bigr),
 \frac{1-t^4}{\omega(1+t/p)}
 \biggr) & \text{if }\lambda=\infty
 \end{cases}
 \in k[[t]]^2,
\end{gather*}
so that for some $h_\lambda\in k[[t]]$,
\begin{align*}
 \varepsilon_\lambda&:=1-\mu_{1,\lambda}(t)\mu_{2,\lambda}(t)
 =
 \begin{cases}
 t^2 & \text{if }\lambda\in k,\\
 t^4 & \text{if }\lambda=\infty,
 \end{cases}\\
 \delta_\lambda&:=\mu_{1,\lambda}(t)^p-1
 =
 \begin{cases}
 \tau t+t^2h_\lambda & \text{if }\lambda\in k,\\
 t+t^2h_\lambda & \text{if }\lambda=\infty.
 \end{cases}
\end{align*}
Then we have the induced quadratic rational function 
\begin{gather*}
 \phi_\lambda(z):=
 z\frac{(1-\mu_{1,\lambda})z+\mu_{1,\lambda}(1-\mu_{2,\lambda})}
 {\mu_{2,\lambda}(1-\mu_{1,\lambda})z+(1-\mu_{2,\lambda})}
 \in F(z)
\end{gather*}
in the normal form in \eqref{eq:MSnormal}; 
indeed, $\deg\phi_\lambda=2$ since
$1-\mu_{1,\lambda},1-\mu_{2,\lambda}\in R^\times$ and
$\varepsilon_\lambda\ne0$.
By the same argument as in Subsection \ref{sec:proofintrinsic}, using
the non-archimedean argument principle \eqref{eq:FKN} and
the slope formula \eqref{eq:slope}, we have
$\xi_g=\xi_{\phi_\lambda}$.
By properness of $\overline{\rM}_{2,k}$, $[\phi_\lambda]$ extends uniquely to
$\gamma_\lambda\in\overline{\rM}_{2,k}(R)$.
Noting, from the multiplier computation in Subsection \ref{sec:proofintrinsic},
that the remaining fixed-point multiplier $\mu_{3,\lambda}$ satisfies
\begin{gather*}
 \mu_{3,\lambda}
 =\frac{2-\mu_{1,\lambda}-\mu_{2,\lambda}}{\varepsilon_\lambda}
 \quad\text{and}\quad
 \mu_{3,\lambda}^{-1}\in\mathfrak{m}_R,
\end{gather*}
we have in the Milnor--Silverman fixed-point multiplier coordinates
\begin{gather*}
 \gamma_\lambda=
 \biggl[
 1+\frac{\mu_{1,\lambda}+\mu_{2,\lambda}}{\mu_{3,\lambda}}
 :\mu_{1,\lambda}+\mu_{2,\lambda}
 +\frac{\mu_{1,\lambda}\mu_{2,\lambda}}{\mu_{3,\lambda}}
 :\frac{1}{\mu_{3,\lambda}}
 \biggr],
\end{gather*}
so that
\begin{gather*}
 \gamma_\lambda|_{\eta_R}=[\phi_\lambda]\in\rM_{2,k}(F)
 \quad\text{and}\quad
 \spe_R(\gamma_\lambda|_{\eta_R})=[1:\omega+\omega^{-1}:0]=x,
\end{gather*}
so $\gamma_\lambda$ is an antispecialization arc degenerating to $x$, and moreover, by
\eqref{eq:dels}, we have
\begin{gather*}
 \as(\gamma_\lambda)
 =\spe_R\!\left(\frac{\delta_\lambda^2}{\varepsilon_\lambda}\right)
 =\begin{cases}
 \tau^2=\lambda & \text{if }\lambda\in k,\\
 \infty & \text{if }\lambda=\infty.
 \end{cases}
\end{gather*}

Fix $n\ge p$, and let us construct a bijection 
\begin{gather*}
 \vartheta_{p,n}:
 \biggl(\im_{\Sch}\Bigl(
 \pr_{p,n}:\pr_1^{-1}(x)
 \to\prod_{j=p}^{n}\overline{\rM}_{2^j,k}
 \Bigr)\biggr)(k)
 \to\bP^1(k);
\end{gather*}
then we will define 
\begin{gather*}
 \Theta_{p,n}:=\vartheta_{p,n}^{-1}. 
\end{gather*}
For this purpose, we note that
since the restriction of $\pr_{p,n}$ to
$\pr_1^{-1}(x)$ is projective and $k$ is algebraically closed,
\begin{gather*}
 \biggl(\im_{\Sch}\Bigl(
 \pr_{p,n}:\pr_1^{-1}(x)
 \to\prod_{j=p}^{n}\overline{\rM}_{2^j,k}
 \Bigr)\biggr)(k)
 =\pr_{p,n}\bigl(\pr_1^{-1}(x)(k)\bigr),
\end{gather*}
so that for every point
\begin{gather*}
 y\in\biggl(\im_{\Sch}\Bigl(
 \pr_{p,n}:\pr_1^{-1}(x)
 \to\prod_{j=p}^{n}\overline{\rM}_{2^j,k}
 \Bigr)\biggr)(k),
\end{gather*}
there is a point $z\in\pr_1^{-1}(x)(k)$ such that
$\pr_{p,n}(z)=y$.  We use the following standard
curve-through-a-point lemma, which follows by taking hyperplane sections
through the given point in an affine neighborhood:
\begin{quote}
\itshape
If $X$ is an irreducible finite-type $k$-scheme of positive
dimension, $U\subset X$ is dense and open, and $z\in X(k)$, then there is an
integral curve in $X$ through $z$ whose generic point lies in $U$.
\end{quote}
Put
\begin{gather*}
 U_{2,n,k}:=
 (\Id,\Phi_{2,2,k},\ldots,\Phi_{2,n,k})(\rM_{2,k})
 =\pr_1^{-1}(\rM_{2,k}),
\end{gather*}
which is open and dense in $\mathcal{G}_{2,n,k}$. By the lemma applied to
$X=\mathcal{G}_{2,n,k}$ and $U=U_{2,n,k}$,
we obtain an integral curve $C$ through $z$ whose generic point lies in
$U_{2,n,k}$.  Letting $\nu:\widetilde{C}\to C$ be the normalization and choosing
$\widetilde{z}\in\nu^{-1}(z)$ by finiteness of normalization so that 
$\widehat{\cO}_{\widetilde{C},\widetilde{z}}\cong k[[t]]=R$, we have
the antispecialization arc
\begin{gather*}
 \gamma_C:
 \Spec R
 =\Spec\widehat{\cO}_{\widetilde{C},\widetilde{z}}
 \to\widetilde{C}
 \overset{\nu}{\to} C
 \to\mathcal{G}_{2,n,k}
 \overset{\pr_1}{\to}\overline{\rM}_{2,k},
\end{gather*}
degenerating to $x$ such that
\begin{gather*}
 \Bigl(\spe_R\bigl(\Phi_{2,j,k}\circ(\gamma_C|_{\eta_R})\bigr)\Bigr)_{j=p}^{n}
 =y
\end{gather*}
(for background in algebraic geometry, see
\cite[I, Exercise~1.8, Theorems~3.9A, 5.4A--5.5A, and~6.2A;
II, Exercise~3.8]{Hartshorne77}).
Now the desired map is defined as
\begin{gather*}
 \vartheta_{p,n}(y):=\as(\gamma_C).
\end{gather*}
By \eqref{eq:tclass}, it is well defined and injective, and
by the construction of the arc $\gamma_\lambda$
from $\lambda\in\bP^1(k)$ above, it is also surjective.
Hence we have the desired bijection $\vartheta_{p,n}$.

Finally, let us realize this bijection by a morphism.
For an affine scheme $U$ and a regular function $f$ on $U$, write
$D(f):=U\setminus V(f)$. Let $\tau$ be an indeterminate and put
\begin{gather*}
 S_\tau:=D\Bigl(1+\frac{\tau}{p}t\Bigr)\subset\bA^2_{\tau,t,k},
 \quad\text{and}\quad
 (\mu_{1,\tau},\mu_{2,\tau}):=
 \biggl(\omega\Bigl(1+\frac{\tau}{p}t\Bigr),
 \frac{1-t^2}{\omega(1+(\tau/p)t)}\biggr).
\end{gather*}
Let $\phi_\tau$ be the normal form \eqref{eq:MSnormal} for 
the pair $(\mu_{1,\tau},\mu_{2,\tau})$ and let
$\Gamma_{\tau,n}$ be the closure of the graph of
\begin{gather*}
 S_\tau\cap D\bigl(t(1-\mu_{1,\tau})(1-\mu_{2,\tau})\bigr)\to
 \prod_{j=p}^n\overline{\rM}_{2^j,k},
 \quad
 (\tau,t)\mapsto\bigl([\phi_\tau^j]\bigr)_{j=p}^n
\end{gather*}
in
$S_\tau\times_k\prod_{j=p}^n\overline{\rM}_{2^j,k}$, with projection
$\pi_\tau:\Gamma_{\tau,n}\to S_\tau$. For every geometric point $a$ of the
line $t=0$, the same curve-through-a-point argument as above and
\begin{gather*}
 1-\mu_{1,\tau}\mu_{2,\tau}=t^2
 \quad\text{and}\quad
 \mu_{1,\tau}^p-1=\tau t+t^2h
\end{gather*}
show that every point of $\pi_\tau^{-1}(a,0)$ is represented by an
antispecialization arc with intrinsic antispecialization parameter $a^2$.
Hence \eqref{eq:tclass} shows that every geometric fiber of $\pi_\tau$ over
$t=0$ consists of one point. After replacing $S_\tau$ by an open neighborhood
of this line and $\Gamma_{\tau,n}$ by its inverse image, $\pi_\tau$ is proper,
birational, and quasi-finite, hence finite; since $S_\tau$ is normal and
$\Gamma_{\tau,n}$ is integral,
\begin{gather*}
 \pi_\tau:\Gamma_{\tau,n}\overset{\sim}{\to}S_\tau.
\end{gather*}
Restricting the inverse graph morphism to $t=0$ therefore defines
\begin{gather*}
 \theta_{\tau,p,n}:\bA^1_{\tau,k}\to
 \im_{\Sch}\Bigl(
 \pr_{p,n}:\pr_1^{-1}(x)
 \to\prod_{j=p}^n\overline{\rM}_{2^j,k}
 \Bigr)
 \quad\text{and}\quad
 \theta_{\tau,p,n}(a)=\Theta_{p,n}(a^2)
\end{gather*}
on $k$-points. Hence $\theta_{\tau,p,n}$ is invariant under
$\tau\mapsto-\tau$ and descends through
$\bA^1_{\tau,k}/\{\pm1\}\cong\bA^1_{\lambda,k}$, $\lambda=\tau^2$, to a morphism
\begin{gather*}
 \theta_{0,p,n}:\bA^1_{\lambda,k}\to
 \im_{\Sch}\Bigl(
 \pr_{p,n}:\pr_1^{-1}(x)
 \to\prod_{j=p}^n\overline{\rM}_{2^j,k}
 \Bigr)
\end{gather*}
whose map on $k$-points is $\Theta_{p,n}|_k$.

For the chart at infinity, let $s$ be an indeterminate and put
\begin{gather*}
 S_s:=D\Bigl(1+\frac{1}{p}t\Bigr)\subset\bA^2_{s,t,k},
 \quad\text{and}\quad
 (\mu_{1,s},\mu_{2,s}):=
 \biggl(\omega\Bigl(1+\frac{1}{p}t\Bigr),
 \frac{1-st^2-t^4}{\omega(1+t/p)}\biggr).
\end{gather*}
The same graph-closure argument, applied to \eqref{eq:MSnormal} for the pair $(\mu_{1,s},\mu_{2,s})$ on
$S_s\cap D\bigl(t(s+t^2)(1-\mu_{1,s})(1-\mu_{2,s})\bigr)$, and using
\begin{gather*}
 1-\mu_{1,s}\mu_{2,s}=st^2+t^4
 \quad\text{and}\quad
 \mu_{1,s}^p-1=t+t^2h,
\end{gather*}
defines, after restricting $S_s$ to an open neighborhood of $t=0$, a morphism
\begin{gather*}
 \theta_{\infty,p,n}:\bA^1_{s,k}\to
 \im_{\Sch}\Bigl(
 \pr_{p,n}:\pr_1^{-1}(x)
 \to\prod_{j=p}^n\overline{\rM}_{2^j,k}
 \Bigr)
\end{gather*}
whose value on $k$-points has parameter $s^{-1}$ for $s\ne0$ and parameter
$\infty$ for $s=0$. On $s=\lambda^{-1}\ne0$, the morphisms
$\theta_{0,p,n}$ and $\theta_{\infty,p,n}$ agree by \eqref{eq:tclass}; hence they
glue to a morphism
\begin{gather*}
 \Theta_{p,n}:\bP^1_k\to
 \im_{\Sch}\Bigl(
 \pr_{p,n}:\pr_1^{-1}(x)
 \to\prod_{j=p}^n\overline{\rM}_{2^j,k}
 \Bigr)
\end{gather*}
whose map on $k$-points is the bijection $\vartheta_{p,n}^{-1}$. Since this morphism is proper
and its image contains every closed point of the target, it is surjective; its
fibers contain at most one closed point by \eqref{eq:tclass}. Hence it is a
bijective morphism.

Since the bijective morphisms $\Theta_{p,n}$, $n\ge p$, are compatible under
truncations of tails, the truncation morphisms between their targets are
bijective. Being projective and bijective, they are finite. Thus the inverse limit in
\eqref{eq:Theta} exists as a scheme, and the morphisms $\Theta_{p,n}$ induce
by the universal property the morphism $\Theta_x$. Since its natural projection
to every finite level is $\Theta_{p,n}$, the morphism $\Theta_x$ is bijective.

Now the proof of (\textit{\ref{item:graph}}) is complete. \qed

\begin{remark}
For uniformity with the proof of (\textit{\ref{item:intrinsic}}), we used
the fixed-point normal form \eqref{eq:MSnormal} above rather than constructing
the arcs directly from the multiplier polynomial \eqref{eq:multiplierpoly}.
\end{remark}

\subsection{Proof of Theorem \ref{th:arc}(\textit{\ref{item:findet}})}
We retain the setting and notation of Subsection \ref{sec:proofintrinsic}
for the fixed $x$ and for each picked $\gamma$.

By (\textit{\ref{item:graph}}), the assertion on specializations for $j\ge p$
follows from the bijectivity of the truncation morphisms.
For $2\le j<p$, by the intrinsic depth computation in
(\textit{\ref{item:finite}}) of Section \ref{sec:stationary} and
$\widetilde{(\phi_{F'}^j)}_{\xi_g}=(\widetilde{\phi_{F'}}_{\xi_g})^j$,
the intrinsic reduction data of $\phi_{F'}^j$ at
$\xi_{\phi_{F'}}=\xi_g$ depend only on $x$ up to projective equivalence, and
$\phi_{F'}^j$ has intrinsic stable reduction at $\xi_g$. Hence, by the same argument as that proving \eqref{eq:modlim} in the proof of
(\textit{\ref{item:limits}}), using the compatibility of the quotient morphism
with specialization and the properness of $\overline{\rM}_{2^j,k}$,
\begin{gather*}
 \bigl[\widehat{\phi_{F'}^j}\bigr]
 =\spe_{R'}\bigl([\phi_{F'}^j]\bigr)
 =\spe_{R'}\Bigl(\Phi_{2,j,k}\bigl(\gamma_{R'}|_{\eta_{R'}}\bigr)\Bigr)
 =\spe_R\bigl(\Phi_{2,j,k}\circ(\gamma|_{\eta_R})\bigr)
 \quad\text{in }\overline{\rM}_{2^j,k}(k).
\end{gather*}
Thus the assertion also holds for $2\le j<p$.
The final assertion follows from Theorem \ref{th:arc}(\textit{\ref{item:intrinsic}})
and Theorem \ref{th:residuemoduli}.

Now the proof of Theorem \ref{th:arc} is complete. \qed

\section{Proofs of Theorems \ref{th:demarco} and \ref{th:compactification}}\label{sec:proofexcfib}

\subsection{Proof of Theorem \ref{th:demarco}}
Retain the notation of Subsection \ref{sec:excfib}, with $k$ algebraically
closed of characteristic $0$. Let $x\in\overline{\rM}_{2,k}(k)$ be an order
$p>1$ boundary point of $\rM_{2,k}$, and fix $n\ge p$.
Choose a primitive $p$th root of unity $\omega\in k^\times$ such that
$P_x(\omega)=0$, so that
\begin{gather*}
 x=[1:\omega+\omega^{-1}:0].
\end{gather*}

In the affine chart $s_1\ne0$ around $x$, put
\begin{gather*}
 (u,v):=\Bigl(\frac{s_2}{s_1}-(\omega+\omega^{-1}),\frac{s_3}{s_1}\Bigr).
\end{gather*}
The fixed-point multiplier polynomial associated to
$[1:\omega+\omega^{-1}+u:v]$ is
\begin{align*}
 P(T):=P_{[1:\omega+\omega^{-1}+u:v]}(T)=-(T-\omega)(T-\omega^{-1})+uT+v(T^3+2);
\end{align*}
we note that $\Phi_p$ is monic and that $\deg P=3$ over $k(u,v)$.

If $p=2$, then $\omega=-1$ and $\Phi_2(T)=T+1$, so
\begin{align*}
 \Res_T(P,\Phi_2)
 =(-1)^3P(-1)
 =-\bigl(-(-1+1)^2-u+v((-1)^3+2)\bigr)
 =-(-u+v)=u-v.
\end{align*}
If $p>2$, then $\deg\Phi_p$ is even and $\omega\ne\omega^{-1}$, so
\begin{align*}
 &\Res_T(P,\Phi_p)
 =(-1)^{3\deg\Phi_p}\prod_{\Phi_p(\zeta)=0}P(\zeta)\notag\\
 =&\bigl(u\omega+v(\omega^3+2)\bigr)
 \bigl(u\omega^{-1}+v(\omega^{-3}+2)\bigr)
 \cdot\prod_{\substack{\Phi_p(\zeta)=0\\
 \zeta\notin\{\omega,\omega^{-1}\}}}
 \bigl(-(\zeta-\omega)(\zeta-\omega^{-1})+u\zeta+v(\zeta^3+2)\bigr)\\
 \equiv& u^2\prod_{\substack{\Phi_p(\zeta)=0\\
 \zeta\notin\{\omega,\omega^{-1}\}}}
 \bigl(-(\zeta-\omega)(\zeta-\omega^{-1})+u\zeta\bigr)
 \pmod v,
\end{align*}
where the last product is a unit at $x$. Consequently, in the local ring
$\cO_{\overline{\rM}_{2,k},x}$,
\begin{gather*}
 \bigl(v,\Res_T(P,\Phi_p)\bigr)
 =\begin{cases}
 (v,u)&\text{if }p=2,\\
 (v,u^2)&\text{if }p>2.
 \end{cases}
\end{gather*}
For each $y=[s_1:s_2:s_3]$ in this affine chart $s_1\ne 0$ around $x$,
we have $P_y(T)=s_1P(T)$, so that $\Res_T(P_y,\Phi_p)=s_1^{\deg\Phi_p}\Res_T(P,\Phi_p)$,
and $s_3=s_1v$.
Hence over some neighborhood of $x$, $\widetilde{\mathcal{G}}_{2,n,k}$ from Subsection
\ref{sec:excfib} is the blow-up of the ideal
$\bigl(v,\Res_T(P,\Phi_p)\bigr)$, so that the rational
map to $\bP^1_k$
\begin{gather}
 (u,v)\mapsto
 -\frac{p^2\Res_T(P_y,\Phi_p)}
 {(2-\omega-\omega^{-1})\Phi_p'(\omega)\Phi_p'(\omega^{-1})
 s_1^{\deg\Phi_p-1}s_3}
 =-\frac{p^2\Res_T(P,\Phi_p)}
 {(2-\omega-\omega^{-1})\Phi_p'(\omega)\Phi_p'(\omega^{-1})v}
 \label{eq:asrationalmap}
\end{gather}
extends to a morphism on this blow-up.
Its value at the closed point of any lifted antispecialization arc is
$\as(\gamma)$ for the projected arc $\gamma$.

By the definition of $\widetilde{\mathcal{G}}_{2,n,k}$, its fiber $E_{p,n}(x)$ over $x$ is identified
with the following projective family, where $(a,b)\in k^2\setminus\{(0,0)\}$
and $b/a=\infty$ if $a=0$. If $p=2$, so $\omega=-1$, it is the family of lines
\begin{gather*}
 L_{b/a}:\quad 2as_1+as_2-bs_3=0, \quad b/a\in\bP^1(k),
\end{gather*}
corresponding to the projectivized tangent directions $[v:u]=[a:b]$. If
$p>2$, it is the family of conics
\begin{gather*}
 C_{(b/a)^2}:\quad
 b^2s_1s_3-a^2\bigl(s_2-(\omega+\omega^{-1})s_1\bigr)^2=0,
 \quad (b/a)^2\in\bP^1(k),
\end{gather*}
corresponding to the quadratic contact data $[v:u^2]=[a^2:b^2]$.

If $p=2$, then $\omega=-1$ so $x=[1:-2:0]$.
For each line $L_{b/a}$ such that $a\ne0$, we have the antispecialization arc
$\gamma_{b/a}:=[1:-2+bt:at]\in\overline{\rM}_{2,k}(R)$
degenerating to $x$, taking $(u,v)=(bt,at)$.
Noting that $\Res_T(P,\Phi_2)=u-v$ and $\Phi_2'(T)=1$, we have
\begin{gather}
 \as(\gamma_{b/a})
 =\spe_R\biggl(-\frac{2^2(bt-at)}
 {1\cdot1\cdot4\cdot1^0\cdot at}\biggr)
 =\spe_R\biggl(\frac{a-b}{a}\biggr)
 =\frac{a-b}{a}.
 \label{eq:excp2}
\end{gather}
Alternatively if $p>2$, then for each conic $C_{(b/a)^2}$ such that $a\ne0$, 
we have the antispecialization arc
$\gamma_{(b/a)^2}:=[1:\omega+\omega^{-1}+bt:a^2t^2]\in\overline{\rM}_{2,k}(R)$
degenerating to $x$, taking $(u,v)=(bt,a^2t^2)$.
Using the above computation of $\Res_T(P,\Phi_p)$,
we have
\begin{align}
 \as(\gamma_{(b/a)^2})
 &=\spe_R\Biggl(-\frac{p^2\cdot (bt\omega+a^2t^2(\omega^3+2))
 (bt\omega^{-1}+a^2t^2(\omega^{-3}+2))}
 {(2-\omega-\omega^{-1})\Phi_p'(\omega)\Phi_p'(\omega^{-1})a^2t^2}\times\notag\\
 &\times\prod_{\substack{\Phi_p(\zeta)=0\\
 \zeta\notin\{\omega,\omega^{-1}\}}}
 \bigl(-(\zeta-\omega)(\zeta-\omega^{-1})+bt\zeta+a^2t^2(\zeta^3+2)\bigr)\Biggr)\notag\\
 &=-\frac{p^2\cdot 1\cdot (b\omega+0)(b\omega^{-1}+0)}
 {(2-\omega-\omega^{-1})a^2}\cdot
 \frac{\prod_{\Phi_p(\zeta)=0,\,\zeta\notin\{\omega,\omega^{-1}\}}
 \bigl(-(\zeta-\omega)(\zeta-\omega^{-1})+0\bigr)}{\Phi_p'(\omega)\Phi_p'(\omega^{-1})}\notag\\
 &=-\frac{p^2}{2-\omega-\omega^{-1}}
 \left(\frac{b}{a}\right)^2\cdot\frac{(-1)^{\deg\Phi_p-2}}{-(\omega-\omega^{-1})^2}
 =\frac{p^2}{\displaystyle -\frac{(\omega-1)^2}{\omega}
 \Bigl(\frac{\omega^2-1}{\omega}\Bigr)^2}
 \Bigl(\frac{b}{a}\Bigr)^2\notag\\
 &=-\frac{p^2\omega^3}{(\omega^2-1)^2(\omega-1)^2}
 \Bigl(\frac{b}{a}\Bigr)^2.
 \label{eq:excp}
\end{align}
Since the above formulas are affine linear with nonzero slope
in the respective parameters $b/a$ and $(b/a)^2$, we have
the isomorphism of projective curves
$\as_{p,n}:E_{p,n}(x)\overset{\sim}{\to}\bP^1_k$,
\begin{gather*}
 \as_{p,n}(e):=
 \begin{cases}
 \displaystyle 1-\frac{b}{a}&\text{if }p=2,\ e=L_{b/a},\\[2pt]
 \displaystyle -\frac{p^2\omega^3}{(\omega^2-1)^2(\omega-1)^2}\Bigl(\frac{b}{a}\Bigr)^2
 &\text{if }p>2,\ e=C_{(b/a)^2},
 \end{cases}
\end{gather*}
which agrees with the morphism extending the rational map \eqref{eq:asrationalmap}
on the dense open subset $E_{p,n}(x)\setminus\{a=0\}$.

We now prove the existence of $\nu_n$ and the desired properties of it.
The scheme $\widetilde{\mathcal{G}}_{2,n,k}$ is normal: above an order $2$
point it is smooth, and above an order $p>2$ point its two charts are open
subsets of
\begin{gather*}
 \Spec k[u,s]\quad(v=u^2s)
 \quad\text{and}\quad
 \Spec\bigl(k[u,v,r]/(u^2-vr)\bigr).
\end{gather*}
The second chart is normal by Serre's criterion, since it is a hypersurface
whose singular locus is the origin.
Let $\Gamma\subset \widetilde{\mathcal{G}}_{2,n,k}\times_k\mathcal{G}_{2,n,k}$ be the closure of the
graph of the rational map induced by the iterate-graph morphism on
$\rM_{2,k}$, and let $\rho:\Gamma\to \widetilde{\mathcal{G}}_{2,n,k}$ be the first projection.
For every geometric point $e$ of $\widetilde{\mathcal{G}}_{2,n,k}$, the curve-through-a-point argument
in Subsection \ref{sec:proofgraph}, over an algebraically closed
extension of $k$ when necessary, represents every closed point of
$\rho^{-1}(e)$ by an arc whose generic point lies over $\rM_{2,k}$.
If $e$ lies over an order $p$ point $x$ with $p\le n$, the extension of
\eqref{eq:tdef} above shows that the projected arc has parameter
$\as_{p,n}(e)$. By Theorem
\ref{th:arc}(\textit{\ref{item:limits}}), all its iterate specializations
with $p\le j\le n$ are determined by $e$, and those with $2\le j<p$
depend only on $x$, by the proof of Theorem
\ref{th:arc}(\textit{\ref{item:findet}}).
If $x$ has order $p>n$, the latter argument applies to every $2\le j\le n$.
At the remaining boundary points, uniqueness follows from Theorem
\ref{th:acyclicmoduli}, and over $\rM_{2,k}$ it is immediate.
Thus every geometric fiber of $\rho$ consists of one point. Since $\rho$
is projective and birational, it is finite, and the normality of $\widetilde{\mathcal{G}}_{2,n,k}$
implies that $\rho$ is an isomorphism. Composing its inverse with the
second projection defines
\begin{gather*}
 \nu_n:\widetilde{\mathcal{G}}_{2,n,k}\to\mathcal{G}_{2,n,k}.
\end{gather*}
This morphism is projective and birational, and its image is all of
$\mathcal{G}_{2,n,k}$ because it contains the dense iterate graph.
On each exceptional fiber it is injective by \eqref{eq:tclass} and the
isomorphism $\as_{p,n}$; away from the blow-up centers,
injectivity follows from the first projection.
The same statements hold over every algebraically closed extension of
$k$, so $\nu_n$ is finite and bijective, and hence a homeomorphism.
Since $\widetilde{\mathcal{G}}_{2,n,k}$ is normal and $\nu_n$ is finite and birational,
$\nu_n$ is the normalization morphism of $\mathcal{G}_{2,n,k}$.
It is unique because it agrees with the iterate-graph morphism on the
dense open subset $\rM_{2,k}$.
In particular, $\Phi_{2,n,k}$ is regular away from the blow-up centers,
whereas \eqref{eq:tclass} with $j=n$ shows that it is indeterminate at each
center. This proves the description of its indeterminacy locus stated
in Subsection \ref{sec:moduli}
(for background in algebraic geometry, see
\cite[II, Exercise~3.5(b), Theorem~8.22A, and Proposition~8.23(b);
III, Exercise~11.2 and the proof of Corollary~11.4]{Hartshorne77}).

For every point $e$ in the dense open subsets considered above, Theorem
\ref{th:arc}(\textit{\ref{item:graph}}) yields
\begin{gather*}
 \pr_{p,n}\bigl(\nu_n(e)\bigr)
 =\bigl(\pi_n\circ\Theta_x\bigr)\bigl(\as_{p,n}(e)\bigr).
\end{gather*}
Hence, since both sides are morphisms, we have
\begin{gather*}
 \pr_{p,n}\circ \nu_n|_{E_{p,n}(x)}
 =\bigl(\pi_n\circ\Theta_x\bigr)\circ\as_{p,n}
\end{gather*}
on all of $E_{p,n}(x)$.
For every $e\in E_{p,n}(x)(k)$ and every antispecialization arc
$\widetilde{\gamma}_n\in \widetilde{\mathcal{G}}_{2,n,k}(R)$ degenerating to $e$, set
$\gamma:=\pr_1\circ \nu_n\circ\widetilde{\gamma}_n$ as in
Subsection \ref{sec:excfib}. Since $\as_{p,n}$ is the
restriction of the extension of \eqref{eq:tdef}, for every Puiseux lift
$\phi\in\bK(z)$ of $\gamma$, Theorem
\ref{th:arc}(\textit{\ref{item:intrinsic}}) implies
\begin{gather*}
 \as_{p,n}(e)=\as(\gamma)=\ph(\phi).
\end{gather*}

Now the proof of Theorem \ref{th:demarco} is complete. \qed

\begin{remark}\label{rem:demarco}
When $k=\bC$, let $B_n$ and $h_n:B_n\to\Gamma_n$ denote DeMarco's
blow-up space and regular homeomorphism in
\cite[\S7, Theorem~7.1]{DeMarco07}, where
$\Gamma_n=\mathcal{G}_{2,n,\bC}$ by the definition of the iterate graphs.
This construction has the same centers, namely the order $p$ boundary points
with $2\le p\le n$, and the same local ideals $(u,v)$ for $p=2$ and
$(u^2,v)$ for $p>2$ (see \cite[\S7, p.~341]{DeMarco07}).
Since the centers are distinct, the successive blow-ups there are canonically
isomorphic to the blow-up defining $\widetilde{\mathcal{G}}_{2,n,\bC}$ here, so that
there is an isomorphism
\begin{gather*}
 \iota_n:\widetilde{\mathcal{G}}_{2,n,\bC}\overset{\sim}{\to}B_n
 \quad\text{over }\overline{\rM}_{2,\bC}.
\end{gather*}
The morphisms $h_n\circ\iota_n$ and $\nu_n$ restrict to the same
iterate-graph morphism on the dense open subset $\rM_{2,\bC}$;
since their source is integral and their target is separated, we have
$h_n\circ\iota_n=\nu_n$.

By the local description of DeMarco's parameter in
\cite[\S7, proof of Theorem~7.1, pp.~343--344]{DeMarco07}\footnote{In the displayed formula
for the case $q>2$ on p.~344 of \cite{DeMarco07}, $(\alpha(t)^2-1)^2$ and the
factor $(\zeta^2-1)$ in the denominator appear to be intended as
$(\alpha(t)^q-1)^2$ and $(\zeta^2-1)^2$, respectively;
cf.~\cite[(5.6)]{DeMarco07}.},
the formulas \eqref{eq:excp2} and \eqref{eq:excp}
yield the coincidence of $\as_{p,n}$ with DeMarco's
complex-dynamical parameter $\tau^2$, so all the three parameters coincide.
\end{remark}

\subsection{Proof of Theorem \ref{th:compactification}}
\label{sec:proofcompactification}

Let $\mathcal{I}_n$ denote the ideal sheaf whose blow-up is
$\widetilde{\mathcal{G}}_{2,n,k}$ in Definition \ref{th:blowupmodels}.
For $n\ge m\ge2$, the ideals $\mathcal{I}_n$ and $\mathcal{I}_m$ agree on
a neighborhood of every order $p$ boundary point with $2\le p\le m$,
and $\mathcal{I}_m$ is the unit ideal away from these points.
Hence the pullback of $\mathcal{I}_m$ to
$\widetilde{\mathcal{G}}_{2,n,k}$ is invertible, so the universal property
of the blow-up defines
\begin{gather*}
 \widetilde{\rho}_{n,m,k}:
 \widetilde{\mathcal{G}}_{2,n,k}\to\widetilde{\mathcal{G}}_{2,m,k}
 \quad\text{over }\overline{\rM}_{2,k}.
\end{gather*}
Every morphism over $\overline{\rM}_{2,k}$ between these models is the
identity on the common dense open subset $\rM_{2,k}$, so it is unique
because the source is integral and the target is separated.
The same argument shows that, for $n\ge m\ge\ell\ge2$,
\begin{gather*}
 \widetilde{\rho}_{m,\ell,k}\circ\widetilde{\rho}_{n,m,k}
 =\widetilde{\rho}_{n,\ell,k},
 \quad
 \widetilde{\rho}_{n,n,k}=\Id,
 \quad\text{and}\quad
 \nu_m\circ\widetilde{\rho}_{n,m,k}
 =\rho_{n,m,k}\circ\nu_n.
\end{gather*}
Since each $\nu_n$ is bijective on $k$-points, their compatibility implies
the bijection between the inverse limits of $k$-points in
(\textit{i}).

For a fixed order $p$ boundary point $x$, the defining ideals agree
near $x$ at every level $n\ge p$, so $\widetilde{\rho}_{n,m,k}$ restricts
to an isomorphism $E_{p,n}(x)\to E_{p,m}(x)$ for $n\ge m\ge p$.
Since the morphisms $\as_{p,n}$ are restrictions of extensions of the
same rational map \eqref{eq:asrationalmap}, we have
\begin{gather*}
 \as_{p,m}\circ\widetilde{\rho}_{n,m,k}|_{E_{p,n}(x)}
 =\as_{p,n}.
\end{gather*}
At the levels $n<p$ the blow-up has no center at $x$, so its fiber over
$x$ is a single point. Thus Theorem \ref{th:demarco} identifies the
inverse limit of the fibers over $x$ with $\bP^1(k)$ through these
compatible parameters. If $x\in\overline{\rM}_{2,k}(k)$ is not an order
$p>1$ boundary point, every blow-up model is an isomorphism over a
neighborhood of $x$, and its fiber over $x$ is a single point.
Together with (\textit{i}), this proves (\textit{ii}).

If $y_1$ is an order $p$ boundary point, Theorem \ref{th:demarco}
with $n=p$ and the bijectivity in
Theorem \ref{th:arc}(\textit{\ref{item:graph}}) show that $y_p$
determines the parameter in (\textit{ii}), so $(y_1,y_p)$ determines
every coordinate. In the other case, (\textit{ii}) shows that $y_1$
alone determines the point. For every fixed $n\ge2$, choose an order
$p>n$ boundary point $x$. The fiber of $\mathcal{G}_{2,n,k}(k)$ over
$x$ is a singleton, whereas the fiber in (\textit{ii}) is $\bP^1(k)$,
so projection onto this finite level is not injective.

Finally, let $\beta:Y\to\overline{\rM}_{2,k}$ be a finite composition
of blow-ups. It is an isomorphism away from a finite set
$S\subset\overline{\rM}_{2,k}(k)$: the non-quasi-finite locus of this
proper birational morphism of integral surfaces consists of finitely
many curves contracted to points; taking $S$ to be their images, the
restriction of $\beta$ over $\overline{\rM}_{2,k}\setminus S$ is finite
and birational with normal target, hence an isomorphism.
There are order $p$ boundary points for every $p>1$, and points of
different orders are distinct, so some such point $x$ lies outside
$S$. By the description of the indeterminacy locus proved above,
$\Phi_{2,p,k}$ is indeterminate at $x$. Since $\beta$ is an isomorphism
near $x$, the rational map $\Phi_{2,p,k}\circ\beta$ cannot extend to a
morphism on $Y$. This completes the proof of (\textit{iii}).

Now the proof of Theorem \ref{th:compactification} is complete. \qed

\appendix
\section{Rigidity in the constant-or-acyclic case}
\label{sec:acyclicappendix}

Unless otherwise stated, we retain the notation and assumptions on $K$ and
$k$ from the beginning of Section \ref{sec:intro}.

We record rigidity results for the constant-or-acyclic case of Theorem \ref{th:quad},
which are used in Section \ref{sec:proofexcfib}.  Within this case, there are three possibilities for the
intrinsic reduction at $\xi_\phi$: it is non-constant of degree two, it is
constant, or it is a projective transformation of infinite order.  The first
case is immediate, the third has full intrinsic rigidity, and the constant
case has an explicit depth-rigidity statement.  We then record the common
moduli-theoretic consequence for the two cases that occur at boundary
points. For convenience, we recall below the field $\bK$ and the notion of a Puiseux lift.
 
\subsection{Intrinsic rigidity}

We use the intrinsic reduction data and simultaneous projective equivalence
introduced in Subsection \ref{sec:residue}.

\begin{appcorollary}[The degree-two reduction case]\label{th:acyc2}
Let $\phi\in K(z)$ be a quadratic rational function such that
$\tilde{\phi}_{\xi_\phi}$ has degree two.  Then for every $j\ge1$,
\begin{gather*}
 \widetilde{(\phi^j)}_{\xi_\phi}
 =(\widetilde{\phi}_{\xi_\phi})^j
 \quad\text{and}\quad
 \dep_{\Vv}\widetilde{(\phi^j)}_{\xi_\phi}=0
 \quad\text{for every }\Vv\in T_{\xi_\phi}\sP^1.
\end{gather*}
In particular, the intrinsic reduction data of $\phi$ at $\xi_\phi$ determine the full indexed
collection of intrinsic reduction data of all iterations.
\end{appcorollary}

\begin{proof}
This is the first case in the proof of Theorem \ref{th:quad}.  There
$\phi^{-j}(\xi_\phi)=\{\xi_\phi\}$ for every $j\ge1$, so every intrinsic
depth is zero.  Since $\phi(\xi_\phi)=\xi_\phi$, by functoriality of the tangent
maps, we have the displayed identity for the intrinsic reductions.
\end{proof}

\begin{appcorollary}[Rigidity in the acyclic projective-transformation case]
\label{th:acyclic}
Let $\phi,\psi\in K(z)$ be quadratic rational functions such that the
intrinsic reductions $\tilde{\phi}_{\xi_\phi}$ and
$\tilde{\psi}_{\xi_\psi}$ are projective transformations of infinite order.
If the intrinsic reduction data of $\phi$
and $\psi$ at $\xi_\phi$ and $\xi_\psi$, respectively, are projectively
equivalent, then the indexed collections of intrinsic
reduction data of $\phi^j$ and $\psi^j$ at $\xi_\phi$ and $\xi_\psi$,
respectively, labeled by $j=1,2,\ldots$, are simultaneously projectively
equivalent.
\end{appcorollary}

\begin{proof}
By case (\textit{\ref{item:infinite}}) in the proof of Theorem \ref{th:quad}
and the composition formula \eqref{eq:depthcomp}, every projective isomorphism
identifying the intrinsic reduction data of $\phi$ and $\psi$ at $\xi_\phi$
and $\xi_\psi$, respectively, also identifies those of all iterations.
\end{proof}

\begin{appcorollary}[When the intrinsic reduction is constant]\label{th:acycconst}
Let $\phi\in K(z)$ be a quadratic rational function such that
$\tilde{\phi}_{\xi_\phi}$ is constant. Then for every $j\ge1$ and every
$\Vv\in T_{\xi_\phi}\sP^1$, we have
\begin{gather*}
 \dep_{\Vv}\widetilde{(\phi^j)}_{\xi_\phi}
 =\deg(\phi^{j-1})\dep_{\Vv}\tilde{\phi}_{\xi_\phi}.
\end{gather*}
\end{appcorollary}

\begin{proof}
This is the depth computation in the constant case of the proof of Theorem
\ref{th:quad}.
\end{proof}

\subsection{Moduli-theoretic consequences}

For the remainder of the appendix, let $k$ be an algebraically closed field
and take
\begin{gather*}
 \bK:=\widehat{\overline{k((t))}},
 \quad F:=k((t)),
 \quad\text{and}\quad R:=F^\circ=k[[t]],
\end{gather*}
where $\bK$ is endowed with the extension of the $t$-adic absolute value on $F$.  Denote by
$\eta_R:=\Spec F$ the generic point of $\Spec R$.
The notation $\spe_R$ means taking the value at the unique closed point
of $\Spec R$ of the unique extension of a rational section of a
proper $R$-scheme.

\begin{appdefinition}\label{th:acyclicboundary}
A boundary point $x\in\overline{\rM}_{2,k}(k)$ of $\rM_{2,k}$ is called
{\em constant-or-acyclic} if the fixed-point multiplier polynomial $P_x$
in \eqref{eq:multiplierpoly} satisfies
\begin{gather*}
 \Res_T(P_x,\Phi_m)\ne0
 \quad\text{for every integer }m>1,
\end{gather*}
where $\Phi_m(T)\in k[T]$ is the image of the $m$th cyclotomic
polynomial in $\bZ[T]$.
\end{appdefinition}

We call an algebraic arc $\gamma\in\overline{\rM}_{2,k}(R)$ an
{\em antispecialization arc degenerating to} a constant-or-acyclic boundary point $x$ if
\begin{gather*}
 \gamma|_{\eta_R}\in\rM_{2,k}(F)\quad\text{and}\quad \spe_R(\gamma|_{\eta_R})=x.
\end{gather*}
For an antispecialization arc $\gamma$ degenerating to $x$,
a {\em Puiseux lift} of $\gamma$ is a quadratic rational function
$\phi\in\bK(z)$ such that, for some finite extension $F'/F$ contained in
$\overline F$, $\phi\in F'(z)$ and $[\phi]\in\rM_{2,k}(F')$ is the base change
of $\gamma|_{\eta_R}$ to $F'$.

By Theorem \ref{th:intrinsicstratification}(\textit{ii})--(\textit{iii}) in Appendix
\ref{sec:djrappendix}, for any quadratic rational function $\phi\in\bK(z)$
with $\spe_{\bK^\circ}([\phi])=x$ a boundary point, $x$ is constant-or-acyclic
if and only if $\tilde{\phi}_{\xi_\phi}$ is constant or a projective transformation
of infinite order. If $\Char k=\ell>0$, the condition for $m=\ell$ excludes
non-trivial unipotent transformations, since $\Phi_\ell(T)=(T-1)^{\ell-1}$ in $k[T]$.

The following is the constant-or-acyclic counterpart of Theorem
\ref{th:arc}(\textit{\ref{item:graph}}).

\begin{apptheorem}[Rigidity and regularity over a constant-or-acyclic boundary point]
\label{th:acyclicmoduli}
Let $x\in\overline{\rM}_{2,k}(k)$ be a constant-or-acyclic boundary point of
$\rM_{2,k}$. Then every $\Phi_{2,j,k}$ is regular at $x$ for $j\ge2$, and
for every antispecialization arc $\gamma$ degenerating to $x$,
\begin{gather*}
 \spe_R\bigl(\Phi_{2,j,k}\circ(\gamma|_{\eta_R})\bigr)
 =\Phi_{2,j,k}(x)
 \quad\text{for every }j\ge2.
\end{gather*}
For every such $\gamma$ and every Puiseux lift $\phi$ of $\gamma$, after a
$\PGL(2,\bK)$-conjugation of $\phi$ if necessary, we also have
\begin{gather*}
 [\widehat{\phi^j}]
 =\spe_R\bigl(\Phi_{2,j,k}\circ(\gamma|_{\eta_R})\bigr)
 =\Phi_{2,j,k}(x)
 \quad\text{in }\overline{\rM}_{2^j,k}(k)
 \quad\text{for every }j\ge2.
\end{gather*}
For every $n\ge2$, the fiber of
$\pr_1:\mathcal{G}_{2,n,k}\to\overline{\rM}_{2,k}$ over $x$
consists of the single reduced point
\begin{gather*}
 (x,\Phi_{2,2,k}(x),\ldots,\Phi_{2,n,k}(x)).
\end{gather*}
\end{apptheorem}

\begin{proof}
Let $\gamma$ be an antispecialization arc degenerating to $x$, and let $\phi$
be a Puiseux lift of $\gamma$.  After a $\PGL(2,\bK)$-conjugation of $\phi$ if necessary, 
we may assume
that $\xi_\phi=\xi_g$.  Let $f$ be the coefficient reduction representing the intrinsic stable
reduction of $\phi$ at $\xi_g$.  By Fact \ref{th:git}(\textit{\ref{item:git}}) and uniqueness of specialization in the proper
scheme $\overline{\rM}_{2,k}$,
\begin{gather*}
 q_{2,k}(f)=\spe_R(\gamma|_{\eta_R})=x.
\end{gather*}
By Fact \ref{th:algred}, Remark \ref{th:algdepth}, and Fact \ref{th:git}(\textit{\ref{item:git}}), the intrinsic
reduction data of $\phi$ are determined by $x$ up to projective equivalence. In particular, by the characterization following Definition \ref{th:acyclicboundary},
$\tilde{\phi}_{\xi_\phi}$ is either constant or a projective transformation of
infinite order, so $\phi$ is in the constant-or-acyclic case of Theorem \ref{th:quad},
so in particular $\MHRL(\phi^j)=\{\xi_\phi\}$ for every $j\ge1$.

If $\tilde{\phi}_{\xi_\phi}$ is a projective transformation of infinite order,
then by Corollary \ref{th:acyclic}, the projective-equivalence class of the
intrinsic reduction data of every $\phi^j$ is determined by that of $\phi$.
If $\tilde{\phi}_{\xi_\phi}$ is constant, then by Corollary \ref{th:acycconst},
for every $j\ge2$, the directions $\Vv$ with $\dep_{\Vv}\widetilde{(\phi^j)}_{\xi_\phi}>0$ are the same two
directions, both of depth $\deg(\phi^j)/2$, while the intrinsic reduction is constant
with value outside these two directions.  Since $\PGL(2,k)$ is transitive on
such triples, the projective-equivalence class of the intrinsic reduction
data of $\phi^j$ is again determined by the intrinsic reduction data of $\phi$. 
Consequently, for every $j\ge2$, by Fact
\ref{th:algred}, Remark \ref{th:algdepth}, and Fact \ref{th:git}(\textit{\ref{item:git}}), the specialization
\begin{gather*}
 y_j(x):=\spe_R\bigl(\Phi_{2,j,k}\circ(\gamma|_{\eta_R})\bigr)
\end{gather*}
depends only on $x$ and is well defined.

Fix $n\ge2$, and let $y$ be a closed point of the fiber of
$\pr_1:\mathcal{G}_{2,n,k}\to\overline{\rM}_{2,k}$ over $x$.
By the curve-through-a-point construction in Subsection \ref{sec:proofgraph},
we obtain an antispecialization arc $\gamma$ degenerating to
$x$, whose iterate graph arc specializes to $y$.  Hence
\begin{gather*}
 y=\bigl(x,y_2(x),\ldots,y_n(x)\bigr).
\end{gather*}
Thus the fiber of $\pr_1$ over $x$ has a single point as its
support, so $\pr_1$ is quasi-finite at the unique point above
$x$.  The quasi-finite locus is open.  Its complement is closed, and its
image under the proper morphism $\pr_1$ is closed and does not
contain $x$.  Hence, after shrinking an open neighborhood of $x$, the
morphism $\pr_1$ is proper and quasi-finite, and therefore
finite.  Since
$\overline{\rM}_{2,k}\cong\bP^2_k$ is normal, this finite birational morphism is
an isomorphism after shrinking the neighborhood.  Therefore every
$\Phi_{2,j,k}$ is regular at $x$, and the fiber above consists of the single reduced point
\begin{gather*}
 \bigl(x,\Phi_{2,2,k}(x),\ldots,\Phi_{2,n,k}(x)\bigr),
\end{gather*}
and $y_j(x)=\Phi_{2,j,k}(x)$.

Finally, the same $\PGL(2,\bK)$-conjugation sending $\xi_\phi$ to $\xi_g$ works for every
iteration by Theorem \ref{th:quad}. By intrinsic stability and
Fact \ref{th:git}(\textit{\ref{item:git}}), we therefore have
\begin{gather*}
 [\widehat{\phi^j}]
 =y_j(x)
 =\spe_R\bigl(\Phi_{2,j,k}\circ(\gamma|_{\eta_R})\bigr)
 =\Phi_{2,j,k}(x)
\end{gather*}
for every $j\ge2$.
\end{proof}

\subsection{Fibers of the blow-up spaces in characteristic zero}

Suppose that $\Char k=0$.

\begin{appcorollary}[Constant-or-acyclic fibers in characteristic zero]
\label{th:acyclicC}
Let $x\in\overline{\rM}_{2,k}(k)$ be a constant-or-acyclic boundary
point of $\rM_{2,k}$. Then for every $n\ge2$, the fiber
$(\pr_1\circ \nu_n)^{-1}(x)$ consists of a single point $e_{n,x}$,
and $\nu_n(e_{n,x})=(x,\Phi_{2,2,k}(x),\ldots,\Phi_{2,n,k}(x))$.
\end{appcorollary}

\begin{proof}
By Definitions \ref{th:orderboundary} and \ref{th:acyclicboundary}, $x$ is not
an order $p$ boundary point for any $p>1$, so the blow-up ideal in
Definition \ref{th:blowupmodels} is the unit ideal near $x$. Hence
$\pr_1\circ\nu_n$ is an isomorphism near $x$, so its fiber consists of a
single point $e_{n,x}$, whose image under $\nu_n$ is as stated by
Theorem \ref{th:acyclicmoduli}.
\end{proof}

\section{Intrinsic stratification in degree two}
\label{sec:djrappendix}

Unless otherwise stated, we retain the notation and assumptions on $K$ and
$k$ from the beginning of Section \ref{sec:intro}.

We relate the intrinsic reduction data of a quadratic rational
function to its specialization in dynamical moduli.  For a quadratic rational
function $\phi\in K(z)$, write $[\phi]\in\rM_2(K)$ for its
$\PGL(2,K)$-conjugacy class.  For every integer $d\ge2$ and every
$y\in\overline{\rM}_d(K)$, write
\begin{gather*}
 \spe_{K^\circ}(y)\in\overline{\rM}_{d,k}(k)
\end{gather*}
for the closed-fiber value of the unique $K^\circ$-point extending $y$, which
exists by properness.

For a quadratic rational function, let $\mu_1,\mu_2,\mu_3$ be the
fixed-point multipliers, counted with multiplicity.  The multiplier-symmetric
coordinates on $\rM_2\cong\bA^2_{\bZ}$ are
\begin{gather*}
 \sigma_1=\mu_1+\mu_2+\mu_3
 \quad\text{and}\quad
 \sigma_2=\mu_1\mu_2+\mu_2\mu_3+\mu_3\mu_1,
\end{gather*}
and we use the resulting homogeneous coordinates
$[\sigma_1:\sigma_2:1]$ on
$\overline{\rM}_2\cong\bP^2_{\bZ}$.  We have the following
explicit stratification.

\begin{apptheorem}[Intrinsic stratification in degree two]
\label{th:intrinsicstratification}
Let $\phi\in K(z)$ be a quadratic rational function.  Then the following
statements hold.
\begin{enumerate}[(i)]
\item $\deg\tilde{\phi}_{\xi_\phi}=2$ if and only if
$\spe_{K^\circ}([\phi])\in\bA^2(k)$.
\item $\tilde{\phi}_{\xi_\phi}$ is a non-trivial projective transformation if
and only if $\spe_{K^\circ}([\phi])=[1:x:0]$ for some $x\in k$.  More
precisely, if $A\in\GL_2(k)$ represents
$\tilde{\phi}_{\xi_\phi}$ under a projective identification
$T_{\xi_\phi}\sP^1\cong\bP^1_k$, then
$x=\tr(A)^2/\det(A)-2$.
\item $\tilde{\phi}_{\xi_\phi}$ is constant if and only if
$\spe_{K^\circ}([\phi])=[0:1:0]$.
\end{enumerate}
\end{apptheorem}

\begin{proof}
Choose $M\in\PGL(2,K)$ such that $M(\xi_\phi)=\xi_g$, and let
$f:=\widehat{M\circ\phi\circ M^{-1}}\in\bP(V_{2,k})(k)$.
A normalized homogeneous lift of $M\circ\phi\circ M^{-1}$ extends its
coefficient point to a $K^\circ$-point of $\bP(V_{2,\bZ})$ with closed fiber $f$. By intrinsic
stability at $\xi_\phi$ and Fact \ref{th:git}(\textit{\ref{item:git}}), composing
with the quotient morphism yields a $K^\circ$-point of
$\overline{\rM}_2$ extending $[\phi]$, so by uniqueness of this extension,
we have
\begin{gather*}
 q_{2,k}(f)=\spe_{K^\circ}([\phi]).
\end{gather*}
If $\deg\tilde{\phi}_{\xi_\phi}=2$, then
$f\in\Rat_{2,k}$, so (i) holds.

Suppose that $\deg\tilde{\phi}_{\xi_\phi}=1$.  By intrinsic stability, there is a unique
direction $\Vv$ with $\dep_{\Vv}\tilde{\phi}_{\xi_\phi}>0$, which is not fixed by $\tilde{\phi}_{\xi_\phi}$;
hence the latter is a non-trivial projective transformation.  Relabel the
fixed-point multipliers so that $\mu_1,\mu_2$ are those of the two fixed-point
branches specializing to its fixed-point divisor, counted with multiplicity,
and $\mu_3$ is the remaining one. Since the fixed-point divisor avoids
$\Vv$, the multipliers $\mu_1$ and $\mu_2$ lie in $K^\circ$ and reduce to
the corresponding multipliers of $\tilde{\phi}_{\xi_\phi}$.
Since $f$ is degenerate, $q_{2,k}(f)=\spe_{K^\circ}([\phi])$ lies on the
boundary, so $|\mu_3|>1$; otherwise both $\sigma_1$ and $\sigma_2$ would
lie in $K^\circ$ and $\spe_{K^\circ}([\phi])$ would lie in $\bA^2(k)$.
After dividing the homogeneous coordinates by $\mu_3$ and reducing, we have
\begin{gather*}
 [\sigma_1:\sigma_2:1]
 =\left[1+\frac{\mu_1+\mu_2}{\mu_3}:
 \mu_1+\mu_2+\frac{\mu_1\mu_2}{\mu_3}:
 \frac{1}{\mu_3}\right]
 \mapsto
 [1:\widehat{\mu_1}+\widehat{\mu_2}:0],
\end{gather*}
so
\begin{gather*}
 \spe_{K^\circ}([\phi])
 =[1:\widehat{\mu_1}+\widehat{\mu_2}:0].
\end{gather*}
If $A\in\GL_2(k)$ represents
$\tilde{\phi}_{\xi_\phi}$, the sum of the multipliers at its two fixed points,
counted with multiplicity, is
$\tr(A)^2/\det(A)-2$.  Hence
$x=\tr(A)^2/\det(A)-2$.  Conversely, since $k$ is algebraically closed, every $x\ne2$ is
$a+a^{-1}$ for some $a\in k^\times\setminus\{1\}$, while
$x=2$ is realized by a non-trivial unipotent projective transformation.
If the direction $\Vv$ with $\dep_{\Vv}\tilde{\phi}_{\xi_\phi}>0$ is chosen away from the fixed directions, the
coefficient point is stable, so the degree-one case fills the points $[1:x:0]$
on the line at infinity.

Finally, if $\tilde{\phi}_{\xi_\phi}$ is constant, by intrinsic stability, there are
two distinct directions $\Vv$ with $\dep_{\Vv}\tilde{\phi}_{\xi_\phi}>0$ and a constant value distinct from both.
Thus all corresponding stable coefficient points are projectively equivalent.
Since the degree-one case fills all points $[1:x:0]$, their common quotient
point is the remaining point $[0:1:0]$.  The three cases exhaust the stable
reductions of quadratic rational functions.  Hence (iii) and the converses in
(i)--(ii) hold.
\end{proof}

If $\Char k=0$, Theorem \ref{th:arc} concerns precisely the boundary points
$[1:\omega+\omega^{-1}:0]$ with $\omega$ a root of unity of order $p>1$.
Every other boundary point is constant-or-acyclic in the sense of Definition
\ref{th:acyclicboundary} and is covered by Theorem \ref{th:acyclicmoduli}.
More precisely, for a quadratic rational function $\phi$ with
$\spe_{K^\circ}([\phi])=x$, the intrinsic reduction $\tilde{\phi}_{\xi_\phi}$
is constant when $x=[0:1:0]$, a non-trivial unipotent projective transformation
of infinite order when $x=[1:2:0]$, and a semisimple projective transformation
of infinite order when $x=[1:a+a^{-1}:0]$ with $a\in k^\times$ not a root of
unity.

\begin{appproposition}[Intrinsic reduction and the fixed-point tree]
\label{th:fixedpointtree}
Let $\phi\in K(z)$ be a quadratic rational function, and let
$\Gamma_{\mathrm{Fix}}$ be the tree in $\sP^1$ spanned by its classical
fixed points. Then the following statements hold.
\begin{enumerate}[(i)]
\item If $\tilde{\phi}_{\xi_\phi}$ is a non-trivial unipotent projective
transformation, then $\xi_\phi\in\Gamma_{\mathrm{Fix}}$ and its valence in
$\Gamma_{\mathrm{Fix}}$ is two.
\item If $\tilde{\phi}_{\xi_\phi}$ is a non-trivial semisimple projective
transformation or is constant, then $\phi$ has three distinct classical
fixed points and $\xi_\phi$ is the unique branch point of
$\Gamma_{\mathrm{Fix}}$.
\end{enumerate}
\end{appproposition}

\begin{proof}
By the discussion following \eqref{eq:stable}, $\phi$ has intrinsic stable
reduction at the unique point $\xi_\phi$ of $\MHRL(\phi)$. For every
$\Vv\in T_{\xi_\phi}\sP^1$, the classical fixed-point divisor satisfies
(cf.\ \cite[Lemma~2.1, (3), and the proof of Lemma~2.2, (10)]{Rumely17})
\begin{gather*}
 \deg\bigl(\Fix_\phi|_{U(\Vv)}\bigr)
 =\dep_{\Vv}\tilde{\phi}_{\xi_\phi}
 +\ord_{\Vv}\bigl(\Fix_{\tilde{\phi}_{\xi_\phi}}\bigr).
\end{gather*}
Indeed, choose a projective identification carrying $\xi_\phi$ to $\xi_g$
and a normalized homogeneous lift $[F:G]$ whose coefficient reduction is
$H[\widetilde F:\widetilde G]$.  By Remark \ref{th:algdepth},
$\ord_{\Vv}(H)=\dep_{\Vv}\tilde{\phi}_{\xi_\phi}$, while the
reduction of the fixed-point equation $YF-XG=0$ is
$H(Y\widetilde F-X\widetilde G)=0$; the second factor is precisely the
fixed-point equation of $\tilde{\phi}_{\xi_\phi}$.  Hence the displayed
identity holds.

By intrinsic stability, in the non-trivial unipotent case, the unique
direction $\Vv$ with $\dep_{\Vv}\tilde{\phi}_{\xi_\phi}>0$ and the unique
fixed direction are distinct and contain, respectively, one and two
classical fixed points counted with multiplicity, so
$\xi_\phi\in\Gamma_{\mathrm{Fix}}$ and its valence is two. In the
non-trivial semisimple case, the unique direction $\Vv$ with
$\dep_{\Vv}\tilde{\phi}_{\xi_\phi}>0$ and the two fixed directions are
distinct. In the constant case, the two directions $\Vv$ with
$\dep_{\Vv}\tilde{\phi}_{\xi_\phi}>0$ and the constant fixed direction are
again distinct. In either case, the displayed identity shows that each of
these three directions contains exactly one classical fixed point, so
$\phi$ has three distinct classical fixed points and $\xi_\phi$ is the
unique branch point of $\Gamma_{\mathrm{Fix}}$.
\end{proof}

This intrinsic stratification recovers the crucial-set stratification of
Doyle--Jacobs--Rumely as follows. Let $\xi$ be the unique point in the
crucial set of $\phi$; see \cite[\S2.2]{DJR16}.

Let $\Gamma_{\mathrm{Fix}}$ be the tree in $\sP^1$ spanned by the classical
fixed points of $\phi$.  Recall from \cite[p.~2]{DJR16} that $\phi$ has
potential good reduction if the reduction at $\xi$ has degree two; potential
additive reduction if it has degree one and is projectively conjugate to
$\zeta\mapsto\zeta+b$ for some $b\in k^\times$, with
$\xi\in\Gamma_{\mathrm{Fix}}$; potential multiplicative reduction if it has
degree one and is projectively conjugate to $\zeta\mapsto a\zeta$ for some
$a\in k^\times\setminus\{1\}$, with $\xi$ a branch point of
$\Gamma_{\mathrm{Fix}}$; and potential constant reduction if it is constant,
with $\xi$ a branch point of $\Gamma_{\mathrm{Fix}}$.
By Proposition \ref{th:fixedpointtree}, these tree conditions hold at
$\xi_\phi$. By the characterization of the crucial point in
\cite[\S2.2.2]{DJR16}, we therefore have $\xi=\xi_\phi$.
Thus the four potential reduction types correspond exactly to the four
intrinsic possibilities at $\xi_\phi$: degree two, non-trivial unipotent,
non-trivial semisimple, and constant, respectively.  The following is therefore
immediate from Theorem \ref{th:intrinsicstratification}.

\begin{appcorollary}[{\cite[Theorem~1.1]{DJR16}}]
\label{th:djr}
Let $\phi\in K(z)$ be a quadratic rational function.  Then
\begin{enumerate}[(i)]
\item $\phi$ has potential good reduction if and only if
$\spe_{K^\circ}([\phi])\in\bA^2(k)$;
\item $\phi$ has potential additive reduction if and only if
$\spe_{K^\circ}([\phi])=[1:2:0]$;
\item $\phi$ has potential multiplicative reduction if and only if
$\spe_{K^\circ}([\phi])=[1:x:0]$ for some $x\in k\setminus\{2\}$, in which
case $x=a+a^{-1}$, where $a$ is the multiplier of the indifferent reduction
at the crucial point; 
\item $\phi$ has potential constant reduction if and only if
$\spe_{K^\circ}([\phi])=[0:1:0]$.
\end{enumerate}
\end{appcorollary}

The proof in \cite[\S3 and \S4.1]{DJR16} uses fixed-point normal forms,
whereas the above argument derives the stratification from
intrinsic reduction data at $\xi_\phi$.

\section{Intrinsic Type II/III normal forms}
\label{sec:normalformsappendix}

Unless otherwise stated, we retain the notation and assumptions on $K$ and
$k$ from the beginning of Section \ref{sec:intro}.

For completeness, we recall the Type II/III local normal forms in
\cite[\S\S4.8.4--4.8.5]{Jonsson15} and their derivation.
Combining these with the intrinsic
definition \eqref{eq:intrinsicsurplus} of the surplus local degrees
$s_{\Vv}(\phi)$ and conservation of
the ring-theoretically defined local degree, we derive the non-archimedean
argument principle \eqref{eq:FKN}, one of the central formulas in this paper.
The local normal forms are also useful for carrying out classical
non-archimedean analysis in the proof of Theorem \ref{th:residuemoduli}.
Their intrinsic interpretation comes from the
pullback on completed residue fields. For
$\xi\in\sH^1_{\mathrm{II}}\cup\sH^1_{\mathrm{III}}$, denote by
$\mathcal{H}(\xi)$ the completion of $K(z)$ with respect to $|\cdot|_\xi$.
For $\phi\in K(z)$ of degree $>0$ and $\xi':=\phi(\xi)$, pullback induces
\begin{gather*}
 \phi^*:\mathcal{H}(\xi')\hookrightarrow\mathcal{H}(\xi)
 \quad\text{and}\quad
 |\phi^*f|_\xi=|f|_{\xi'}\quad(f\in\mathcal{H}(\xi')).
\end{gather*}
When $\xi'\ne\xi$, the intrinsic reduction $\tilde{\phi}_\xi$ is constant,
whereas this pullback still records the local map from $\xi$ to $\xi'$.
Define the graded reduction
\begin{gather*}
 \gr\mathcal{H}(\xi)
 :=\bigoplus_{t>0}
 \frac{\{f\in\mathcal{H}(\xi):|f|_\xi\le t\}}
 {\{f\in\mathcal{H}(\xi):|f|_\xi<t\}},
\end{gather*}
and denote by $\initial_\xi(f)$ the class of a nonzero $f$ in
degree $|f|_\xi$. The isometric pullback induces
\begin{gather*}
 \gr(\phi^*):
 \gr\mathcal{H}(\xi')\longrightarrow
 \gr\mathcal{H}(\xi)\quad\text{and}\quad
 \initial_{\xi'}(f)\longmapsto
 \initial_\xi(\phi^*f).
\end{gather*}

\subsection{Type II points}

For $\xi$ of type II, the pullback on ordinary residue fields is that of
the tangent map between the natural projective $k$-lines
$T_\xi\sP^1$ and $T_{\xi'}\sP^1$, whose degree is $\deg_\xi(\phi)$ by
\cite[Proposition~4.13 and \S4.8.5]{Jonsson15}.
Define $m_{\Vv}(\phi):=\deg_{\Vv}(\phi_{*,\xi})$ for
$\Vv\in T_\xi\sP^1$. Then for every $\Vw\in T_{\xi'}\sP^1$,
\begin{gather*}
 (\phi_{*,\xi})^*[\Vw]
 =\sum_{\Vv\in T_\xi\sP^1:\,\phi_{*,\xi}(\Vv)=\Vw}
 m_{\Vv}(\phi)[\Vv].
\end{gather*}
The factors in \eqref{eq:typeIInormal} lift the zero and pole divisors of
this tangent map in the chosen coordinates; the factor $1+h$ records that
the quotient by this lift has residue $1$.

\begin{appproposition}[Intrinsic Type II normal form]
\label{th:typeIInormal}
Let $\phi\in K(z)$ be of degree $>0$, let $\xi\in\sH^1_{\mathrm{II}}$, and put
$\xi':=\phi(\xi)$. Fix $\Vv\in T_\xi\sP^1$ and put
$\Vw:=\phi_{*,\xi}(\Vv)$. Choose $M,N\in\PGL(2,K)$ such that
\begin{gather*}
 M(\xi)=N(\xi')=\xi_g\quad\text{and}\quad
 M_*(\Vv)=N_*(\Vw)=\widevec{\xi_g0}.
\end{gather*}
Then there are integers $m\ge1$, $\ell\ge m$, and $k\ge0$,
elements $a,\alpha_i,\beta_j\in K$, and $h\in K(z)$ such that
\begin{gather}
 (N\circ\phi\circ M^{-1})(z)=a z^m
 \frac{\prod_{i=1}^{\ell-m}(z-\alpha_i)}
 {\prod_{j=1}^{k}(z-\beta_j)}(1+h(z))\quad\text{in }K(z)
 \quad\text{and}\quad |h|_{\xi_g}<1,
 \label{eq:typeIInormal}
\end{gather}
where
\begin{gather*}
 m=m_{\Vv}(\phi),\quad
 |a|=|\alpha_i|=|\beta_j|=1,\quad
 \widehat{\alpha_i}\ne\widehat{\beta_j}\text{ for every }i,j,\quad\text{and}\quad
 \max\{\ell,k\}=\deg_\xi(\phi).
\end{gather*}
Moreover,
\begin{gather}
 \widetilde{(N\circ\phi\circ M^{-1})}(\zeta)
 =\widehat a\,\zeta^m
 \frac{\prod_{i=1}^{\ell-m}(\zeta-\widehat{\alpha_i})}
 {\prod_{j=1}^{k}(\zeta-\widehat{\beta_j})}\quad\text{in }k(\zeta).
 \label{eq:typeIIintrinsic}
\end{gather}
In particular, the rational function in \eqref{eq:typeIIintrinsic} and its degree
are respectively the tangent map in these projective coordinates and
$\deg_\xi(\phi)$.
\end{appproposition}

\begin{proof}
By Fact \ref{th:algred} and \cite[\S1.4.2]{Kiwi14}, in these projective
coordinates the tangent map is
$\zeta\mapsto\widetilde{(N\circ\phi\circ M^{-1})}(\zeta)$ and its degree is
$\deg_\xi(\phi)$. Let $q$ be its local degree at $0$. Then
\begin{gather*}
 \widetilde{(N\circ\phi\circ M^{-1})}(\zeta)
 =\widehat a\,\zeta^q
 \frac{\prod_{i=1}^{\ell-q}(\zeta-\widehat{\alpha_i})}
 {\prod_{j=1}^{k}(\zeta-\widehat{\beta_j})}\quad\text{in }k(\zeta)
\end{gather*}
for some $\ell\ge q$ and $k\ge0$, where the factors are counted with
multiplicity and $\max\{\ell,k\}=\deg_\xi(\phi)$. Lift the nonzero zeros,
poles, and the scalar to elements of absolute value $1$ in $K$, and denote the
resulting rational function in $K(z)$ by $G$. Then
$\widetilde G=\widetilde{(N\circ\phi\circ M^{-1})}$ in $k(\zeta)$, so
\begin{gather*}
 (N\circ\phi\circ M^{-1})/G=1+h
 \quad\text{and}\quad |h|_{\xi_g}<1
\end{gather*}
for $h:=(N\circ\phi\circ M^{-1})/G-1\in K(z)$. By the continuity of
$|h|_{\cdot}$ on $\sP^1$, for every $s<1$ sufficiently close to $1$ we have
$|h|_{\xi_{0,s}}<1$, so by \cite[Proposition~4.13 and Lemma~4.14]{Jonsson15},
\begin{gather*}
 (N\circ\phi\circ M^{-1})(\xi_{0,s})=\xi_{0,s^q}
 \quad\text{and}\quad
 \deg_{\xi_{0,s}}(N\circ\phi\circ M^{-1})=q.
\end{gather*}
Thus the slope at $\xi_g$ in $\widevec{\xi_g0}$ is
$q=m_{\Vv}(\phi)$ by the definition above. Setting $m:=q$ proves
\eqref{eq:typeIInormal} and \eqref{eq:typeIIintrinsic}.
\end{proof}

\subsection{Type III points}

For $\xi$ of type III, the ordinary residue field of $\mathcal{H}(\xi)$ is
$k$. Define $m_{\Vv}(\phi):=\deg_\xi(\phi)$ for each of the two directions
$\Vv\in T_\xi\sP^1$. The graded pullback retains both the exponent and the leading
coefficient of the local map. In the notation of Proposition
\ref{th:typeIIInormal} below, its coordinate expression is
\begin{gather*}
 \gr(\phi^*)\bigl(\initial_{\xi'}(N)\bigr)
 =\initial_\xi(a)\initial_\xi(M)^m.
\end{gather*}
The exponent $m=\deg_\xi(\phi)$ is the local metric expansion factor:
on a sufficiently short segment containing $\xi$ in its interior, we have
\begin{gather*}
 \rho\bigl(\phi(\eta),\phi(\xi)\bigr)
 =\deg_\xi(\phi)\,\rho(\eta,\xi).
\end{gather*}
The factor $1+h$ in \eqref{eq:typeIIInormal} has initial form $1$.
In the fixed-point form below, the graded pullback also retains the
coefficient residue $\widehat a$ used in the proof of Proposition
\ref{th:multiple}.

\begin{appproposition}[Intrinsic Type III normal form]
\label{th:typeIIInormal}
Let $\phi\in K(z)$ be of degree $>0$, let $\xi\in\sH^1_{\mathrm{III}}$, and put
$\xi':=\phi(\xi)$. Fix $\Vv\in T_\xi\sP^1$ and put
$\Vw:=\phi_{*,\xi}(\Vv)$. Choose $M,N\in\PGL(2,K)$ such that
\begin{gather*}
 M(\xi)=\xi_{0,r},\quad N(\xi')=\xi_{0,r'},\quad
 M_*(\Vv)=\widevec{\xi_{0,r}0},\quad\text{and}\quad
 N_*(\Vw)=\widevec{\xi_{0,r'}0}
\end{gather*}
for some $r,r'\notin|K^\times|$. Then
\begin{gather}
 (N\circ\phi\circ M^{-1})(z)=a z^m(1+h(z))\quad\text{in }K(z),
 \quad |h|_{\xi_{0,r}}<1,
 \quad\text{and}\quad |a|r^m=r',
 \label{eq:typeIIInormal}
\end{gather}
for some $a\in K^\times$ and $h\in K(z)$, where
\begin{gather*}
 m=m_{\Vv}(\phi)=\deg_\xi(\phi).
\end{gather*}
\end{appproposition}

\begin{proof}
On a sufficiently small annulus about $\xi_{0,r}$, write the Laurent expansion of
$N\circ\phi\circ M^{-1}$. There is a unique dominant Laurent monomial. Indeed, if two distinct
terms $c_i z^i$ and $c_j z^j$ had the same norm at $\xi_{0,r}$, then
\begin{gather*}
 r^{i-j}=|c_j/c_i|\in|K^\times|,
\end{gather*}
which would imply $r\in|K^\times|$ since $K$ is algebraically closed. Thus
\begin{gather*}
 (N\circ\phi\circ M^{-1})(z)=a z^m(1+h(z))\quad\text{in }K(z)
 \quad\text{and}\quad |h|_{\xi_{0,r}}<1
\end{gather*}
for some $m\in\bZ$. Since
$M_*(\Vv)=\widevec{\xi_{0,r}0}$ and
$N_*(\Vw)=\widevec{\xi_{0,r'}0}$, we have $m>0$. On either side of $\xi_{0,r}$, the radius is
therefore transformed by $s\mapsto |a|s^m$, and
$(N\circ\phi\circ M^{-1})(\xi_{0,r})=\xi_{0,r'}$ implies
$|a|r^m=r'$. The ordinary residue fields of $\mathcal{H}(\xi)$ and
$\mathcal{H}(\xi')$ are both $k$, and their value groups in these coordinates
are $|K^\times|r^{\bZ}$ and $|K^\times|r^{m\bZ}$, respectively.
Since $|K^\times|$ is divisible and $r\notin|K^\times|$,
\cite[Proposition~4.13, (17)]{Jonsson15} yields
\begin{gather*}
 \deg_\xi(\phi)
 =[\mathcal{H}(\xi):\phi^*\mathcal{H}(\xi')]
 =[|K^\times|r^{\bZ}:|K^\times|r^{m\bZ}]
 =m=m_{\Vv}(\phi).
\end{gather*}
Moreover, by \cite[Lemma~4.14]{Jonsson15}, for every $s$ sufficiently close to
$r$, we have
$\deg_{\xi_{0,s}}(N\circ\phi\circ M^{-1})=m$.
\end{proof}

\begin{appcorollary}[Fixed-point form]
\label{th:fixednormal}
Let $\phi\in K(z)$ be of degree $>0$, let
$\xi\in\sH^1_{\mathrm{II}}\cup\sH^1_{\mathrm{III}}$ be fixed by $\phi$, and
let $\Vv_0,\Vv_\infty\in T_\xi\sP^1$ be two distinct directions.
Suppose that both are fixed by $\phi_{*,\xi}$.
If $\xi\in\sH^1_{\mathrm{II}}$, choose $M\in\PGL(2,K)$ so that
\begin{gather*}
 M(\xi)=\xi_g,\quad
 M_*(\Vv_0)=\widevec{\xi_g0},\quad\text{and}\quad
 M_*(\Vv_\infty)=\widevec{\xi_g\infty}.
\end{gather*}
Then there are integers $m\ge1$, $\ell\ge m$, and $k\ge0$,
elements $a,\alpha_i,\beta_j\in K$, and $h\in K(z)$ such that
\begin{gather*}
 (M\circ\phi\circ M^{-1})(z)
 =a z^m
 \frac{\prod_{i=1}^{\ell-m}(z-\alpha_i)}
 {\prod_{j=1}^{k}(z-\beta_j)}(1+h(z))\quad\text{in }K(z)
 \quad\text{and}\quad |h|_{\xi_g}<1,
\end{gather*}
where
\begin{gather*}
 m=m_{\Vv_0}(\phi),\quad
 |a|=|\alpha_i|=|\beta_j|=1,\quad
 \widehat{\alpha_i}\ne\widehat{\beta_j}\text{ for every }i,j,\quad
 \ell=\deg_\xi(\phi),\quad\text{and}\quad
 \ell-k=m_{\Vv_\infty}(\phi).
\end{gather*}
and
\begin{gather*}
 \widetilde{(M\circ\phi\circ M^{-1})}(\zeta)
 =\widehat a\,\zeta^m
 \frac{\prod_{i=1}^{\ell-m}(\zeta-\widehat{\alpha_i})}
 {\prod_{j=1}^{k}(\zeta-\widehat{\beta_j})}
 \quad\text{in }k(\zeta).
\end{gather*}

If $\xi\in\sH^1_{\mathrm{III}}$, choose $M\in\PGL(2,K)$ and
$r\notin|K^\times|$ so that
\begin{gather*}
 M(\xi)=\xi_{0,r},\quad
 M_*(\Vv_0)=\widevec{\xi_{0,r}0},\quad\text{and}\quad
 M_*(\Vv_\infty)=\widevec{\xi_{0,r}\infty}.
\end{gather*}
Then necessarily
$\deg_\xi(\phi)=m_{\Vv_0}(\phi)=m_{\Vv_\infty}(\phi)=1$,
and for some $a\in K^\times$ and some $h\in K(z)$, 
\begin{gather*}
 (M\circ\phi\circ M^{-1})(z)=a z(1+h(z))\quad\text{in }K(z),
 \quad |a|=1,
 \quad\text{and}\quad |h|_{\xi_{0,r}}<1.
\end{gather*}
\end{appcorollary}

\begin{proof}
Suppose first that $\xi\in\sH^1_{\mathrm{II}}$. Apply Proposition
\ref{th:typeIInormal} with $N=M$ and $\Vv=\Vv_0$. Since
$M_*(\Vv_\infty)$ is fixed, the numerator degree $\ell$ in
\eqref{eq:typeIInormal} is strictly greater than the denominator degree $k$.
Hence $\ell=\deg_\xi(\phi)$, while the local degree at $\infty$ is
$\ell-k=m_{\Vv_\infty}(\phi)$, which proves the asserted form.

Suppose next that $\xi\in\sH^1_{\mathrm{III}}$. Apply Proposition
\ref{th:typeIIInormal} with $N=M$. Then
\begin{gather*}
 (M\circ\phi\circ M^{-1})(z)=a z^{\deg_\xi(\phi)}(1+h(z))\quad\text{in }K(z),
 \quad |h|_{\xi_{0,r}}<1,
 \quad\text{and}\quad |a|r^{\deg_\xi(\phi)}=r.
\end{gather*}
If $\deg_\xi(\phi)>1$, then
$r^{\deg_\xi(\phi)-1}=|a|^{-1}\in|K^\times|$, hence
$r\in|K^\times|$ since $K$ is algebraically closed, a contradiction. Therefore
$\deg_\xi(\phi)=1$, and then $|a|=1$; moreover
$m_{\Vv_0}(\phi)=m_{\Vv_\infty}(\phi)=1$ by Proposition \ref{th:typeIIInormal}.
\end{proof}

\subsection{Derivation of the non-archimedean argument principle}

We derive \eqref{eq:FKN} from the local forms above and conservation of
degree. Let $\phi\in K(z)$ have degree $d\ge1$, let $\xi\in\sH^1$, and set
$\xi':=\phi(\xi)$. For $\Vv\in T_\xi\sP^1$, define
\begin{gather}
 s_{\Vv}(\phi)
 :=(\phi^*\delta_{\xi'})(U(\Vv))
 =\sum_{x\in U(\Vv):\,\phi(x)=\xi'}\deg_x(\phi).
 \label{eq:intrinsicsurplus}
\end{gather}

For $\xi$ of type II, this definition also has a ring-theoretic interpretation.
Choose $M$ and $N$ in $\PGL(2,K)$ such that $M(\xi)=N(\xi')=\xi_g$,
and denote by $[\bar{F}:\bar{G}]\in\bP^{2d+1}(k)$ the coefficient reduction of
$N\circ\phi\circ M^{-1}$. Let $a\in\bP^1(k)$ correspond to $M_*\Vv$ and set
$R_a:=\cO_{\bP^1_k,a}$. Using a local trivialization of $\cO(d)$
at $a$, regard $\bar{F}$ and $\bar{G}$ as elements of $R_a$. Since
$N\circ\phi\circ M^{-1}$ fixes $\xi_g$, Remark \ref{th:algdepth} and the fact
that $R_a$ is a discrete valuation ring imply
\begin{gather*}
 s_{\Vv}(\phi)
 =\operatorname{length}_{R_a}\bigl(R_a/(\bar{F},\bar{G})\bigr)
 =\ord_a[\gcd(\bar{F},\bar{G})=0].
\end{gather*}
Changes of source and target coordinates and of the local trivialization
preserve this colength of the base ideal. The common factor disappears from
the induced map on residue fields, so that map alone does not determine this
colength.

For type II and III points, use the definitions of $m_{\Vv}(\phi)$ above;
for a type IV point, define $m_{\Vv}(\phi):=\deg_\xi(\phi)$ in its unique
direction. Set $n_{\Vv}(y):=(\phi^*\delta_y)(U(\Vv))$ for $y\in\sP^1$.
For a connected simple neighborhood $V$ of $y\ne\xi'$ disjoint from
$\{\xi'\}$, every component of $\phi^{-1}(V)$ lies entirely in $U(\Vv)$
or entirely outside it, since $\partial U(\Vv)=\{\xi\}$.
By \cite[Proposition~4.5]{Jonsson15}, its degree over $V$ is constant.
Thus $n_{\Vv}$ is locally constant on $\sP^1\setminus\{\xi'\}$ and hence
constant on each $U(\Vw)$ with $\Vw\in T_{\xi'}\sP^1$.

Suppose first that $\xi$ is of type II or III. Since $\phi$ is finite,
we can choose a sufficiently small connected simple neighborhood $V$ of
$\xi'$ so that the components $W_x$ of $\phi^{-1}(V)$ are indexed by
$x\in\phi^{-1}(\xi')$, with $W_x\cap\phi^{-1}(\xi')=\{x\}$ and
$W_x\subset U(\Vv)$ or $W_x\cap U(\Vv)=\emptyset$ for every $x\ne\xi$.
Then \cite[Proposition~4.5]{Jonsson15} implies
\begin{gather*}
 \sum_{z\in W_x:\,\phi(z)=y}\deg_z(\phi)=\deg_x(\phi)
 \quad(y\in V).
\end{gather*}
Fix $\Vw\in T_{\xi'}\sP^1$ and choose $y\in V\cap U(\Vw)$ on a sufficiently
short segment issuing from $\xi'$. For every $\Vu\in T_\xi\sP^1$ with
$\phi_{*,\xi}(\Vu)=\Vw$, the local forms in Propositions
\ref{th:typeIInormal} and \ref{th:typeIIInormal}, including the local degrees
along the corresponding segments computed in their proofs, produce
\begin{gather*}
 x_{\Vu}(y)\in W_\xi\cap U(\Vu),\quad
 \phi(x_{\Vu}(y))=y,\quad\text{and}\quad
 \deg_{x_{\Vu}(y)}(\phi)=m_{\Vu}(\phi).
\end{gather*}
There are finitely many such directions, so $y$ can be chosen to satisfy
these conditions for all of them. The degree of the tangent map in type II
and the bijection between the two directions in type III, imply
\begin{gather*}
 \sum_{z\in W_\xi:\,\phi(z)=y}\deg_z(\phi)
 =\deg_\xi(\phi)
 =\sum_{\Vu\in T_\xi\sP^1:\,\phi_{*,\xi}(\Vu)=\Vw}m_{\Vu}(\phi),
\end{gather*}
so the displayed points exhaust the preimages of $y$ in $W_\xi$.
Adding the contributions from the other $W_x$, we obtain
\begin{gather*}
 n_{\Vv}(y)
 =\sum_{x\in U(\Vv):\,\phi(x)=\xi'}\deg_x(\phi)
   +m_{\Vv}(\phi)1_{\{\phi_{*,\xi}(\Vv)=\Vw\}}
 =s_{\Vv}(\phi)+m_{\Vv}(\phi)1_{U(\phi_{*,\xi}(\Vv))}(y).
\end{gather*}
Since $\Vw$ was arbitrary, constancy on $U(\Vw)$ proves this equality for
every $y\ne\xi'$, and the case $y=\xi'$ is the definition of
$s_{\Vv}(\phi)$.

If $\xi$ is of type IV, then so is $\xi'$ by \cite[Lemma~4.2]{Jonsson15}.
Both points have only one direction, so
\begin{gather*}
 n_{\Vv}(y)
 =d-\deg_\xi(\phi)1_{\{\xi'\}}(y)
 =(d-\deg_\xi(\phi))
   +\deg_\xi(\phi)1_{\sP^1\setminus\{\xi'\}}(y)
 =s_{\Vv}(\phi)+m_{\Vv}(\phi)1_{U(\phi_{*,\xi}(\Vv))}(y).
\end{gather*}
This proves \eqref{eq:FKN} for every $\xi\in\sH^1$.
Finally, the partition of $\sP^1\setminus\{\xi\}$ into its directions yields
\begin{gather*}
 \sum_{\Vv\in T_\xi\sP^1}s_{\Vv}(\phi)
 =(\phi^*\delta_{\xi'})(\sP^1\setminus\{\xi\})
 =d-\deg_\xi(\phi),
\end{gather*}
which is \eqref{eq:surplus}.

\begin{acknowledgement}
 This research was partially supported by JSPS Grant-in-Aid for Scientific Research (C), 23K03129, and
by the Research Institute for Mathematical Sciences,
an International Joint Usage/Research Center at Kyoto University.
\end{acknowledgement}

\end{document}